\documentclass{amsart}    
\usepackage{amsmath} 
\usepackage{paralist}
\usepackage{graphics} 
\usepackage{epsfig} 
\usepackage{graphicx}  \usepackage{epstopdf} 
\usepackage[colorlinks=true]{hyperref}
\hypersetup{urlcolor=blue, citecolor=red}

\newtheorem{Definition}{Definition}[section] 

\newtheorem{Proposition}{Proposition}[section]

\newtheorem{Lemme}{Lemma}[section]
\newtheorem{Theoreme}{Theorem}[section]

\newtheorem{Remarque}{Remark}

\def \vf{\vec{f}} 
\def \vg{\vec{g}}
\def \vu{\vec{u}}
\def \vv{\vec{v}}
\def \vw{\vec{w}}
\def \vvarphi{\vec{\varphi}}
\def \P{\mathbb{P}}
\def \U{\vec{U}}

\def \R{\mathbb{R}}
\def \Rt{\mathbb{R}^3}
\def \T{\mathbb{T}^3}
\def \Z{\mathbb{{\bf Z}}^3 \setminus \{0\}} 
\def \finpv{\hfill $\blacksquare$}
\def \pv{{\bf{Proof.}}~} 

\def \ds{\displaystyle}

\usepackage{amssymb}

\title[\bf New Navier-Stokes-alpha model with nonlinear filter] %
{ Mathematical study of a new Navier-Stokes-alpha model with  nonlinear filter equation - Part I} 
\author[ Manuel Fernando Cortez  and Oscar Jarr\'in]{}


\email{oscar.jarrin@udla.edu.ec}
\email{manuel.cortez@epn.edu.ec} 

\thanks{$^*$Corresponding author:  Oscar Jarr\'in}

\begin{document}
	\maketitle
	
 	\centerline{\scshape Manuel Fernando Cortez}
	\medskip
	{\footnotesize
		\centerline{Departamento de Matem\'aticas}
		\centerline{Escuela Politécnica Nacional} 
		\centerline{Ladr\'on de Guevera E11-253, Quito, Ecuador} 
	}

 \medskip
	
	\centerline{\scshape Oscar Jarr\'in$^*$}
	\medskip
	{\footnotesize
		\centerline{Escuela de Ciencias Físicas y Matemáticas}
		\centerline{Universidad de Las Américas}
		\centerline{V\'ia a Nay\'on, C.P.170124, Quito, Ecuador}
	} 
	
\bigskip
\begin{abstract} This article is devoted to the mathematical study of a new Navier-Stokes-alpha model with a nonlinear filter equation. For a given indicator function, this filter equation was first considered by W. Layton, G. Rebholz, and C. Trenchea in [Modular Nonlinear Filter Stabilization of Methods for Higher Reynolds Numbers Flow, J. Math. Fluid Mech. 14: 325–354 (2012)] to select eddies for damping based on the understanding of how nonlinearity acts in real flow problems. Numerically, this nonlinear filter equation was applied to the nonlinear term in the Navier-Stokes equations to provide a precise analysis of numerical diffusion and error estimates.

Mathematically, the resulting alpha-model is described by a doubly nonlinear parabolic-elliptic coupled system. We therefore undertake the first theoretical study of this system by considering periodic boundary conditions in the spatial variable. Specifically, we address the existence and uniqueness of weak Leray-type solutions, their rigorous convergence to weak Leray solutions of the classical Navier-Stokes equations, and their long-time dynamics through the concept of the global attractor and some upper bounds for its fractal dimension.

Handling the nonlinear filter equation together with the well-known nonlinear transport term makes certain estimates delicate, particularly when deriving upper bounds on the fractal dimension. For the latter, we adapt techniques developed for hyperbolic-type equations. 
\end{abstract}
{\footnotesize \tableofcontents}
\section{Introduction and motivation of the model}
The study of fluid dynamics provides us with several nonlinear fluid-type models of relevance. Among them are the well-known Navier-Stokes equations and some related equations known as the alpha-models. Alpha-models have been developed in the mathematical and experimental literature as physically relevant approximations of the Navier-Stokes equations. Examples include Bardina's model \cite{Bardina1,Bardina}, the viscous Camassa-Holm model \cite{Chen}, the Clark-alpha model \cite{Cao}, and the Leray-alpha model \cite{CheskidovHolm}. 

\medskip

Numerical solutions of the Navier–Stokes equations for problems of physical and engineering relevance are currently not feasible, as the mathematical theory concerning the uniqueness and regularity of Leray's solutions remains one of the most challenging open questions \cite{Dubois,Pope}. Therefore, alpha-models serve as \emph{regularized} approximations of the classical Navier-Stokes equations, for which the global well-posedness of Leray's solutions can be established. In addition, this regularization approach demonstrates appealing robustness at high Reynolds numbers, suggesting that alpha-models are suitable for large eddy simulations of turbulence \cite{Berselli,Geurts,John,Lesieur,Sagaut}.  

\medskip

Despite some structural differences in the equations of the  alpha-models, their primary common characteristic is the inclusion of a \emph{filtered velocity} in their nonlinear transport term. For a parameter $\alpha > 0$, the filtered velocity $\vu_\alpha$ is obtained from the fluid velocity $\vu$ by solving the Helmholtz filter equation:
\begin{equation}\label{Helmholtz-Filter}
    -\alpha^2 \Delta \vu_\alpha + \vu_\alpha = \vu, \quad \text{div}(\vu_\alpha)=0.  
\end{equation}
In the alpha-models, or any turbulence model, the filtered velocity $\vu_\alpha$ enables us to describe the large-scale motion of the fluid, effectively filtering or averaging fluid motions at scales smaller than the specified length-scale $\alpha$. 

\medskip

Nevertheless, recent experimental studies \cite{Geurts0, Geurts1, Layton} provide evidence that this \emph{linear} filter equation numerically truncates scales uniformly in space. Consequently, it over-regularizes the laminar parts of the solution and often removes critical flow structures. To overcome this problem, a \emph{nonlinear} filter equation was proposed in \cite[Definition $2.2$]{Layton2}:
\begin{equation*}
-\alpha^2 \text{div} \left( A (\vu) (\vec{\nabla} \otimes \vu_\alpha)^T \right) + \vu_\alpha + \vec{\nabla}P_\alpha= \vu, \quad
\text{div}(\vu_\alpha) = 0,
\end{equation*}
where $(\vec{\nabla} \otimes \vu_\alpha)^T$ denotes the transpose matrix of $\vec{\nabla} \otimes \vu_\alpha$; and where the given \emph{indicator function} $A:\mathbb{R}^3 \to \mathbb{R}$ selects regions requiring local filtering. 
Roughly speaking, in simulations of high Reynolds number flows, nonlinear terms break down eddies into smaller ones until they are small enough for molecular viscosity to dissipate them rapidly. This transfer of energy occurs due to nonlinearity in the actual flow. In this context, the nonlinear filter adjusts the amount and location of eddy viscosity based on the local flow structures. For a more detailed explanation, please refer to \cite{Layton2}.  On the other hand, $P_\alpha$ is a pressure term formally related to the velocity $\vu$ and the filtered velocity $\vu_\alpha$ by the equation
\[ - \Delta P_\alpha = - \alpha^2 \text{div}\Big( \text{div} \left( A (\vu) (\vec{\nabla} \otimes \vu_\alpha)^T \right)\Big),  \]
due to the divergence-free property of $\vu$ and $\vu_\alpha$. When comparing with the linear Helmholtz filter equation (\ref{Helmholtz-Filter}), we observe that any pressure term $P_\alpha$ in the equation $-\alpha^2 \Delta \vu_\alpha + \vu_\alpha + \vec{\nabla}P_\alpha=\vu$, is formally
related to the divergence-free vector fields $\vu$ and $\vu_\alpha$ as $-\Delta P_\alpha= \text{div}\Big( -\alpha^2 \Delta \vu_\alpha + \vu_\alpha - \vu \Big)=0$. Consequently,  mathematical frameworks such as the whole space $\Rt$ or periodic boundary conditions lead to the conclusion that $\vec{\nabla}P_\alpha = 0$. In contrast, on bounded domains, this is not true in general.  

\medskip

This nonlinear filter was applied to the nonlinear transport term in the classical Navier-Stokes equations in \cite{Bowers1, Bowers2, Layton2}, and more recently in \cite{Takhirov-Trenchea}, to conduct numerical studies on stabilizing algorithms and error analysis. Mathematically, for the fluid velocity $\vu=\vu(t,x)\in \Rt$ and the pressure $P(t,x)\in \R$, one gets the following new Navier-Stokes-alpha model:
\begin{equation}\label{Alpha-model-Intro}
\begin{cases}\vspace{2mm}
\partial_t \vu - \nu \Delta \vu + (\vu_\alpha \cdot \vec{\nabla})\vu + \vec{\nabla} P = \vf, \qquad \ \nu>0,  \ \  \text{div}(\vu)=0,\\
  -\alpha^2 \text{div} \left( A (\vu) (\vec{\nabla} \otimes \vu_\alpha)^T \right)+\vu_\alpha +\vec{\nabla}P_\alpha = \vu, \qquad \alpha>0, \ \  \text{div}(\vu_\alpha)=0.
\end{cases}
\end{equation}
We are dealing with a doubly nonlinear parabolic-elliptic coupled system, where $\nu$ denotes the viscosity constant, $\vf=\vf(t,x)\in \Rt$ is a given external force acting on the fluid, and $A(\cdot)\in \R$ is a specified indicator function.

\medskip

To the best of our knowledge, the coupled system (\ref{Alpha-model-Intro}) has not been studied theoretically. The main objective of this article is to mathematically validate this model. Therefore, we address the following program for the qualitative study of this system:
\begin{enumerate}
\item[(i)] The existence of global-in-time weak solutions in the natural energy space.
\item[(ii)] The uniqueness of these solutions.
\item[(iii)] The asymptotic behavior of solutions as the parameter $\alpha$ approaches 0 and the indicator function $A(\cdot)$ approaches 1 in the nonlinear filter equation.
\item[(iv)] The long-term behavior of solutions using the concept of the global attractor and upper bounds on its fractal dimension.
\end{enumerate}
The existence of finite-energy weak solutions is a somewhat standard result, owing to Leray's well-known method developed for the classical Navier-Stokes equations. In contrast, uniqueness of these solutions is a fundamental property expected for any alpha-model. Experimentally, the uniqueness of solutions holds great importance in numerical algorithms. Moreover, the uniqueness of the solution is a key element for successfully applying tools of dynamical systems to study their long-term behavior.  

\medskip

Uniqueness of solutions to the coupled system (\ref{Alpha-model-Intro}) strongly depends on the appropriate continuity and regularity properties of the indicator function $A(\cdot)$. However, this information alone may not be sufficient to rigorously prove this fact. In the following lines, we provide a brief explanation. 

\medskip

Returning to the first equation in the system (\ref{Alpha-model-Intro}), both the velocity $\vu$ and the filtered velocity $\vu_\alpha$ satisfy a divergence-free property. Standard \emph{a priori} energy estimates show that $\vu$ shares the same energy controls as Leray's solutions in the classical Navier-Stokes equations, specifically $\vu\in L^\infty_t L^2_x \cap L^2_t \dot{H}^{1}_x$. According to the well-known theory of the Navier-Stokes equations, this alone is insufficient to prove uniqueness\footnote{Interestingly, there is emerging research suggesting non-uniqueness. See, for instance, \cite{Albritton} and its recent citations.} of solutions. Therefore, we need additional information about the filtered velocity $\vu_\alpha$. Specifically, we require $H^2$-estimates on $\vu_\alpha$, which can be obtained from the elliptic nonlinear filter equation in (\ref{Alpha-model-Intro}) as long as the term $A(\vu)$ and its derivatives $\vec{\nabla}(A(\vu))$ are continuous and bounded. However, the classical chain rule implies that this necessary information about $\vec{\nabla}(A(\vu))$ enables us to control $\vec{\nabla}\otimes \vu$ in functional spaces of continuous and bounded functions. Achieving this seems difficult with only the information $\vu \in L^\infty_t L^2_x \cap L^2_t \dot{H}^{1}_x$, thus resulting in $\vec{\nabla}\otimes \vu \in L^2_t L^2_x$.

\medskip

It is worth mentioning that suitable continuity and bounded properties of the terms $A(\vu)$ and $\vec{\nabla}(A(\vu))$ are essential tools for proving the regularity of solutions to elliptic equations, such as the filter equation in (\ref{Alpha-model-Intro}). For further insights, refer to \cite[Chapter 13]{Gilbarg},\cite{LaManna,Zhang} and the references therein. Consequently, obtaining the required $H^2$-estimates for the filtered velocity $\vu_\alpha$ appears challenging within the framework where velocities $\vu$ belong to $L^\infty_t L^2_x \cap L^2_t \dot{H}^{1}_x$.

\medskip

As highlighted in  \cite[Page $328$]{Layton2}, certain regularizations of the filter equation may be necessary for practical computations. Thus, to ensure the uniqueness of weak solutions in the natural energy space $L^\infty_t L^2_x \cap L^2_t \dot{H}^{1}_x$, we consider the following \emph{regularized} nonlinear filter:
\[ -\alpha^2 \text{div} \left( A (\varphi \ast \vu) (\vec{\nabla} \otimes \vu_\alpha)^T \right)+\vu_\alpha + \vec{\nabla}P_\alpha= \vu,\]
where $\varphi=\varphi(x)\in \mathbb{R}$ is a given convolution term with certain assumed useful properties. We are addressing the initial value problem of the following Navier-Stokes-alpha model:
\begin{equation}\label{Alpha-model}
    \begin{cases}\vspace{2mm}
    \partial_t \vu - \nu\Delta \vu + (\vu_\alpha \cdot \vec{\nabla})\vu + \vec{\nabla} P  =  \vf, \quad \text{div}(\vu)=0, \\ \vspace{2mm}
     -\alpha^2 \text{div} \left( A (\varphi\ast\vu) (\vec{\nabla} \otimes \vu_\alpha)^T \right)+\vu_\alpha + \vec{\nabla}P_\alpha = \vu, \quad \text{div}(\vu_\alpha)=0,  \\ 
\vu(0,\cdot)=\vu_0, \quad \text{div}(\vu_0)=0,
    \end{cases}
\end{equation}
where $\vu_0=\vu_0(x)$ denotes the initial data. The regularized term $A(\varphi\ast\vu)$ in the nonlinear filter equation will allow us to control $\vec{\nabla}(A(\varphi\ast\vu))$ suitably. Please refer to Remark \ref{Rmk1} and Remark \ref{Rmk2} for more details. Consequently,   for this alpha-model  we can  rigorous prove all the points $(i)-(iv)$ in the program stated above.  

\medskip

Although the convolution product with the function $\varphi$ is a technical requirement for our computations to work, note that, in point (iii) above, we will be able to show convergence to the classical Navier-Stokes equations as $\alpha \to 0$. In this sense, the mollification applied to the nonlinear filter equation still makes this model a satisfactory \emph{theoretical approximation} of the Navier-Stokes equations, which is one of the main features of alpha-models (see \cite[Chapter $17$]{PLe}). Additionally, in Section \ref{Sec:Example} below,  we provide a simple example of a convolution $\varphi$ formally satisfying
\[\varphi \ast \vu \approx \vu, \]
This example suggests  that one could construct (more sophisticated) explicit convolution terms $\varphi$ such that
\[ A(\varphi \ast \vu ) \approx A(\vu). \]
Consequently, the local filtering effects of the indicator function $A(\cdot)$ with respect to the velocity $\vu$ are not strongly affected by $\varphi$.

\medskip

On the other hand, it is natural to ask for the asymptotic behavior of solutions to the system (\ref{Alpha-model}) when we set $\varphi=\varphi_\varepsilon$ in the mollified filter equation, where $(\varphi_\varepsilon)_{\varepsilon>0}$ is a suitable family of approximations of the identity. This study is addressed at the end of the article in Appendix \ref{Appendix}. Under stronger continuity assumptions on the indicator function $A(\cdot)$, we are able to obtain a finite-energy weak solution $\vu$ to the system (\ref{Alpha-model-Intro}) as the limit of solutions $\vu_\varepsilon$ to the system (\ref{Alpha-model}). Nevertheless, as explained, uniqueness of this solution $\vu$ is an open question without smoothing effects of the convolution term $\varphi$. Consequently, the remaining points $(ii)-(iv)$ in our research program for the model (\ref{Alpha-model-Intro}) are also open questions.

\section{Main results}
In a first theoretical study of the coupled system (\ref{Alpha-model}), we shall consider these equations with \emph{periodic boundary conditions} on the torus $\T=[0, L]^3$, where $L > 0$ is the period. Periodic conditions provide us with a useful framework to study fluid-type equations, leveraging Fourier series as a key tool. Consequently, they offer a rich set of tools for developing rigorous mathematical analysis and are also suitable for computer simulations. Following the standard periodic setting, we denote by $\mathcal{P}$ the set of vector trigonometric polynomials with periodic boundary conditions on $\T$. Then, we define $\mathcal{V} = \left\{ \vvarphi \in \mathcal{P} : \text{div}(\vvarphi) = 0 \text{ and } \int_{\T} \vvarphi(x)dx = 0 \right\}$. Thus, we define the \emph{Hilbert space} of \emph{solenoidal and zero-mean vector fields} $L^2_{\text{s}}(\T)$ as the closure of $\mathcal{V}$ in $L^2(\T)$. 

\medskip

For clarity, we have divided our research program $(i)-(iv)$ into three parts. In the section named \emph{Global well-posedness}, we address points $(i)$ and $(ii)$. Then, in the section \emph{Validity of the model}, we study point $(iii)$, and finally, in the section \emph{Long-time dynamics}, we examine point $(iv)$. 

\medskip

As mentioned, all these qualitative properties of the coupled system (\ref{Alpha-model}) strongly depend on suitable conditions for the indicator function $A(\cdot)$ and the convolution term $\varphi$ in the nonlinear filter equation. In the statement of our results, we explicitly state them separately. Additionally, for a given parameter $0 < \beta < 1$, we assume that the indicator $A(\cdot):\mathbb{R}^3 \to \mathbb{R}$ satisfies:
\begin{equation}\label{Conditions-a-existence-1}
\beta \leq A(x) \leq 1, \quad \text{for all } x \in \mathbb{R}^3.
\end{equation}
The upper bound $A(\cdot) \leq 1$ is defined in \cite[Definition 2.1]{Layton2} to ensure that the nonlinear filter equation in (\ref{Alpha-model}) locally behaves like the linear Helmholtz filter equation (\ref{Helmholtz-Filter}) in selected regions where $A(\cdot) = 1$. On the other hand, the lower bound $\beta < A(\cdot)$ was considered in \cite{Takhirov-Trenchea} to enforce the positivity of the indicator function. Experimentally, this property is useful for the numerical stabilization of the alpha-model. Mathematically, we will demonstrate that this property provides suitable controls on the filtered velocity $\vu_\alpha$.

\medskip

{\bf Global well-posedness}. We start by introducing the following:
\begin{Definition}\label{Def-Leray-Solutions} Let $\alpha>0$ fixed. We shall say that $(\vu, \vu_\alpha, P, P_\alpha)$ is a weak Leray solution to the coupled parabolic-elliptic system (\ref{Alpha-model}) if the following conditions hold:
\begin{enumerate}
	\item[(A)] We have 
 \[\vu \in L^{\infty}_{loc}\big([0,+\infty[, L^2_{\text{s}}(\T)\big) \cap L^{2}_{loc}\big( [0,+\infty[, \dot{H}^1(\T)\big), \quad \vu_\alpha \in L^{\infty}_{loc}\big([0,+\infty[, H^1(\T)\big),\]
 and 
 \[ P\in L^2_{loc}\big( [0,+\infty[, L^2(\T)\big), \quad P_\alpha \in L^\infty_{loc}\big( [0,+\infty[, L^2(\T)\big). \]
 \item[(B)] The set $(\vu,\vu_\alpha,P, P_\alpha)$  solves the  system (\ref{Alpha-model}) in the weak sense. 
	\item[(C)] For all $t \geq 0$, the velocity $\vu$ verifies the energy control: 
	\begin{equation}\label{Energy-Control}
	\| \vu(t,\cdot)\|^{2}_{L^2} + 2 \nu \int_{0}^{t} \| \vu(s,\cdot)\|^{2}_{\dot{H}^1}ds \leq \| \vu_0 \|^{2}_{L^{2}} + 2 \int_{0}^{t} \big\langle f(s,\cdot), \vu(s,\cdot)\big\rangle_{\dot{H}^{-1}\times\dot{H}^1}\, ds,
	\end{equation}
     provided that $\vu_0 \in L^2_{\text{s}}(\T)$,  and  $\vf \in L^{2}_{loc}\big([0,+\infty[, \dot{H}^{-1}(\T)\big)$ with $\text{div}(\vf)=0$. 
\end{enumerate}		
\end{Definition}	
\begin{Remarque} For the sake of simplicity, we assume an external force that is divergence-free. However, with minor technical modifications, all our definitions and results hold in the more general setting involving any external force. 
\end{Remarque}

As mentioned in the introduction, \emph{a priori} energy estimates on the first equation in the system (\ref{Alpha-model}) show that the velocity $\vu$ formally satisfies the energy control (\ref{Energy-Control}). This fact enables us to impose conditions $(A)$, $(B)$, and $(C)$ above to naturally define weak Leray solutions. 

\medskip

For $m\in \mathbb{N}$, recall that $\mathcal{C}^m(\Rt)$ denotes the space of functions that are $m$-times continuously differentiable with bounded derivatives. Moreover, for $0<s\leq 1$, we define the H\"older space $\mathcal{C}^{m,s}(\Rt)$  as the space of $\mathcal{C}^m$-functions whose derivatives of order $m$ are $s$-H\"older continuous. In the particular case when $s=1$, the derivatives of order $m$ are Lipschitz continuous.

\begin{Theoreme}[Existence of weak Leray solutions]\label{Th-Existence} Let $\vu_0 \in L^2_{\text{s}}(\T)$ be an initial datum, and let $\vf \in L^{2}_{loc}([0,+\infty[, \dot{H}^{-1}(\T))$ be a divergence-free and zero-mean  external force. In the nonlinear filter equation of the system (\ref{Alpha-model}) assume that the indicator function $A(\cdot)$ verifies (\ref{Conditions-a-existence-1}) with $0<\beta<1$, and 
\begin{equation}\label{Conditions-a-existence-2}
    A(\cdot)\in \mathcal{C}^{0,1}(\Rt). 
\end{equation}
Moreover, assume  that the convolution term $\varphi$ verifies: 
\begin{equation}\label{Conditions-varphi-existence}
\varphi \in L^2(\T). 
\end{equation}

Then, for every $\alpha>0$ the coupled system (\ref{Alpha-model}) has a  weak Leray's solution  $(\vu,\vu_\alpha,P, P_\alpha)$ in the sense of Definition \ref{Def-Leray-Solutions}. In addition, for any time $t>0$, the filtered velocity $\vu_\alpha$ satisfies the following control:
\[ \| \vu_\alpha(t,\cdot)\|^2_{H^1} \lesssim \frac{1}{\min(\alpha^2 \beta,1)}  \|\vu(t,\cdot)\|^2_{L^2}. \]
\end{Theoreme}

We observe that the existence of weak Leray solutions is guaranteed when the indicator function $A(\cdot)$ is Lipschitz continuous. This condition on $A(\cdot)$ is often utilized in the study of elliptic equations of divergence form \cite{Zhang}, such as the nonlinear filter equation in (\ref{Alpha-model}). On the other hand, we also require the (technical) condition of square-integrability on the convolution term $\varphi$. As mentioned earlier, in Appendix \ref{Appendix}, we demonstrate that this latter condition can be omitted when proving the existence of weak Leray solutions 

\medskip

In contrast,  stronger regularity assumptions on $A(\cdot)$ and $\varphi$ are the crucial tools to prove their uniqueness.  

\begin{Theoreme}[Uniqueness]\label{Th-Uniqueness} Within the framework of Theorem \ref{Th-Existence}, let $\vu_{0,1}, \vu_{0,2} \in L^2_{\text{s}}(\T)$ be two initial data, and let $\vf_1, \vf_2 \in L^{2}_{loc}([0,+\infty[,\dot{H}^{-1}(\T))$ be two divergence-free and zero-mean external forces. Let  $(\vu_1, \vu_{1,\alpha}, P_1, P_{1,\alpha})$ and $(\vu_2, \vu_{2,\alpha}, P_2, P_{2,\alpha})$ be  the corresponding  weak Leray solutions of the coupled system (\ref{Alpha-model}). For any time $T>0$, define the energy norm
\[\| \vu_1 - \vu_2 \|_{E_T}=\sup_{0\leq t \leq T}\| \vu_1(t,\cdot)-\vu_2(t,\cdot)\|_{L^2}+\sqrt{\nu} \left( \int_{0}^{T}\|\vu_1(t,\cdot)-\vu_2(t,\cdot)\|^2_{\dot{H}^1}ds  \right)^{1/2}.\]

\medskip

In the filter equation of the system (\ref{Alpha-model}), assume additionally that  
\begin{equation}\label{conditions-a-varphi-uniqueness}
A(\cdot) \in \mathcal{C}^{1,1}(\Rt) \quad \mbox{and} \quad \varphi \in H^1(\T).
\end{equation}Then, the following estimates hold:
\begin{equation}\label{Dependence-date}
   \| \vu_1- \vu_2\|^{2}_{E_T} \leq  e^{\mathfrak{C}_0 T^{1/2}} \left( \| \vu_{0,1}-\vu_{0,2}\|^{2}_{L^2}\,  +\frac{4}{\nu}\int_{0}^{T} \| \vf_1(t,\cdot)-\vf_2(t,\cdot)\|^{2}_{\dot{H}^{-1}}dt \right),
\end{equation}
\begin{equation}\label{Dependence-data-Filter-equation}
\| \vu_{1,\alpha} - \vu_{2,\alpha}\|^2_{L^\infty_t H^1_x} \leq  \mathfrak{D}_0 \| \vu_1-\vu_2 \|^2_{L^\infty_t L^2_x},
\end{equation}
with   constants  $\mathfrak{C}_0>0$  and $\mathfrak{D}_0 >0$,  depending on the parameters $\alpha,\beta,\nu$, the functions $A(\cdot)$, $\varphi$, the data  $\vu_{0,1}, \vu_{0,2},\vf_1, \vf_2$, and the time $T$.
\end{Theoreme} 

In estimates (\ref{Dependence-date}) and (\ref{Dependence-data-Filter-equation}), we demonstrate a stronger result regarding the continuous dependence on data of two weak Leray solutions. These estimates ensure the uniqueness of the weak Leray solution associated with the same data. The uniqueness of the pressure terms $P$ and $P_\alpha$ follows from the uniqueness of the velocity  $\vu$, proven in (\ref{Dependence-date}); the uniqueness of the filtered velocity $\vu_\alpha$, proven in (\ref{Dependence-data-Filter-equation}); and the well-known characterization of the pressure:
\begin{equation}\label{Characterization-pressions}
\begin{split}
 &P= (-\Delta)^{-1} \text{div}\Big( \text{div} (\vu\otimes \vu_\alpha)\Big),\\
 &P_\alpha = -\alpha^2 (-\Delta)^{-1} \text{div}\Big( \text{div} \left( A (\varphi\ast\vu) (\vec{\nabla} \otimes \vu_\alpha)^T \right)\Big),
\end{split}
\end{equation}
where the operator $(-\Delta)^{-1}$ is classically defined in the Fourier variable as in expression (\ref{Lap-Frac}) below.

\medskip

{\bf Validity of the model: convergence properties to classical models}.  Now, we are interested in studying the asymptotic behavior of weak Leray solutions to the coupled system (\ref{Alpha-model}) with respect to variations of the parameter $\alpha>0$ and the indicator function $A(\cdot)$ in the  nonlinear  filter equation:
\[   -\alpha^2  \text{div} \left( A (\varphi\ast\vu) (\vec{\nabla} \otimes \vu_\alpha)^T \right)+\vu_\alpha +\vec{\nabla}P_\alpha = \vu, \quad \text{div}(\vu_\alpha)=0. \]
Specifically, on one hand, setting $\alpha=0$, we \emph{formally} observe that the system (\ref{Alpha-model}) reduces to the classical Navier-Stokes equations:
\begin{equation}\label{Navier-Stokes}
     \partial_t \vu - \nu\Delta \vu + (\vu \cdot \vec{\nabla})\vu + \vec{\nabla} P = \vf, \quad \text{div}(\vu)=0. 
\end{equation}
Indeed, when $\alpha = 0$, the filter equation yields the identity $\vu_\alpha + \vec{\nabla}P_\alpha = \vu$. Furthermore, by the second expression in (\ref{Characterization-pressions}), we (formally) deduce that $P_\alpha = 0$. Consequently, $\vu_\alpha = \vu$.  

\medskip

 On the other hand, for fixed $\alpha>0$, and recalling that the indicator function $A(\cdot)$  satisfies the lower and upper bounds in (\ref{Conditions-a-existence-1}) with a parameter $0<\beta<1$, as $\beta\to 1$, we formally obtain $A(\cdot)\equiv 1$, and we observe that the system (\ref{Alpha-model}) \emph{formally} becomes the classical Leray-alpha model:
\begin{equation}\label{Alpha-model-classical}
    \begin{cases}\vspace{2mm}
    \partial_t \vu - \nu\Delta \vu + (\vu_\alpha \cdot \vec{\nabla})\vu + \vec{\nabla} P = \vf, \quad \text{div}(\vu)=0, \\
     -\alpha^2 \Delta \vu_\alpha +\vu_\alpha = \vu, \quad \text{div}(\vu_\alpha)=0.
    \end{cases}
\end{equation}
In fact, when $A(\cdot)\equiv 1$,  the filter equation becomes  $-\alpha^2 \Delta \vu_\alpha+ \vu_\alpha + \vec{\nabla}P_\alpha=\vu$. Applying the divergence operator to both sides and using the divergence-free property of $\vu$ and $\vu_\alpha$, it follows that $\Delta P_\alpha=0$. Moreover, as $P_\alpha$ is a periodic function we deduce that it is a constant. Consequently, $\vec{\nabla}P_\alpha=0$.

\medskip

Thus, in our next results we rigorously  study the convergence of solutions of the coupled system (\ref{Alpha-model}) to solutions of the classical models (\ref{Navier-Stokes}) and (\ref{Alpha-model-classical}), as $\alpha \to 0$ and $\beta \to 1$, respectively. For simplicity, for any time $0<T<+\infty$ we denote the  spaces $E_T= L^\infty\big([0,T],L^2_\text{s}(\T)\big)\cap L^2\big([0,T],\dot{H}^1(\T)\big)$, $L^\infty_tH^1_x =L^\infty([0,T], H^1(\T))$,  $L^2_t L^2_x = L^2([0,T],L^2(\T))$ and $L^\infty([0,T],L^2(\T))=L^\infty_t L^2_x$.

\medskip

We start by proving a convergence property to solutions of the classical Navier-Stokes equations.  Under the same hypotheses as in Theorem \ref{Th-Existence}, let 
\[ \Big(\vu_{(\alpha)}, \vu_{(\alpha),\alpha},P_{(\alpha)}, P_{(\alpha),\alpha}\Big)_{\alpha>0}\subset E_T \times L^\infty_tH^1_x \times L^2_t L^2_x \times L^\infty_t L^2_x,\]
be the family of weak Leray solutions to the coupled system (\ref{Alpha-model}), which we rewrite here as:
\begin{equation}\label{Alpha-model-alpha} 
    \begin{cases}\vspace{2mm}
    \partial_t \vu_{(\alpha)} - \nu\Delta \vu_{(\alpha)} +  (\vu_{(\alpha),\alpha} \cdot \vec{\nabla})\vu_{(\alpha)} + \vec{\nabla}P_{(\alpha)}  =  \vf, \quad \text{div}(\vu_{(\alpha)})=0, \\ \vspace{2mm}
     -\alpha^2\, \text{div} \left( A (\varphi\ast\vu_{(\alpha)}) (\vec{\nabla} \otimes \vu_{(\alpha),\alpha})^T \right)+\vu_{(\alpha),\alpha} + \vec{\nabla}P_{(\alpha),\alpha}= \vu_{(\alpha)}, \quad \text{div}(\vu_{(\alpha),\alpha})=0, \\ 
\vu_{(\alpha)}(0,\cdot)=\vu_0.
    \end{cases}
\end{equation}
\begin{Proposition}[Convergence to the Navier-Stokes equations]\label{Prop1-Conv-NS} 
There exists a sub-sequence $(\alpha_k)_{k\in \mathbb{N}}$, where $\alpha_k \to 0$ as $k\to +\infty$, and  there exists $(\vu,P) \in E_T \times L^2_t \dot{H}^{-1/2}_x$, such that:
\begin{enumerate}
\item The velocities  $\vu_{(\alpha_k)}$ converge to $\vu$ in the \emph{weak-$*$ topology} of the space $E_T$, and in the \emph{strong topology} of the space $L^2_t L^2_x$. 
    \item The filtered velocities   $\vu_{(\alpha_k),\alpha_k}$ also  converge to $\vu$ in the \emph{strong topology} of the space $L^p_t L^2_x$, with $1\leq p <+\infty$.  
    \item The pressures  $P_{(\alpha_k)}$ converge to $P$ in the \emph{weak topology} of the space $L^2_t \dot{H}^{-1/2}_x$.
    \item  The pressures $P_{(\alpha_k),\alpha_k}$ converge to zero in the \emph{strong topology} of the space $L^\infty_t L^2_x$.
\end{enumerate}
Moreover, $(\vu,P)$ is a weak Leray solution to the classical Navier-Stokes equations (\ref{Navier-Stokes}).
\end{Proposition}

Now, we shall prove a convergence property to solutions of the  Leray-alpha model.  Within the setting  of Theorem \ref{Th-Existence} and for fixed $\alpha>0$,  let 
\[ \Big(\vu_{\beta}, \vu_{\beta,\alpha}, P_\beta, P_{\beta,\alpha}\Big)_{0<\beta<1}\subset E_T\times L^\infty_t H^1_x \times L^2_t L^2_x \times L^\infty_t L^2_x, \]
be the family of  weak  Leray solutions to the coupled system (\ref{Alpha-model}), which we also rewrite in the form:
\begin{equation}\label{Alpha-model-beta}
    \begin{cases}\vspace{2mm}
    \partial_t \vu_{\beta} - \nu\Delta \vu_{\beta} +  (\vu_{\beta,\alpha} \cdot \vec{\nabla})\vu_{\beta} + \vec{\nabla} P_{\beta}  =  \vf, \quad \text{div}(\vu_{\beta})=0, \\ \vspace{2mm}
     -\alpha^2\, \text{div} \left( A_\beta  (\varphi\ast\vu_{\beta}) (\vec{\nabla} \otimes \vu_{\beta,\alpha})^T \right)+\vu_{\beta,\alpha} + \vec{\nabla}P_{\beta,\alpha} = \vu_{\beta}, \quad \text{div}(\vu_{\beta,\alpha})=0, \quad \beta \leq A_\beta(\cdot)\leq 1, \\ 
\vu_{\beta}(0,\cdot)=\vu_0.
    \end{cases}
\end{equation}
\begin{Proposition}[Convergence to the Leray-alpha model]\label{Prop2-Conv-NSalpha}
There exists a sub-sequence $(\beta_k)_{k \in \mathbb{N}}$, where $\beta_k \to 1$ as $k\to +\infty$, and there exists $(\vu, \vu_\alpha,P) \in E_T\times L^\infty_t H^1_x \times L^2_t L^2_x$, such that:
\begin{enumerate}
\item The velocities  $\vu_{(\beta_k)}$ converge to $\vu$ in the \emph{strong topology} of the space $E_T$.
    \item The filtered velocities   $\vu_{\beta_k,\alpha}$  converge to the filtered velocity $\vu_\alpha$ in the \emph{strong topology} of the space $L^\infty_t H^1_x$.
    \item The pressures  $P_{\beta_k}$ converge to $P$ in the \emph{strong topology} of the space $L^2_t L^2_x$.
    \item The terms $\vec{\nabla}P_{\beta_k,\alpha}$ converge to zero in the \emph{strong topology}  of the space $L^\infty_t \dot{H}^{-1}_x$.
 \end{enumerate}
Moreover, $(\vu,\vu_\alpha,P)$ is a weak Leray  solution of the  Leray-alpha model  (\ref{Alpha-model-classical}). 
\end{Proposition}

{\bf Long-time dynamics}. We are concerned with the study of the long-time behavior of weak Leray solutions to the coupled system (\ref{Alpha-model}). For this, we set our framework as follows: first, we fix a stationary (time-independent) divergence-free external force $\vf \in \dot{H}^{-1}(\T)\cap L^2_{\text{s}}(\T)$. Next, since the pressure terms $P$ and $P_\alpha$ do not play any substantial role in this study, from now on, we shall consider the system
\begin{equation}\label{Alpha-model-Proy}
    \begin{cases}\vspace{2mm}
    \partial_t \vu - \nu\Delta \vu + \P \big( (\vu_\alpha \cdot \vec{\nabla})\vu \big)  =  \vf, \quad \text{div}(\vu)=0,  \\ \vspace{2mm}
     -\alpha^2 \, \P \Big( \text{div} \left( A (\varphi\ast\vu) (\vec{\nabla} \otimes \vu_\alpha)^T \right)\Big)+\vu_\alpha = \vu, \quad \text{div}(\vu_\alpha)=0, \quad \alpha >0,  \\ 
\vu(0,\cdot)=\vu_0, \quad \text{div}(\vu_0)=0. 
    \end{cases}
\end{equation}
Moreover, in the nonlinear filter equation above,    we shall assume the parameter $\alpha>0$ fixed, and we shall assume that  the indicator function $A(\cdot)$ and the convolution term $\varphi$ verify the conditions (\ref{conditions-a-varphi-uniqueness}):
\begin{equation*}
A(\cdot) \in \mathcal{C}^{1,1}(\Rt) \quad \mbox{and} \quad \varphi \in H^1(\T).
\end{equation*}
Thus, by Theorems \ref{Th-Existence} and \ref{Th-Uniqueness},  the system (\ref{Alpha-model-Proy}) defines a  \emph{continuous semi-group} acting on the Hilbert space $L^2_{\text{s}}(\T)$ as follows:  for any time $t\geq 0$, we define the semi-group 
 \begin{equation}\label{Semi-group}
     S(t):L^2_{\text{s}}(\T)\to L^2_{\text{s}}(\T), \quad \vu_0 \mapsto S(t)\vu_0= \vu(t,\cdot),
 \end{equation}
 where $\vu(t,\cdot)$ denotes the \emph{unique} weak  Leray solution to (\ref{Alpha-model-Proy}) arising from the initial datum $\vu_0$.  

 \medskip

The study of the long-time dynamics of weak Leray solutions of the coupled system (\ref{Alpha-model-Proy}) can be approached through the analysis of the semi-group $S(t)$ as $t\to +\infty$.  More precisely, our aim is to prove the existence of a global attractor, the definition of which we recall as follows. See \cite{Babin,Chepyzhov,Temam} for more details. 
\begin{Definition}\label{Def-Global-Attractor}  A global attractor for the semi-group $(S(t))_{t\geq 0}$ defined in (\ref{Semi-group}) is a set $\mathcal{A}\subset L^2_{\text{s}}(\T)$ which satisfies:
 \begin{enumerate}
     \item[(A)] The set $\mathcal{A}$ is compact in the strong topology of $L^2_\sigma(\T)$, which is the same as the strong topology of $L^2(\T)$.
     \item[(B)] The set $\mathcal{A}$ is strictly invariant under $S(t)$: for any time $t\geq 0$, we have $S(t)\mathcal{A}=\mathcal{A}$.
     \item[(C)] $\mathcal{A}$ is an attracting set: for any bounded set $B \subset L^2_{\text{s}}(\T)$ and for any neighborhood $\mathcal{O}$ of $\mathcal{A}$, there exists a time $T=T(B,\mathcal{O})>0$ (which depends on the set $B$ and the neighborhood $\mathcal{O}$) such that for all $t>T$,   it holds $S(t)B\subset \mathcal{O}$. 
 \end{enumerate}
 \end{Definition}
By point (C), roughly speaking, a global attractor attracts the image through the semi-group $S(t)$ of all bounded sets $B \subset L^2_{\text{s}}(\T)$ as $t\to +\infty$. This property allows us to gain a better understanding of the long-time behavior of weak Leray solutions for the alpha-model (\ref{Alpha-model-Proy}). Indeed, we observe first that for all initial data $\vu_0 \in L^2_{\text{s}}(\T)$, setting the bounded set $B=\{ \vu_0 \}$, and moreover, for any given neighborhood $\mathcal{O}$ of the attractor $\mathcal{A}$, we find that the weak Leray solution $\vu(t,\cdot)$  of (\ref{Alpha-model-Proy}) (arising from $\vu_0$)  lies in the neighborhood $\mathcal{O}$ from some time $ T=T(\vu_0, \mathcal{O})$. Consequently, for any initial datum $\vu_0$, the solution $\vu(t,\cdot)$ approaches the attractor $\mathcal{A}$ as closely as desired as $t\to +\infty$.

\medskip

With this last idea in mind, it is also interesting to characterize the global attractor $\mathcal{A}$. Precisely, one asks for the elements belonging to $\mathcal{A}$. These elements are defined through a particular type of solutions of our alpha-model, so-called the \emph{eternal solutions}.
\begin{Definition}\label{Def-Eternal-Solutions} Let $\vf \in \dot{H}^{-1}(\T)\cap L^2_{\text{s}}(\T)$ be the stationary external force. We shall say that $\vu_{\text{e}}$  is an eternal solution to the  Navier-Stokes-alpha model if:  
 \begin{enumerate}
    \item[(A)] We have $\ds{\vu_{\text{e}} \in L^{\infty}_{loc}\big(\R, L^2_{\text{s}}(\T)\big) \cap L^{2}_{loc}\big(\R,\dot{H}^1(\T)\big)}$. 
 \item[(B)] The vector field $\vu_\text{e}$ solves the following   coupled system  in the weak sense:
 \begin{equation}\label{Alpha-model-eternal}
\begin{cases}\vspace{2mm}
\partial_t \vu_\text{e} - \nu\Delta \vu_\text{e} + \P\big( (\vu_{\text{e},\alpha} \cdot \vec{\nabla})\vu_\text{e}\big)  = \vf, \quad \text{div}(\vu_\text{e})=0, \\ 
-\alpha^2 \, \P\Big( \text{div} \left( A (\varphi \ast \vu_\text{e}) (\vec{\nabla} \otimes \vu_{\text{e},\alpha})^T \right)\Big)+\vu_{\text{e},\alpha} = \vu_\text{e}, \quad \text{div}(\vu_{\text{e},\alpha})=0.
\end{cases}
\end{equation} 
 \end{enumerate}
 Moreover, when $\vu_{\text{e}}\in L^{\infty}\big(\R, L^2_{\text{s}}(\T)\big)$, it is called a bounded eternal solution. 
 \end{Definition}
As observed, eternal solutions only depend on the given external force $\vf$ and they are defined on the whole real line $\R$. Consequently, they do not arise from any initial data. In Lemma \ref{Lemma-Eternal-Sol} below, we prove the existence of these solutions. 

\medskip

Once we have introduced these definitions, our next result states as follows: 
\begin{Theoreme}[Existence of a unique global attractor]\label{Th-Global-Attractor} Assume (\ref{conditions-a-varphi-uniqueness}), and let $\vf \in \dot{H}^{-1}\cap L^2_{\text{s}}(\T)$ be the stationary external force. Then, the semi-group $(S(t)){t\geq 0}$ defined in (\ref{Semi-group}) has a unique \emph{global attractor} $\mathcal{A}=\mathcal{A}(\vf)$ in the sense of Definition \ref{Def-Global-Attractor}. Moreover, we have the following characterization:
\begin{equation}\label{Characterization-global-attractor}
\mathcal{A} = \left\{ \vu_\text{e}(0,\cdot) \in L^2_{\text{s}}(\T) : \ \ \text{where the \emph{bounded} eternal solution $\vu_\text{e}$ is given in Definition \ref{Def-Eternal-Solutions}} \ \right\}.
\end{equation}
\end{Theoreme}
The following comments are in order here. First, by expression (\ref{Characterization-global-attractor}), we observe that the global attractor $\mathcal{A}$ is explicitly characterized as the set of bounded eternal solutions at time $t=0$. Additionally, recall that all these eternal solutions depend on the given external force $\vf$. Thus, we write $\mathcal{A}=\mathcal{A}(\vf)$ to emphasize the dependence of  $\mathcal{A}$ on $\vf$. 

\medskip

The proof of this theorem relies on earlier results from dynamical systems theory, detailed in Section \ref{Sec:Existence-Global-Attractor}. Essentially, the key approach to establishing the existence of the global attractor $\mathcal{A}$ involves demonstrating that the semi-group $S(t)$ exhibits favorable asymptotic behavior for sufficiently large $t$ (refer to Theorem \ref{Th-Tech-Existence-Global-Attractor}). One of the primary challenges is ensuring that this semi-group is asymptotically compact, as defined in Definition \ref{Def-Asymp-Compact}. To establish this property, we adapt energy methods previously employed in \cite{Ilyn} for the 2D damped Navier-Stokes equations and in \cite{CortezJarrin} for the 3D damped Navier-Stokes-Bardina  model. 

\medskip

Our last result studies an explicit upper bound on the \emph{fractal dimension} of the global attractor $\mathcal{A}$. The fractal dimension of the compact set $\mathcal{A}$ offers a precise measure of its size. For further details, please refer to \cite{Babin,Temam}.  \begin{Definition} The fractal dimension of $\mathcal{A}$ is defined by the quantity: 
\begin{equation}
  \text{dim}_f (\mathcal{A})=  \limsup_{\varepsilon\to 0^{+}} \frac{N_\varepsilon(\mathcal{A})}{\ln(1/\varepsilon)},
\end{equation}
where $N_\varepsilon(\mathcal{A})>0$ denotes the minimal number of $\varepsilon$-balls in $L^2_\sigma(\T)$ needed to cover $\mathcal{A}$. 
\end{Definition}

Our goal is to demonstrate that the quantity $\text{dim}_f (\mathcal{A})$  is finite and establish an explicit upper bound for it. This upper bound depends on a non-trivial set of quantities defined by all the parameters of our model. Please refer to estimate (\ref{Estim-Tech-Fractal-Dim}) below. To simplify the presentation of our next result, and following some concepts from \cite{Cao, CheskidovHolm, Chen}, we will focus solely on quantities related to the norm of the external force $\| \vf \|_{\dot{H}^{-1}}$,  as well as the parameters $\alpha>0$ and $0<\beta<1$  in the nonlinear filter equation of system (\ref{Alpha-model-Proy}). 

\medskip

We thus introduce the following quantity: 
\begin{equation}\label{K-alpha-beta}
K(\alpha,\beta)= \begin{cases}\vspace{2mm}
    \frac{1}{\alpha^5 \beta^{5/2}}, \quad 0< \alpha \leq 1, \\ \vspace{2mm}
    \frac{1}{\alpha  \beta^{5/2}}, \quad 1<\alpha \leq \frac{1}{\sqrt{\beta}}, \\
    \alpha^4, \quad \frac{1}{\sqrt{\beta}}<\alpha, 
\end{cases}    
\end{equation} 
which is defined solely by $\alpha$ and $\beta$. Together with  the quantity $\| \vf\|_{\dot{H}^{-1}}$, we can state:
\begin{Theoreme}[Upper bounds on the fractal dimension]\label{Th-Fractal-Dimension} Let $\mathcal{A}\subset L^2_{\text{s}}(\T)$ be the unique global attractor obtained  by Theorem \ref{Th-Global-Attractor}.  The following estimate holds:
\begin{equation}\label{Estimate-fractal-dimension}
\begin{split}
    \text{dim}_{f}(\mathcal{A})\lesssim  &\,  \left(1+\sqrt{K(\alpha,\beta)}(1+\| \vf \|_{\dot{H}^{-1}})^3 e^{K(\alpha,\beta)(1+\| \vf \|_{\dot{H}^{-1}})^{3/2}} \right)^3 \\
    &\, \times \, \ln\left( 1+ K(\alpha,\beta)(1+\| \vf \|_{\dot{H}^{-1}})^6\, e^{K(\alpha,\beta)(1+\| \vf \|_{\dot{H}^{-1}})^3} \right).
    \end{split}
\end{equation}
\end{Theoreme}

In the case of the \emph{linear} Helmholtz filter equation (\ref{Helmholtz-Filter}), upper bounds on the fractal dimension of the global attractor associated with the classical Leray-alpha model (\ref{Alpha-model-classical}) were derived in \cite{Cao} using a \emph{linearized} version of (\ref{Alpha-model-classical}). The key approach involves linearizing the nonlinear transport term in (\ref{Alpha-model-classical}) around a trajectory within the global attractor. This linearization is possible by the identity $(\vu_1+\vu_2)_\alpha = \vu_{1,\alpha} + \vu_{2,\alpha}$, which applies to two filtered velocities obtained from (\ref{Helmholtz-Filter}). However, due to the \emph{nonlinear} nature of the filter equation in the system (\ref{Alpha-model-Proy}), this method does not appear applicable in our scenario. To overcome this problem, we pursue a different approach to derive an upper bound for the fractal dimension $ \text{dim}_{f}(\mathcal{A})$. The methodology developed in \cite{Chueshov} for \emph{hyperbolic} second-order evolution equations can be adapted to the \emph{parabolic-elliptic} framework of the system (\ref{Alpha-model-Proy}), leading to the estimate (\ref{Estimate-fractal-dimension}). For detailed technical information, please refer to Section \ref{Sec:Fractal-Dimension}. 

\medskip

In estimate (\ref{Estimate-fractal-dimension}), as observed, we note that the quantity $ \text{dim}_{f}(\mathcal{A})$ is controlled by the magnitude of the external force $\| \vf \|_{\dot{H}^{-1}}$ and the parameters $\alpha$ and $\beta$. Specifically, by fixing $\alpha$ and $\beta$, where for simplicity we denote $K(\alpha,\beta)\simeq 1$, we obtain
\begin{equation*}
\begin{split}
    \text{dim}_{f}(\mathcal{A})\lesssim  &\,  (1+\| \vf \|_{\dot{H}^{-1}})^9\,  e^{3(1+\| \vf \|_{\dot{H}^{-1}})^{3/2}} \\
    &\, \times \, \ln\left((1+\| \vf \|_{\dot{H}^{-1}})^6\, e^{(1+\| \vf \|_{\dot{H}^{-1}})^3} \right),
    \end{split}
\end{equation*}
which indicates that the right-hand side of (\ref{Estimate-fractal-dimension}) exponentially grows for large values of  $\| \vf \|_{\dot{H}^{-1}}$. Conversely, when fixing $\| \vf \|_{\dot{H}^{-1}}$, where for simplicity  we  write $1+\| \vf \|_{\dot{H}^{-1}} \simeq 1$, once again,  by (\ref{Estimate-fractal-dimension}) we derive
\begin{equation*}
\begin{split}
    \text{dim}_{f}(\mathcal{A})\lesssim  &\,  \left(1+\sqrt{K(\alpha,\beta)} \, e^{K(\alpha,\beta)} \right)^3 \\
    &\, \times \, \ln\left( 1+ K(\alpha,\beta)\, e^{K(\alpha,\beta)} \right).
    \end{split}
\end{equation*}
Here, according to the definition of the quantity $K(\alpha,\beta)$ given in (\ref{K-alpha-beta}), we make the following observations.  From the first expression in (\ref{K-alpha-beta}), for fixed $0<\beta<1$ and the fact that $K(\alpha,\beta)\, e^{K(\alpha,\beta)}=\frac{1}{\alpha^5 \beta^{5/2}} e^{\frac{1}{\alpha^5 \beta^{5/2}}}$, we note that the upper bound of $ \text{dim}_{f}(\mathcal{A})$ increases like $1/\alpha^5$ for small values of $\alpha$.  Next, from the second expression in (\ref{K-alpha-beta}), for fixed $1<\alpha\leq \frac{1}{\sqrt{\beta}}$ and the fact that $K(\alpha,\beta)\, e^{K(\alpha,\beta)}=\frac{1}{\alpha \beta^{5/2}} e^{\frac{1}{\alpha \beta^{5/2}}}$, we obtain that the upper bound of  $ \text{dim}_{f}(\mathcal{A})$ increases like $1/\beta^{5/2}$ for small values of $\beta$. Finally, from the third expression in (\ref{K-alpha-beta}), we see that $ \text{dim}_{f}(\mathcal{A})$ is controlled solely by the parameter $\alpha$, and its upper bound also grows exponentially for large values of  $\alpha$.

\medskip


{\bf Future research and open questions}.  In a forthcoming article, we extend the research program $(i)-(iv)$ on the new Navier-Stokes-Alpha model (\ref{Alpha-model}) in the following directions. First, we study the higher regularity properties of weak Leray solutions, provided that the indicator function $A(\cdot)$ verifies suitable regularity conditions. By following some of the ideas from \cite{Ilyn}, we also aim to show that these regularity properties of weak Leray solutions yield higher regularity estimates on the global attractor $\mathcal{A}$. For instance, we expect to prove that $\mathcal{A}$ is also bounded in the Sobolev space $H^s(\T)$ with $s > 0$. Thereafter, we  address the long-time behavior of weak Leray solutions through the notion of stationary (time-independent) solutions. These solutions depend solely on the spatial variable and consequently formally solve the doubly elliptic nonlinear problem:
\begin{equation*}
\begin{cases}
- \nu\Delta \U + (\U_\alpha \cdot \vec{\nabla})\U + \vec{\nabla} P = \vf , \quad \text{div}(\U)=0, \\
-\alpha^2 \text{div} \left( A (\varphi \ast\U) (\vec{\nabla} \otimes \U_\alpha)^T \right)+\U_\alpha +\vec{\nabla}P_\alpha = \U, \quad \text{div}(\U_\alpha)=0.
\end{cases}
\end{equation*}
We thus prove the existence of these solutions in a natural energy space and study their orbital and asymptotic stability. 

\medskip

On the other hand, as mentioned, uniqueness of weak Leray solutions to the alpha-model (\ref{Alpha-model-Intro}) with a non-regularized filter equation cannot be obtained by following the proof of Theorem \ref{Th-Uniqueness}, since the involved quantities $\mathfrak{C}_0 > 0$ and $\mathfrak{D}_0 > 0$ explicitly depend on $\| \varphi \|_{H^1}$. In this context, a natural first approach is to prove a weak-strong uniqueness result by considering, for instance, the classical Prodi-Serrin criterion \cite{Prodi,Serrin} or more sophisticated \emph{a priori} conditions \cite{Fernandez-Jarrin,PLe0}. 

\medskip

{\bf Organization of the article}. In Section \ref{Sec:Example}, we provide an example of a convolution term $\varphi$ that formally satisfies $\varphi \ast \vu \approx \vu$. Section \ref{Sec:Filter-Equation} is devoted to studying  the nonlinear filter equation in the coupled system (\ref{Alpha-model}). In  Section \ref{Sec:GWP}, we address  the global well-posedness of  (\ref{Alpha-model}), and in Section \ref{Sec:Convergences}, we prove some convergence properties related to classical models. In Section \ref{Sec:Global-Attractor}, we study the long-time behavior of weak Leray solutions to (\ref{Alpha-model}). Finally, in Appendix \ref{Appendix}, we prove the existence of these solutions for the coupled system (\ref{Alpha-model-Intro}). 

\medskip


{\bf Preliminaries and notation}. We summarize here some standard definitions and notation in the framework of periodic boundary conditions $\T=[0,  L]^3$. As usual, we shall assume that both the initial datum $\vu_0$ and the external force $\vf$ are \emph{zero-mean} vector fields: $\ds{\int_{\T}\vu_0(x)dx=\int_{\T}\vf(x)dx=0}$ which yields that the velocity $\vu(t,x)$ is also a \emph{zero-mean} vector field for all $t\geq 0$. For simplicity, we shall denote ${\bf Z}^3:=\frac{2\pi}{L}\mathbb{Z}^3$,  and we we will frequently use the Fourier decomposition 
\begin{equation}\label{Fourier-decompositio}
\vu(t,x)= \sum_{k \in \Z} e^{i k \cdot x}\,  \widehat{\vu}_k(t,\cdot),    \end{equation}
with  spatially dependent Fourier coefficients 
\begin{equation}\label{Fourier-coefficients}
\widehat{\vu}_k(t,\cdot) := \frac{1}{L^3} \int_{\T} e^{-i k\cdot x} \vu(t,x)dx.    
\end{equation}

We denote the usual Lebesgue spaces by $L^p(\T)$ for $1\leq p \leq +\infty$, and  $(\cdot , \cdot)_{L^2}$ stands for the $L^2-$inner product. We define $\dot{L}^2(\T)$ as the space of zero-mean $L^2$-functions.  By Helmholtz decomposition (see \cite[Theorem $2.6$]{Wiedemann}), we have $\dot{L}^2(\T)=L^2_{\text{s}}(\T) \bigoplus (L^{2}_{\text{s}})^{\perp}(\T)$, where $(L^{2}_{\text{s}})^{\perp}(\T)$ is the space of $\dot{L}^2$- vector fields $\vg$ such that $\vg = \vec{\nabla}\phi$ with $\phi \in \dot{H}^1(\T)$.

\medskip

For $s\in \R$,   homogeneous Sobolev spaces  of zero-mean functions are classically characterized by the Fourier series as follows: 
\begin{equation*}
\dot{H}^s(\T)=\Big\{  g(x)=\sum_{k\in \Z } e^{i k \cdot x}\,\widehat{g}_k :  \ \overline{\widehat{g}_k}=\widehat{g}_{-k}, \ \widehat{g}_0=0, \ \mbox{and}\ \| g \|^2_{\dot{H}^s}:= L^3\sum_{k \in \Z } |k|^{2s} |\widehat{g}_k|^2<+\infty \Big\}.
\end{equation*}
Here, the condition $\overline{\widehat{g}_k}=\widehat{g}_{-k}$ means that we only work with real-valued functions. Remark that $\dot{H}^0(\T)=\dot{L}^2(\T)$, and for any $s_1<s_2$ we have the \emph{compact embedding} $\dot{H}^{s_2}(\T)\subset \dot{H}^{s_1}(\T)$.
Moreover, non-homogeneous Sobolev spaces are similarly defined as: 
\begin{equation*}
H^s(\T)=\Big\{  g(x)=\sum_{k\in \Z} e^{i k \cdot x}\,\widehat{g}_k :  \ \overline{\widehat{g}_k}=\widehat{g}_{-k}, \ \widehat{g}_0=0, \ \mbox{and}\ \| g \|^2_{H^s}:=L^3\sum_{k \in \Z} (1+|k|^2)^s |\widehat{g}_k|^2<+\infty \Big\}.
\end{equation*}

Thereafter, $\P:\dot{L}^2(\T)\to L^2_{\text{s}}(\T)$ is the well-known Leray's projector, which is  defined by the following expression   $\ds{\P(\vg)(x)=\sum_{k\in \Z} e^{i k\cdot x} \left(\widehat{\vg}_k - \frac{\widehat{\vg}_k\cdot k}{|k|^2} k \right)}$. Finally, fractional powers $s \in \R$
of the Laplacian operator are defined by the expression
\begin{equation}\label{Lap-Frac}
(-\Delta)^{\frac{s}{2}}g(x)=\sum_{k\in \Z}  e^{i k \cdot x}\,|k|^s\widehat{g}_k.
\end{equation}

\section{An example of a convolution term}\label{Sec:Example}
We provide a simple example of a convolution term $\varphi$ that formally satisfies $\varphi \ast \vu \approx \vu$. We emphasize that, in contrast to the other estimates in this article, this approximation is not rigorously proven. As already explained in the introduction, our aim is to shed light on the fact that, although the convolution term $\varphi$ is a technical requirement in this theoretical study of the coupled system (\ref{Alpha-model}), it seems possible to design more advanced convolution terms $\varphi$ in which the local filtering effects of the indicator function $A(\cdot)$ on the velocity $\vu$ remain relatively unaffected by $\varphi$.

\medskip

Using the Fourier decomposition introduced above, we define
\begin{equation*}
\varphi(x)= \sum_{k \in \Z} e^{i k\cdot x} \widehat{\varphi}_k,
\end{equation*}
where, for a fixed parameter $\kappa > 0$, we have
\begin{equation*}
\widehat{\varphi}_k = \begin{cases}
1, & \text{if} \quad  |k| \leq \kappa, \\
0, & \text{if} \quad  |k| > \kappa.
\end{cases}
\end{equation*}

Mathematically, $\varphi$ corresponds to an approximation of the periodic Dirac mass $\delta_0$. Therefore, using the Fourier decomposition of the velocity $\vu$ given in (\ref{Fourier-decompositio})-(\ref{Fourier-coefficients}), and by well-known properties of the Fourier coefficients, we obtain
\begin{equation*}
\varphi \ast \vu (t,x)= \sum_{k \in \Z} e^{i k\cdot x} \widehat{(\varphi \ast \vu(t,\cdot))_k} = \sum_{k \in \Z} e^{i k\cdot x} \widehat{\varphi}_{k} \widehat{\vu}_{k}(t,\cdot) = \sum_{k \in \Z, \ |k| \leq \kappa} e^{i k\cdot x} \widehat{\vu}_{k}(t,\cdot).
\end{equation*}
We thus note that $\varphi \ast \vu$ cuts off the Fourier decomposition of $\vu$ at high frequencies $|k| > \kappa$.

\medskip

The cut-off frequency $\kappa > 0$ can be set using different criteria. On the one hand, some numerical studies of turbulence models \cite{Geurts0,Kuerten} apply discrete Fourier transform techniques to truncate the velocity $\vu$ at high-frequency wave numbers. In this way, the cut-off frequency $\kappa > 0$ can be chosen based on machine accuracy criteria so that the truncated velocity $\ds{\sum_{k \in \Z, \ |k| \leq \kappa}} e^{i k\cdot x} \widehat{\vu}_{k}(t,\cdot)$ retains a significant amount of information about the velocity $\vu$. In this context, we formally write $\varphi \ast \vu \approx \vu$.

\medskip

On the other hand, the cut-off frequency $\kappa>0$ can be also set by following some ideas of the  conventional Kolmogorov's  theory of turbulence  \cite{Jacquin}. To begin with, let's clarify the basic concept behind Kolmogorov's dissipation law, known as the \emph{energy cascade model}. This model describes how kinetic energy is introduced into the fluid by an external force $\vf$ at a length scale $\ell_0 > 0$, referred to as the \emph{energy input scale}. Given that we can control the external force, this input scale $\ell_0$ is always fixed. In a turbulent regime characterized by a large Reynolds number $Re$ or Grashof number $Gr$ (defined below), the energy dissipation due to viscosity is ineffective at this input scale $\ell_0$, leading to a transfer of energy to progressively smaller scales. This transfer of energy is physically realized through vortex stretching, where eddies at a scale $\ell_1 < \ell_0$ break down into even smaller eddies at scales $\ell_2 < \ell_1$. This \emph{cascade} of energy continues until reaching the Kolmogorov dissipation scale $\ell_D$, below which the kinetic energy from larger scales is eventually dissipated by molecular viscosity. Consequently, the \emph{inertial range} is defined as the interval of length scales $]\ell_D,\ell_0[$ where kinetic energy is transferred. For more details, please refer to \cite{Doering-Foias2}.

\medskip

In the Fourier domain, we define $\kappa_0=\frac{1}{\ell_0}$ the energy input frequency, and  we define $\kappa_D = \frac{1}{\ell_D}$ as the frequency above which the kinetic energy is dissipated. Specifically, considering the Fourier decomposition of the velocity $\vu$ given in expression (\ref{Fourier-decompositio}), in a turbulent regime with a large Reynolds number $Re$, it is experimentally observed \cite{Foias,Jacquin} that the Fourier coefficients $\widehat{\vu}_k(t,\cdot)$ defined in expression (\ref{Fourier-coefficients})) behave as 
\begin{equation}\label{Energy-dissipation}
|\widehat{\vu}_k(t,\cdot)| \simeq e^{-|k|}, \quad    |k| > \kappa_D. 
\end{equation}

\medskip

The frequency $\kappa_D$ is so-called the Kolmogorov's dissipation wave number, and it is  characterized by the expression  $\kappa_D \simeq Re^{1/2}\kappa_0$. Here,  the Reynolds number $Re$ is defined as follows: from the velocity $\vu$, which is a Leray weak solution of the coupled system (\ref{Alpha-model}), and the period $L$,  following \cite{Doering-Foias} we introduce  the fluid characteristic velocity $U:=\left(\ds{\lim_{T\to +\infty}}\frac{1}{T}\int_{0}^{T}\| \vu(t,\cdot)\|^2_{L^2} \frac{dt}{L^3}\right)^{1/2}$. Then, with the viscosity constant $\nu>0$, we define  $\ds{Re= \frac{UL}{\nu}}$. In addition,  in a flow driven by a fixed body force, the characteristic velocity $U$, and  consequently, the Reynolds number $Re$ are not directly controllable. Instead, what can be controlled is the amplitude of the force $\ds{F := \frac{\| \vf \|_{L^2}}{L^{3/2}}}$ and its dimensionless counterpart, the Grashof number $\ds{Gr = \frac{F L^3}{\nu^2}}$, which can be predetermined in advance. In addition, in the turbulent regime one has $(Re)^2\simeq Gr$, hence we obtain $\kappa_D \simeq Gr\kappa_0$.

\medskip

In this setting, for a given energy input frequency $\kappa_0 > 0$ and the Grashof number $Gr$ defined above, we set $\kappa = Gr \kappa_0$. Since a turbulent fluid, characterized by large values of the Grashof number $Gr$, is expected to dissipate energy rapidly at high frequencies, we formally write, by (\ref{Energy-dissipation})
\begin{equation*}
    \begin{split}
  \vu(t,x) =&\,  \sum_{k \in \Z}e^{i k\cdot x} \widehat{\vu}_k(t,\cdot)=\sum_{k \in \Z, \ |k|\leq \kappa }e^{i k\cdot x} \widehat{\vu}_k(t,\cdot) + \sum_{k \in \Z, \ |k|> \kappa }e^{i k\cdot x} \widehat{\vu}_k(t,\cdot)\\
  =&\, \varphi \ast \vu (t,x)+ \sum_{k \in \Z, \ |k|> \kappa }e^{i k\cdot x} \widehat{\vu}_k(t,\cdot) \approx \varphi\ast \vu(t,x) + \sum_{k \in \Z, \ |k|> \kappa }e^{i k\cdot x} e^{-|k|}\approx \varphi \ast \vu(t,x).      
    \end{split}
\end{equation*}

\section{The regularized nonlinear filter equation with $L^2$-data}\label{Sec:Filter-Equation} 
In this section, for any fixed $\alpha>0$ and $0<\beta <1$, we study the elliptic nonlinear filter equation in the system (\ref{Alpha-model}):
\begin{equation}\label{Filter-Regularized}
-\alpha^2 \text{div} \left( A (\varphi * \vu) (\vec{\nabla} \otimes \vu_\alpha)^T \right) + \vu_\alpha + \vec{\nabla}P_\alpha = \vu, \quad \text{div}(\vu_\alpha) = 0, \quad \beta \leq A(\cdot) \leq 1,
\end{equation}
with periodic boundary conditions in the torus $\T$, and where, for any time $t>0$, the velocity $\vu(t,\cdot)\in L^2_{\text{s}}(\T)$  is considered as a \emph{given datum}. In forthcoming propositions, we prove the existence of weak  $H^1$-solutions to equation (\ref{Filter-Regularized}), continuous dependence of data, and higher order $H^2$-estimates, which are the key tools in the next sections.

\medskip

We emphasize that in the particular case when $A(\cdot)\equiv 1$, these are well-known results for the linear Helmholtz filter equation (\ref{Helmholtz-Filter}) in bounded smooth domains $\Omega$  and the whole space $\Rt$. See, for instance, \cite[Chapter IX]{Brezis} and \cite[Chapter 6.3]{Evans}. Additionally, these are also well-known results in the case when the indicator function $A(\cdot)$ does not depend on the datum $\vu$, and we deal with the equation
\[-\alpha^2 \text{div}(A(x)(\nabla \otimes \vu_\alpha))^T)+\vu_\alpha+\vec{\nabla}P_\alpha=\vu,\]
with additional suitable regularity properties assumed on $A(\cdot)$. See, for instance, \cite[Chapter 9]{Gilbarg} and \cite[Chapter 8]{Ya}.

\medskip

In our case, the presence of the datum $\vu$ in the expression $A(\varphi\ast \vu)$ of the divergence term in equation (\ref{Filter-Regularized}) and the nonlinear character of on the indicator function $A(\cdot)$ make this equation nonlinear, and it does not allow us to directly apply these previous results to this equation. Additionally, as already explained, in the setting of the whole model (\ref{Alpha-model}), our datum $\vu$ is constrained to being only an $L^2$-function, and any supplementary hypotheses (such as continuity and regularity) cannot be assumed. 

\medskip

For the completeness of this article, we revisit classical results from \cite{Brezis,Evans,Gilbarg,Ya} to study the existence and regularity of solutions to equation (\ref{Filter-Regularized}) (as given in Propositions \ref{Existence-Filter} and \ref{Reg-H2-Filter}). We also establish new estimates concerning the continuous dependence of data (as stated in Propositions \ref{Dependence-data-Filter} and \ref{Cont-Dep-H2-Filter}). We emphasize the strong dependence of these results on suitable assumptions for the indicator function $A(\cdot)$ and the convolution term $\varphi$.

\medskip

All the estimates that we will prove below involve quantities that depend on the parameters $\alpha,\beta$  and some norms of the functions $A(\cdot)$ and $\varphi$. We will explicitly write out expressions for these quantities since they will be useful in what follows. However, for simplicity in our exposition, we will omit explicit mention of generic constants, employing a slight abuse of notation. 

 \begin{Proposition}[Existence of $H^1$-solutions]\label{Existence-Filter} Let fixed $\alpha>0$. For a time $0<T<+\infty$,  let $\vu \in L^{\infty}([0,T],L^2_{\text{s}}(\T))$.  Assume that the function $A(\cdot)$ verifies (\ref{Conditions-a-existence-1}) with $0<\beta<1$. Then, there exists 
 \[ (\vu_\alpha, P_\alpha) \in L^\infty([0,T],H^1(\T)) \times L^\infty([0,T],L^2(\T)), \]
  a  weak solution to equation  (\ref{Filter-Regularized}). Moreover,  for all $0\leq t \leq T$, this solution verifies the control:
\begin{equation}\label{Control-Filter}
  \| \vu_{\alpha}(t,\cdot)\|^{2}_{H^1} \leq  K_0\| \vu(t,\cdot)\|^{2}_{L^2}, \quad \mbox{with} \quad  K_{0} \simeq  \frac{1}{\min(\alpha^2 \beta, 1)}>0. 
\end{equation} 
\end{Proposition} 
\pv First of all, recall that $\mathcal{V}$ is defined as the set vector trigonometric polynomials verifying divergence-free and zero-mean conditions. We then define the space ${H}^{1}_{\text{s}}(\T)$ as the  closure of  $\mathcal{V}$ in ${H}^1(\T)$.  

\medskip

In equation (\ref{Filter-Regularized}), we apply the Leray's projector $\P$  to obtain $\ds{ -\alpha^2 \P\, \text{div}\Big( A(\varphi \ast \vu)(\vec{\nabla}\otimes \vu_\alpha)^T \Big)+\vu_\alpha=\vu}$. Integrating by parts this equation, for fixed $ 0\leq t \leq T$, and for all $\vec{\phi} \in {H}^1_{\text{s}}(\T)$ (recall that $\text{div}(\vec{\phi})=0$ and that $\P$ is an auto-adjoint operator) we obtain the following variational formulation: 
\begin{equation}\label{Variational}
\begin{split}
&\alpha^2 \int_{\T} A((\varphi\ast \vu)(t,x)) (\vec{\nabla} \otimes \vu_\alpha)^T (t,x): (\vec{\nabla} \otimes \vec{\phi})^T (x) dx + \int_{\T} \vu_\alpha(t,x) \cdot \vec{\phi} (x) dx \\
=& \int_{\T} \vu(t,x)\cdot \vec{\phi}(x) dx. 
\end{split}
\end{equation}
As  the term $A((\varphi\ast \vu))$ does not depend on $\vu_\alpha$ nor $\vec{\phi}$, the first and the second expression on the left-hand side in (\ref{Variational}) defines bilinear form
\[ a\left(\vu_\alpha\,,\vec{\phi}\right):= \alpha^2 \int_{\T} A((\varphi\ast \vu)(t,x)) (\vec{\nabla} \otimes \vu_\alpha)^T (t,x): (\vec{\nabla} \otimes \vec{\phi})^T (x) dx+\int_{\T} \vu_\alpha(t,x) \cdot \vec{\phi} (x) dx.\]
Lower and upper bounds on the function $A(\cdot)$ given in (\ref{Conditions-a-existence-1}) yield that $a(\cdot , \cdot)$ is 
 continuous and coercive  in the Banach space ${H}^1_{\text{s}}(\T)$.
By a direct application of the Lax-Milgram theorem, for any fixed time $ 0\leq t \leq T$, there exists $\vu_\alpha=\vu_\alpha(t,\cdot) \in {H}^1_{\text{s}}(\T)$ a  solution to equation (\ref{Variational}).  Recall that by definition of space ${H}^1_{\text{s}}(\T)$ we directly have $\vu_\alpha(t,\cdot)\in {H}^1(\T)$ and $\text{div}(\vu_\alpha(t,\cdot))=0$.

\medskip

Moreover, setting $\vec{\phi}=\vu_\alpha(t,\cdot)$ in (\ref{Variational}), and using again the lower bound in (\ref{Conditions-a-existence-1}), we have 
\begin{equation}\label{Control-Filter-L2}
    \alpha^2 \beta \| \vec{\nabla}\otimes \vu_\alpha (t,\cdot)\|^{2}_{L^2}+\| \vu_\alpha(t,\cdot)\|^{2}_{L^2}\leq \int_{\T} \vu(t,x)\cdot \vu_\alpha(t,x) dx. 
\end{equation} 
Applying the Cauchy-Schwarz  inequality in the last term, we get 
\begin{equation*}
  \int_{\T} \vu(t,\cdot)\cdot \vu_\alpha(t,\cdot) dx \leq \| \vu(t,\cdot)\|_{L^2} \| \vu_\alpha(t,\cdot)\|_{L^2} \leq \frac{1}{2}\| \vu(t,\cdot)\|^{2}_{L^2}+\frac{1}{2}\|\vu_\alpha(t,\cdot)\|^{2}_{L^2}.  
\end{equation*}
Then, we can write
\begin{equation*}
\min(\alpha^2 \beta, 1/2)\| \vu_\alpha(t,\cdot)\|^{2}_{H^1} \leq    \alpha^2 \beta \| \vec{\nabla}\otimes \vu_\alpha (t,\cdot)\|^{2}_{L^2}+ \frac{1}{2}\| \vu_\alpha(t,\cdot)\|^{2}_{L^2}\leq  \frac{1}{2}\| \vu(t,\cdot)\|^{2}_{L^2},  \end{equation*}
which yields the wished estimate (\ref{Control-Filter}) with $\ds{K_0=\frac{1}{2\min(\alpha^2 \beta, 1/2)}\simeq \frac{1}{\min(\alpha^2 \beta,1)}}$, and we thus have   $\vu_\alpha \in L^\infty([0,T],H^1(\T))$.

\medskip

Finally,  we use classical arguments to recover the pressure term $P_\alpha$. Since $\vu_\alpha$ solves the   equation $\ds{ -\alpha^2 \P\, \text{div}\Big( A(\varphi \ast \vu)(\vec{\nabla}\otimes \vu_\alpha)^T \Big)+\vu_\alpha=\vu}$,  and $\text{div}(\vu_\alpha) = \text{div}(\vu) = 0$, we have 
\[ \P \left (  -\alpha^2  \text{div}\Big( A(\varphi \ast \vu)(\vec{\nabla}\otimes \vu_\alpha)^T \Big)+\vu_\alpha - \vu\right)=0.\]
By well-known properties of the Leray's projector $\P$ (see \cite[Lemma $6.3$]{PLe}) there exists a distribution $P_\alpha$ such that 
\[-\alpha^2  \text{div}\Big( A(\varphi \ast \vu)(\vec{\nabla}\otimes \vu_\alpha)^T \Big)+\vu_\alpha=\vu -\vec{\nabla}P_\alpha.\]
Moreover, the pressure $P_\alpha$ is related to $\vu$ and $\vu_\alpha$ by the second expression in (\ref{Characterization-pressions}). Since $\vu_\alpha \in L^\infty_t H^1_x$ and,  by (\ref{Conditions-a-existence-1}),   $A(\varphi \ast \vu)\in L^\infty_t L^\infty_x$,  it is straightforward  to verify that the term on the right-hand side belongs to the space $L^\infty([0, T], L^2(\T))$. Consequently,  $P_\alpha \in L^\infty([0, T], L^2(\T))$.  Proposition \ref{Existence-Filter} is proven.  \finpv

\begin{Remarque}\label{Rmk0} The proof follows a standard variational argument, where  the convolution term $\varphi$ in equation (\ref{Filter-Regularized}) does not play any substantial role. Consequently, this result also  holds for the non-regularized nonlinear filter equation:
\[   -\alpha^2 \text{div} \left( A (\vu) (\vec{\nabla} \otimes \vu_\alpha)^T \right)+\vu_\alpha+\vec{\nabla}P_\alpha = \vu, \quad  \text{div}(\vu_\alpha)=0. \]
\end{Remarque}

Conversely, additional information on $\varphi$ is required to prove the following: 
\begin{Proposition}[$H^1$-continuous dependence of data]\label{Dependence-data-Filter} Under the same hypothesis of  Proposition \ref{Existence-Filter}, for $i=1,2$, denote $\vu_{i,\alpha}\in L^\infty([0,T], H^1(\T))$ the corresponding solutions  of equation (\ref{Filter-Regularized}) arising from the data $\vu_i \in L^\infty([0,T], L^2_{\text{s}}(\T))$. 

\medskip

Additionally, assume that the indicator function $A(\cdot)$ verifies   (\ref{Conditions-a-existence-2}) and the convolution term $\varphi$ verifies (\ref{Conditions-varphi-existence}). Moreover,  define the quantity
\begin{equation}\label{L1}
    L_1= \alpha^2 C_A\, \| \varphi \|_{L^2},
\end{equation}
where $C_A>0$ is the  Lipschitz constant of $A(\cdot)$. Then, for any time $0\leq t \leq T$, it holds: 
\begin{equation}\label{Dependence-data}
\| \vu_{1,\alpha}(t,\cdot) - \vu_{2,\alpha}(t,\cdot)\|_{H^1} \leq K_1\Big(\| \vu_2(t,\cdot)\|_{L^2}+1\Big) \| \vu_1(t,\cdot)-\vu_2(t,\cdot)\|_{L^2}, \quad \mbox{with} \ \  K_1\simeq \max\left( K^{3/2}_{0}L1, K_0 \right),
\end{equation}  
where the quantity $K_0$ is given in (\ref{Control-Filter}).
\end{Proposition} 
\pv In order to simplify our writing, we shall omit the time dependence of functions. Note that $\vu_{\alpha,1}-\vu_{\alpha,2}$ solves the equation
\begin{equation*}
-\alpha^2 \P\, \text{div}\left(A(\varphi \ast \vu_1)(\vec{\nabla}\otimes \vu_{1,\alpha})^T -A(\varphi \ast \vu_2)(\vec{\nabla}\otimes \vu_{2,\alpha})^T\right)+(\vu_{1,\alpha}-\vu_{2,\alpha})= \vu_1 - \vu_2,
\end{equation*}
hence, we write
\begin{equation}\label{PDE-Filter-difference}
\begin{split}
-\alpha^2 \P\,\text{div}\Big(A(\varphi \ast \vu_1) &(\vec{\nabla}\otimes (\vu_{1,\alpha}-\vu_{2,\alpha}))^T +\big[A(\varphi\ast \vu_1)-A(\varphi \ast \vu_2)\big](\vec{\nabla}\otimes \vu_{2,\alpha})^T\Big)\\
&+(\vu_{1,\alpha}-\vu_{2,\alpha})= \vu_1 - \vu_2.
\end{split}
\end{equation}
Multiplying this equation by $(\vu_{1,\alpha}-\vu_{2,\alpha})$, and integrating by parts, we obtain
\begin{equation*}
\begin{split}
&\,\alpha^2 \int_{\T} A(\varphi \ast \vu_1) | \vec{\nabla}\otimes (\vu_{1,\alpha}-\vu_{2,\alpha})|^2 dx \\
&+\,  \alpha^2 \int_{\T} \Big(A(\varphi\ast \vu_1)-A(\varphi \ast \vu_2)\Big)(\vec{\nabla}\otimes \vu_{2,\alpha})^T :  (\vec{\nabla}\otimes (\vu_{1,\alpha}-\vu_{2,\alpha}))^T \, dx + \| \vu_{1,\alpha}-\vu_{2,\alpha} \|^{2}_{L^2}\\
=&\, \int_{\T} (\vu_1 - \vu_2)\cdot (\vu_{1,\alpha}-\vu_{2,\alpha}) dx.
\end{split}
\end{equation*}
By the lower bound in (\ref{Conditions-a-existence-1}),   setting $C_{\alpha,\beta}=\min(\alpha^2 \beta,1)$, we have
\begin{equation}\label{Estim-Filter-01}
\begin{split}
C_{\alpha ,\beta}\, \|\vu_{1,\alpha}-\vu_{2,\alpha}\|^{2}_{H^1} \leq &\,  \,\alpha^2 \beta \|\vec{\nabla}\otimes (\vu_{1,\alpha}-\vu_{2,\alpha})\|^{2}_{L^2}+ \| \vu_{1,\alpha}-\vu_{2,\alpha} \|^{2}_{L^2} \\
\leq &\, \alpha^2 \int_{\T} \big|A(\varphi\ast \vu_1)-A(\varphi \ast \vu_2)\big|\, \big|\vec{\nabla}\otimes \vu_{2,\alpha}\big| \, \big|\vec{\nabla}\otimes (\vu_{1,\alpha}-\vu_{2,\alpha})\big| \, dx \\
&\, + \int_{\T} |\vu_1 - \vu_2|\,|\vu_{1,\alpha}-\vu_{2,\alpha}| dx,
\end{split}
\end{equation}
where we must estimate each term on the right-hand side. For the first term, recall that by (\ref{Conditions-a-existence-2}) the function $A(\cdot)$ is Lipschitz-continuous, and for the constant $C_A>0$  we can write
\begin{equation*}
\begin{split}
 &\, \alpha^2 \int_{\T} \big|A(\varphi\ast \vu_1)-A(\varphi \ast \vu_2)\big|\, \big|\vec{\nabla}\otimes \vu_{2,\alpha}\big|\,   \big|\vec{\nabla}\otimes (\vu_{1,\alpha}-\vu_{2,\alpha})\big| \, dx \\
\leq &\, \alpha^2 C_A\, \int_{\T} \big| \varphi\ast ( \vu_1-\vu_2)\big|\, \big|\vec{\nabla}\otimes \vu_{2,\alpha}\big| \, \big|\vec{\nabla}\otimes (\vu_{1,\alpha}-\vu_{2,\alpha})\big| \, dx.
 \end{split}    
\end{equation*}
\begin{Remarque}\label{Rmk1} Note that with the only available information that $\vu_{i,\alpha} \in H^1(\T)$, the second and third terms mentioned above are only controlled in the $L^2$-norm. This necessitates controlling  the first term in the $L^\infty$-norm. However,  considering only the information that $\vu_i \in L^2_{\text{s}}(\T)$,  the convolution term $\varphi$ assumes significant importance.
\end{Remarque}
By (\ref{Conditions-varphi-existence})  we  have $\varphi \in L^{2}(\T)$. We thus  apply H\"older inequalities and Young inequalities to get
\begin{equation}\label{Estim-Filter-02}
    \begin{split}
&\,  \alpha^2 C_A\, \int_{\T} \big| \varphi\ast ( \vu_1-\vu_2)\big|\,\big|\vec{\nabla}\otimes \vu_{2,\alpha}\big|\, \big|\vec{\nabla}\otimes (\vu_{1,\alpha}-\vu_{2,\alpha})\big| \, dx \\
\leq &\, \alpha^2 C_A\, \| \varphi\ast ( \vu_1-\vu_2)\|_{L^\infty}\| \vec{\nabla}\otimes \vu_{2,\alpha} \|_{L^2}\,\| \vec{\nabla}\otimes (\vu_{1,\alpha}-\vu_{2,\alpha})\|_{L^2}\\
\leq &\, \alpha^2 C_A\, \| \varphi \|_{L^2}\| \vu_1-\vu_2 \|_{L^2}\| \vec{\nabla}\otimes \vu_{2,\alpha} \|_{L^2}\,\| \vec{\nabla}\otimes (\vu_{1,\alpha}-\vu_{2,\alpha})\|_{L^2}\\
\leq &\, \alpha^2 C_A\, \| \varphi \|_{L^2}\| \vu_1-\vu_2 \|_{L^2}\| \vec{\nabla}\otimes \vu_{2,\alpha} \|_{L^2}\,\| \vu_{1,\alpha}-\vu_{2,\alpha}\|_{H^1}.
    \end{split}
\end{equation}
For the second term, using  the Cauchy-Schwarz  inequality, we write
\begin{equation}\label{Estim-Filter-03}
 \int_{\T} |\vu_1 - \vu_2|\,|\vu_{1,\alpha}-\vu_{2,\alpha}| dx\leq \| \vu_1 - \vu_2 \|_{L^2}\|\vu_{1,\alpha}-\vu_{2,\alpha}\|_{L^2}\leq \| \vu_1 - \vu_2 \|_{L^2}\|\vu_{1,\alpha}-\vu_{2,\alpha}\|_{H^1}. \end{equation}

With estimates (\ref{Estim-Filter-02}) and (\ref{Estim-Filter-03}) at hand, we get back to the inequality (\ref{Estim-Filter-01}) to obtain
\begin{equation*}
C_{\alpha,\beta}\|\vu_{1,\alpha}-\vu_{2,\alpha}\|^{2}_{H^1} \leq  \Big( \alpha^2 C_A\, \| \varphi \|_{L^2}\|\vec{\nabla}\otimes \vu_{2,\alpha}\|_{L^2}+1\Big)\|\vu_1-\vu_2 \|_{L^2}\,\| \vu_{1,\alpha}-\vu_{2,\alpha}\|_{H^1},   
\end{equation*}
hence, we get
\begin{equation*}
C_{\alpha,\beta}\|\vu_{1,\alpha}-\vu_{2,\alpha}\|_{H^1} \leq  \Big( \alpha^2 C_A\, \| \varphi \|_{L^2}\|\vec{\nabla}\otimes \vu_{2,\alpha}\|_{L^2}+1\Big)\|\vu_1-\vu_2 \|_{L^2}.   
\end{equation*}
Moreover, remark that by estimate (\ref{Control-Filter}) we have $\ds{\|\vec{\nabla}\otimes \vu_{2,\alpha}\|_{L^2} \leq \sqrt{K_0} \| \vu_2\|_{L^2}}$.  Then, we write
\begin{equation*}
    \begin{split}
    \|\vu_{1,\alpha}-\vu_{2,\alpha}\|_{H^1} \leq &\, \frac{1}{C_{\alpha,\beta}} \Big( \alpha^2 C_A\, \| \varphi \|_{L^2}\sqrt{C_0}\| \vu_2\|_{L^2}+1\Big)\|\vu_1-\vu_2 \|_{L^2}\\
    \leq &\, \frac{1}{C_{\alpha,\beta}} \max\left(\alpha^2 C_A\, \| \varphi \|_{L^2}\sqrt{K_0},1\right) \Big(\| \vu_2\|_{L^2}+1\Big)\|\vu_1-\vu_2 \|_{L^2}\\
  =: &\, K_1 \Big(\| \vu_2\|_{L^2}+1\Big)\|\vu_1-\vu_2 \|_{L^2}.
    \end{split} 
\end{equation*} 
Finally, recalling that $C_{\alpha,\beta}=\min(\alpha^2 \beta,1)$,  $K_0\simeq \frac{1}{\min(\alpha^2 \beta,1)}$ and $L_1=:\alpha^2 C_A \|\varphi \|_{L^2}$, we have  
\[ K_1 =: \frac{1}{C_{\alpha,\beta}} \max\left(\alpha^2 C_A\, \| \varphi \|_{L^2}\sqrt{K_0},1\right)   \simeq \max\left( K^{3/2}_{0}L1, K_0 \right). \]
Proposition \ref{Dependence-data-Filter} is proven. \finpv

\medskip

\begin{Proposition}[$H^2$-regularity]\label{Reg-H2-Filter} Within the setting of   Proposition \ref{Existence-Filter}, additionally  assume that the indicator function $A(\cdot)$ and the convolution term $\varphi$ satisfy (\ref{conditions-a-varphi-uniqueness}), and define the quantity 
\begin{equation}\label{L2}
    L_2= \alpha^2 \| \vec{\nabla}A(\cdot)\|_{L^\infty} \| \varphi \|_{H^1}.
\end{equation}
Then, the solution $\vu_\alpha \in L^\infty([0,T],H^1(\T))$ to equation (\ref{Filter-Regularized}) verifies $\vu_\alpha \in L^\infty([0,T], H^2(\T))$, and for any time $0\leq t \leq T$, the following estimate holds:
\begin{equation}\label{Filter-Control-H2}
    \| \vu_\alpha(t,\cdot)\|_{H^2} \leq K_2 \big(\| \vu(t,\cdot)\|_{L^2}+1\big)\|\vu(t,\cdot)\|_{L^2}, \quad \mbox{with} \ \ K_2\simeq \max\left( K^{3/2}_0 L_2, K_0\right), 
\end{equation}
where $K_0$ is  defined in (\ref{Control-Filter}).
\end{Proposition}
\pv We shall use the method of differential quotients. Let $h\in \Rt$ be such that $h\neq 0$.  For every $\vec{\phi}\in \mathcal{V}$ define the operator 
\[ D_h \vec{\phi}(x)= \frac{\vec{\phi}(x+h)-\vec{\phi}(x)}{|h|},\]
which verifies the  useful properties:
\begin{equation}\label{Iden-Diff-quocient}
\int_{\T} \vec{\phi}_1 \cdot D_h \vec{\phi}_2\, dx = \int_{\T} D_{-h} \vec{\phi}_1 \cdot \vec{\phi}_2 \, dx, \quad D_h(\vec{\phi}_1 \cdot \vec{\phi}_2)= (D_h \vec{\phi}_1)\cdot \vec{\phi}_2 + \vec{\phi}_1\cdot (D_h \vec{\phi}_2).
\end{equation}

As before, we omit the time dependence of functions. Coming back to the variational formulation (\ref{Variational}), and setting  $\vec{\phi}=D_{-h}D_h \vu_\alpha \in {H}^1_{\text{s}}(\T)$, we obtain
\begin{equation*}
\begin{split}
&\alpha^2 \int_{\T} A(\varphi\ast \vu) (\vec{\nabla} \otimes \vu_\alpha)^T : (\vec{\nabla} \otimes D_{-h}D_h \vu_\alpha)^T  dx + \int_{\T} \vu_\alpha \cdot D_{-h}D_h \vu_\alpha dx \\
=& \int_{\T} \vu \cdot D_{-h}D_h \vu_\alpha dx.
\end{split}
\end{equation*}
In the first and the second expressions on the left,  we apply the identities in   (\ref{Iden-Diff-quocient}) to write
\begin{equation*}
\begin{split}
   &\,  \alpha^2 \int_{\T} A(\varphi\ast \vu) (\vec{\nabla} \otimes \vu_\alpha)^T : (\vec{\nabla} \otimes D_{-h}D_h \vu_\alpha)^T  dx  =  \alpha^2 \int_{\T} A(\varphi\ast \vu) | \vec{\nabla}\otimes D_h \vu_\alpha |^2 dx \\
    &\, + \alpha^2 \int_{\T} D_h (A(\varphi\ast \vu)) (\vec{\nabla} \otimes \vu_\alpha)^T:(\vec{\nabla} \otimes D_h \vu_\alpha)^T dx,  
 \end{split}
\end{equation*}
and
\begin{equation*}
  \int_{\T} \vu_\alpha \cdot D_{-h}D_h \vu_\alpha dx=  \int_{\T}| D_h \vu_\alpha|^2 dx. 
\end{equation*}
Then, we rearrange  terms to write
\begin{equation*}
\begin{split}
&\alpha^2 \int_{\T} A(\varphi\ast \vu) | \vec{\nabla}\otimes D_h \vu_\alpha |^2 dx + \int_{\T}| D_h \vu_\alpha|^2 dx\\
=&\, - \alpha^2 \int_{\T} D_h (A(\varphi\ast \vu)) (\vec{\nabla} \otimes \vu_\alpha)^T:(\vec{\nabla} \otimes D_h \vu_\alpha)^T dx+ \int_{\T} \vu \cdot D_{-h}D_h \vu_\alpha dx.
\end{split}
\end{equation*}
By the lower bound in (\ref{Conditions-a-existence-1}) and using the constant $C_{\alpha,\beta}=\min(\alpha^2\beta,1)>0$, we obtain 
\begin{equation}\label{Estim-Filter-03}
    C_{\alpha,\beta} \| D_h \vu_\alpha \|^{2}_{H^1} \leq \,  \alpha^2 \int_{\T} \left| D_h (A(\varphi\ast \vu)) (\vec{\nabla} \otimes \vu_\alpha)^T:(\vec{\nabla} \otimes D_h \vu_\alpha)^T  \right| dx + \int_{\T} | \vg  \cdot D_{-h}D_h \vu_\alpha| dx, 
\end{equation}
where we must estimate each term on the right-hand side. For the first term, note that  by \cite[Proposition $IX.3$]{Brezis} and the Young inequalities (with $1+1/\infty=1/2+1/2$), we have 
\begin{equation*}
 \| D_h (A(\varphi\ast \vu)) \|_{L^\infty} \leq \,  C \| \vec{\nabla}\big(A(\varphi\ast \vu)\big)\|_{L^\infty} \leq C \sum_{i=1}^{3} \| \partial_i A (\varphi\ast \vu) (\partial_i \varphi) \ast \vu \|_{L^\infty}\leq C \| \vec{\nabla}A(\cdot)\|_{L^\infty} \| \varphi \|_{H^1}\| \vu \|_{L^2}.      
\end{equation*}
\begin{Remarque}\label{Rmk2} Note that the information  $\vu \in L^2(\T)$ is not sufficient to ensure that $\vec{\nabla}(A(\vu))\in L^\infty(\T)$, even  if the indicator function  $A(\cdot)$  is sufficiently regular. 
\end{Remarque}
Continuing with the proof, by H\"older inequalities and by the estimate above and the control (\ref{Control-Filter}), we write
\begin{equation}\label{Estim-Filter-04}
 \begin{split}
  &\, \alpha^2 \int_{\T} \left| D_h (A(\varphi\ast \vu)) (\vec{\nabla} \otimes \vu_\alpha)^T:(\vec{\nabla} \otimes D_h \vu_\alpha)^T  \right| dx \\
 \leq &\,  \alpha^2 \| D_h (A(\varphi\ast \vu)) \|_{L^\infty} \| \vec{\nabla} \otimes \vu_\alpha \|_{L^2} \, \| \vec{\nabla} \otimes D_h \vu_\alpha \|_{L^2}\\
 \leq &\, \Big( \alpha^2 C \|\vec{\nabla}A(\cdot)\|_{L^\infty} \| \varphi \|_{H^1}\| \vu \|_{L^2} \Big) \, \Big(\sqrt{K_0} \| \vu \|_{L^2} \Big) \, \| \vec{\nabla} \otimes D_h \vu_\alpha \|_{L^2}\\
 \leq &\, \Big( \alpha^2 C \sqrt{K_0}  \|\vec{\nabla}A(\cdot)\|_{L^\infty} \| \varphi \|_{H^1}\Big) \, \| \vu \|^{2}_{L^2} \| D_h \vu_\alpha \|_{H^1}.
 \end{split}   
\end{equation}
For the second term,  applying  the Cauchy-Schwarz inequality and \cite[Proposition $IX.3$]{Brezis},  we write
\begin{equation}\label{Estim-Filter-05}
\int_{\T} | \vu  \cdot D_{-h}D_h \vu_\alpha| dx \leq C \| \vu \|_{L^2} \| D_{-h}D_h \vu_\alpha \|_{L^2} \leq C  \| \vg \|_{L^2} \| D_h \vu_\alpha \|_{L^2}  \leq C\| \vu \|_{L^2} \|D_h \vu_\alpha \|_{H^1}. 
\end{equation}
With estimates (\ref{Estim-Filter-04}) and (\ref{Estim-Filter-05}) at hand, we come back to estimate (\ref{Estim-Filter-03}), to get 
\begin{equation*}
 C_{\alpha,\beta} \| D_h \vu_\alpha \|^{2}_{H^1}   \leq C\left( \alpha^2 C \sqrt{K_0} \|\vec{\nabla}A(\cdot)\|_{L^\infty} \| \varphi \|_{H^1}\| \vu \|_{L^2}+ 1 \right) \| \vu \|_{L^2} \|D_h \vu_\alpha \|_{H^1},
\end{equation*} 
hence, 
\begin{equation*}
 C_{\alpha,\beta} \| D_h \vu_\alpha \|_{H^1}   \leq C\left( \alpha^2 C\sqrt{K_0}  \|\vec{\nabla}A(\cdot)\|_{L^\infty} \| \varphi \|_{H^1}\| \vu \|_{L^2}+ 1 \right)\| \vu \|_{L^2}.
\end{equation*}
With this estimate and letting $h\to 0$, by  \cite[Proposition $IX.3$]{Brezis} we obtain that  $\vu_\alpha \in H^2(\T)$. In addition, we following estimates holds:  
\begin{equation}\label{Control-Filter-H2}
\begin{split}
\| \vu_\alpha \|_{H^2} \leq &\,  \frac{C}{C_{\alpha,\beta}}\left( \alpha^2 C\sqrt{C_0}\|\vec{\nabla}A(\cdot)\|_{L^\infty} \| \varphi \|_{H^1}\| \vu \|_{L^2}+ 1 \right)\| \vu \|_{L^2}\\
\leq &\, \frac{1}{C_{\alpha,\beta}} \max\left( \alpha^2 C\sqrt{C_0}\|\vec{\nabla}A(\cdot)\|_{L^\infty} \| \varphi \|_{H^1} , 1 \right) \left( \|  \vu \|_{L^2}+ 1 \right)\| \vu \|_{L^2} \\
=:&\, K_2 \left( \|  \vu \|_{L^2}+ 1 \right)\| \vu \|_{L^2}.
\end{split}
\end{equation}
Finally, recalling that $C_{\alpha,\beta}=\min(\alpha^2 \beta,1)$,  $K_0 \simeq \frac{1}{\min(\alpha^2 \beta,1)}$ and $L_2:= \alpha^2 \| \vec{\nabla}A(\cdot)\|_{L^\infty} \| \varphi \|_{H^1}$, we have
\[ K_2:= \frac{1}{C_{\alpha,\beta}} \max\left( \alpha^2 C\sqrt{C_0}\|\vec{\nabla}A(\cdot)\|_{L^\infty} \| \varphi \|_{H^1} , 1 \right)  \simeq \max\left( K^{3/2}_{0} L_2, K_0 \right). \]
Proposition \ref{Reg-H2-Filter} is proven.  \finpv

\begin{Proposition}[$H^2$-continuous dependence of data]\label{Cont-Dep-H2-Filter} Under the same hypothesis   of Proposition \ref{Reg-H2-Filter}, define the quantity
\begin{equation}\label{L3}
    L_3= \alpha^2 \| \varphi \|_{H^1}\max\left( C^{'}_{A} \| \varphi\|_{H^1}, \| \vec{\nabla}A(\cdot)\|_{L^\infty}\right),
\end{equation}
where $C^{'}_{A}$ is the Lipschitz constant of $\vec{\nabla}A(\cdot)$. 

\medskip

For $i=1,2$, let $\vu_{i,\alpha} \in L^\infty([0,T],H^2(\T)) $ be two solutions of equation (\ref{Filter-Regularized}) arising from the data $\vu_i \in L^\infty([0,T],L^2_{\text{s}}(\T))$. Then, for any time $0\leq t \leq T$, the following estimate holds:
\begin{equation}\label{Control-Filter-Diff-H2}
 \| \vu_{1,\alpha}(t,\cdot)-\vu_{2,\alpha}(t,\cdot)\|_{H^2}\leq K_3\Big(\| \vu_1(t,\cdot) \|_{L^2}+\| \vu_2(t,\cdot) \|_{L^2}+1\Big)^2  \| \vu_1(t,\cdot)-\vu_2(t,\cdot) \|_{L^2},
\end{equation}
where
\begin{equation}\label{K3}
K_3 \simeq \max\left( K_0K_2L_1, K_0 K_1 L_2, K^{3/2}_{0}L_3\right),
\end{equation}
with $K_0$ given in (\ref{Control-Filter}), $K_1$ given in (\ref{Dependence-data}), $K_2$ defined in (\ref{Filter-Control-H2}), $L_1$ defined in (\ref{L1}) and $L_2$ given in (\ref{L2}).
\end{Proposition}
\pv For simplicity, we shall denote $\vw_\alpha = \vu_{1,\alpha}- \vu_{2,\alpha}$ and $\vw=\vu_1-\vu_2$. Remark that  $\vw_\alpha$ verifies the equation (\ref{PDE-Filter-difference}), which  writes down as follows:
\begin{equation*}
-\alpha^2 \P\, \text{div}\Big(A(\varphi \ast \vu_1) (\vec{\nabla}\otimes \vw_\alpha  )^T +\big(A(\varphi\ast \vu_1)-A(\varphi \ast \vu_2)\big)(\vec{\nabla}\otimes \vu_{2,\alpha})^T\Big) + \vw_\alpha= \vw.
\end{equation*}

As  the indicator function $A(\cdot)$ and the convolution term $\varphi$ satisfy (\ref{conditions-a-varphi-uniqueness}), and as  $\vu_1, \vu_2 \in L^2_{\text{s}}(\T)$, by Proposition \ref{Reg-H2-Filter} we have  $\vu_{1,\alpha}, \vu_{2,\alpha}, \vw_\alpha \in H^2(\T)$. Then,  for fixed $k=1,2,3$,  we can  multiply the equation above  by $\partial^2_k \vw_\alpha$, and  integrating by parts each term (first moving $\P$ and the divergence operator to the right and thereafter moving the derivative $\partial_k$ to the left) we find:
\begin{equation*}
\begin{split}
&\,-\alpha^2 \int_{\T} \partial_k \left(A(\varphi \ast \vu_1)(\vec{\nabla} \otimes \vw_\alpha)^T\right): (\vec{\nabla} \otimes \partial_k \vw_\alpha)^T \, dx \\
&\, - \alpha^2 \int_{\T} \partial_k\left(\big(A(\varphi\ast \vu_1)-A(\varphi \ast \vu_2)\big)(\vec{\nabla} \otimes \vu_{2,\alpha})^T \right)  : (\vec{\nabla} \otimes \partial_k \vw_{\alpha})^T \, dx \\
&\,- \| \partial_k \vw_\alpha\|^2_{L^2}= \int_{\T} \vw \cdot \partial^2_k \vw_\alpha\, dx.
\end{split}
\end{equation*}
We develop the terms involving the derivative of a product, and  we  multiply by $-1$ each side of this identity  to obtain 
\begin{equation*}
\begin{split}
&\,\alpha^2 \int_{\T} \partial_k (A(\varphi \ast \vu_1)) (\vec{\nabla} \otimes \vw_\alpha)^T : (\vec{\nabla} \otimes \partial_k \vw_\alpha)^T\, dx\\
+&\, \alpha^2 \int_{\T}  A(\varphi \ast \vu_1) | \vec{\nabla} \otimes \partial_k \vw_\alpha |^2 \, dx\\
+&\, \alpha^2 \int_{\T} \partial_k\big(A(\varphi\ast \vu_1)-A(\varphi \ast \vu_2)\big) (\vec{\nabla} \otimes \vu_{2,\alpha})^T : (\vec{\nabla} \otimes \partial_k \vw_{\alpha})^T\, dx \\
+&\, \alpha^2 \int_{\T} \big(A(\varphi\ast \vu_1)-A(\varphi \ast \vu_2)\big) (\vec{\nabla}\otimes \partial_k \vu_{2,\alpha})^T : (\vec{\nabla}\otimes \partial_k \vw_{\alpha})^T \, dx\\
+&\, \| \partial_k \vw_\alpha\|^2_{L^2}= -\int_{\T} \vw \cdot \partial^2_k \vw_\alpha \, dx.
\end{split}
\end{equation*}
According to the lower bound in (\ref{Conditions-a-existence-1}), the second term above is estimated from below as:
\[ \alpha^2 \beta \| \vec{\nabla}\otimes \partial_k \vw_\alpha \|^{2}_{L^2}\leq  \alpha^2 \int_{\T}  A(\varphi \ast \vu_1) | \vec{\nabla} \otimes \partial_k \vw_\alpha |^2 \, dx.\]
Moreover, always  with the constant $C_{\alpha,\beta}=\min(\alpha^2\beta,1)>0$,  we write
\[ C_{\alpha,\beta} \| \partial_{k} \vw_\alpha \|^{2}_{H^1} \leq \alpha^2 \beta  \| \vec{\nabla}\otimes \partial_k \vw_\alpha \|^{2}_{L^2} +  \| \partial_k \vw_\alpha \|^{2}_{L^2}.\]
Rearranging  terms, we obtain
\begin{equation}\label{Estim-Cont-Dep-01}
\begin{split}
 \,C_{\alpha,\beta} \| \partial_{k} \vw_\alpha \|^{2}_{H^1} 
 \leq &\, \alpha^2 \int_{\T} \left|\partial_k (A(\varphi \ast \vu_1)) (\vec{\nabla} \otimes \vw_\alpha)^T : (\vec{\nabla} \otimes \partial_k \vw_\alpha)^T \right|\, dx\\
 &\, + \alpha^2 \int_{\T} \left|\partial_k\big(A(\varphi\ast \vu_1)-A(\varphi \ast \vu_2)\big) (\vec{\nabla} \otimes \vu_{2,\alpha})^T : (\vec{\nabla} \otimes \partial_k \vw_{\alpha})^T \right|\, dx \\
 &\,+ \alpha^2 \int_{\T} \left|\big(A(\varphi\ast \vu_1)-A(\varphi \ast \vu_2)\big) (\vec{\nabla}\otimes \partial_k \vu_{2,\alpha})^T : (\vec{\nabla}\otimes \partial_k \vw_{\alpha})^T \right| \, dx\\
 &\, + \int_{\T} \left| \vw \cdot \partial^2_k \vw_\alpha \, \right| dx.
 \end{split}
\end{equation}
Now, we need to estimate each term on the right-hand side. For the reader's convenience, we will provide some estimates in the forthcoming technical lemmas.
\begin{Lemme} For the first term, we have
\begin{equation}\label{Estim-Cont-Dep-02}
\begin{split}
    &\, \alpha^2 \int_{\T} \left|\partial_k (A(\varphi \ast \vu_1)) (\vec{\nabla} \otimes \vw_\alpha)^T : (\vec{\nabla} \otimes \partial_k \vw_\alpha)^T \right|\, dx \\
    \leq &\, K_1 L_2    \| \vu_1 \|_{L^2} \big(\| \vu_2\|_{L^2}+1\big)\|\vw \|_{L^2}\, \|  \partial_k \vw_\alpha\|_{H^1}.
 \end{split}
\end{equation}
\end{Lemme}
\pv As $A(\cdot)$ and $\varphi$ verify (\ref{conditions-a-varphi-uniqueness}), using H\"older and Young inequalities we can write
\begin{equation*}
\begin{split}
 &\,  \alpha^2 \int_{\T} \left|\partial_k (A(\varphi \ast \vu_1)) (\vec{\nabla} \otimes \vw_\alpha)^T : (\vec{\nabla} \otimes \partial_k \vw_\alpha)^T \right|\, dx\\
 \leq &\, \alpha^2 \| \partial_k (A(\varphi \ast \vu_1)) (\vec{\nabla} \otimes \vw_\alpha) \|_{L^2}\, \| \vec{\nabla} \otimes \partial_k \vw_\alpha\|_{L^2} \\
 \leq &\, \alpha^2 \| \partial_k (A(\varphi \ast \vu_1))  \|_{L^\infty}\, \|  \vec{\nabla} \otimes \vw_\alpha \|_{L^2}\, \| \vec{\nabla} \otimes \partial_k \vw_\alpha\|_{L^2}\\
 = &\, \alpha^2  \| \partial_k A(\varphi \ast \vu_1) (\partial_k \varphi) \ast \vu_1\|_{L^\infty}\, \|  \vec{\nabla} \otimes \vw_\alpha \|_{L^2}\, \| \vec{\nabla} \otimes \partial_k \vw_\alpha\|_{L^2}\\
 \leq & \alpha^2  \| \vec{\nabla}A(\cdot)\|_{L^\infty}\, \| (\partial_k\varphi) \ast \vu_1 \|_{L^\infty}\, \|  \vec{\nabla} \otimes \vw_\alpha \|_{L^2}\, \| \vec{\nabla} \otimes \partial_k \vw_\alpha\|_{L^2}\\
 \leq &\, \alpha^2 \| \vec{\nabla}A(\cdot)\|_{L^\infty} \| \partial_k\varphi \|_{L^2}\, \| \vu_1 \|_{L^2}\, \|  \vec{\nabla} \otimes \vw_\alpha \|_{L^2}\, \| \vec{\nabla} \otimes \partial_k \vw_\alpha\|_{L^2}\\
  \leq &\, \alpha^2 \| \vec{\nabla}A(\cdot)\|_{L^\infty} \| \varphi \|_{H^1}\, \| \vu_1 \|_{L^2}\, \| \vw_\alpha \|_{H^1}\, \|  \partial_k \vw_\alpha\|_{H^1}. 
\end{split}
\end{equation*}
Moreover, recalling that $\vw_\alpha=\vu_{\alpha,1}-\vu_{\alpha,2}$ and $\vw=\vu_1 - \vu_2$, by Proposition {\ref{Dependence-data-Filter}} we have the control
\[ \| \vw_\alpha \|_{H^1} \leq K_1(\| \vu_2\|_{L^2}+1)\|\vw \|_{L^2}.\]
With these estimates we obtain
\begin{equation*}
\begin{split}
&\, \alpha^2 \int_{\T} \left|\partial_k (A(\varphi \ast \vu_1)) (\vec{\nabla} \otimes \vw_\alpha)^T : (\vec{\nabla} \otimes \partial_k \vw_\alpha)^T \right|\, dx \\
\leq & \, K_1 \big(\alpha^2  \| \vec{\nabla}A(\cdot)\|_{L^\infty} \| \varphi \|_{H^1}\big) \,  \| \vu_1 \|_{L^2} \big(\| \vu_2\|_{L^2}+1\big)\|\vw \|_{L^2}\, \|  \partial_k \vw_\alpha\|_{H^1}.
\end{split}
\end{equation*}
Finally, by the quantity $L_2$ given in (\ref{L2}),  we obtain the wished estimate (\ref{Estim-Cont-Dep-02}). \finpv

\medskip

\begin{Lemme} For the second term, we have
\begin{equation}\label{Estim-Cont-Dep-03}
\begin{split}
&\alpha^2 \int_{\T} \left|\partial_k\big(A(\varphi\ast \vu_1)-A(\varphi \ast \vu_2)\big) (\vec{\nabla} \otimes \vu_{2,\alpha})^T : (\vec{\nabla} \otimes \partial_k \vw_{\alpha})^T \right|\, dx \\
\leq &\, \sqrt{K_0} L_3 (\| \vu_1 \|_{L^2}+1)\| \vu_{2} \|_{L^2}\, \| \vw \|_{L^2}\,  \| \partial_k \vw_{\alpha} \|_{H^1}.   
\end{split}
\end{equation}
\end{Lemme}
\pv We start by writing the following pointwise computations:
\begin{equation*}
    \begin{split}
 &\, \left|\partial_k\big(A(\varphi\ast \vu_1)-A(\varphi \ast \vu_2)\big) \right| \\
 =&\, \left| \partial_k A (\varphi\ast \vu_1)(\partial_k \varphi)\ast \vu_1 -  \partial_k A (\varphi\ast \vu_2)(\partial_k \varphi)\ast \vu_2\right|\\
=&\, \left| \big( \partial_k A (\varphi\ast \vu_1) - \partial_k A (\varphi\ast \vu_2) \big) (\partial_k \varphi)\ast \vu_1 + \partial_k A (\varphi \ast \vu_2) (\partial_k \varphi)\ast (\vu_1-\vu_2)\right| \\
\leq &\, \left|  \partial_k A (\varphi\ast \vu_1) - \partial_k A (\varphi\ast \vu_2)\right|| (\partial_k \varphi)\ast \vu_1| + |\partial_k A (\varphi \ast \vu_2) | |(\partial_k \varphi)\ast (\vu_1-\vu_2)|.
    \end{split}
\end{equation*}
Moreover,  by the first assumption in  (\ref{conditions-a-varphi-uniqueness}), the function $\partial_k A(\cdot)$ is Lipschitz continuous and there exists  a constant $C^{'}_{A}>0$ such that we have
\begin{equation*}
    \begin{split}
 &\, \left|  \partial_k A (\varphi\ast \vu_1) - \partial_k A (\varphi\ast \vu_2)\right|| (\partial_k \varphi)\ast \vu_1| + |\partial_k A (\varphi \ast \vu_2) | |(\partial_k \varphi)\ast (\vu_1-\vu_2)|\\
 \leq &\, C^{'}_{A} | \varphi \ast (\vu_1-\vu_2)| | (\partial_k \varphi)\ast \vu_1|+|\partial_k A (\varphi \ast \vu_2) | |(\partial_k \varphi)\ast (\vu_1-\vu_2)|\\
 =&\, C^{'}_{A}  | \varphi \ast \vw| | (\partial_k \varphi)\ast \vu_1|+|\partial_k A (\varphi \ast \vu_2) | |(\partial_k \varphi)\ast \vw|
    \end{split}
\end{equation*}
We thus obtain
\begin{equation*}
 \begin{split}
&\,\alpha^2 \int_{\T} \left|\partial_k\big(A(\varphi\ast \vu_1)-A(\varphi \ast \vu_2)\big) (\vec{\nabla} \otimes \vu_{2,\alpha})^T : (\vec{\nabla} \otimes \partial_k \vw_{\alpha})^T \right|\, dx\\
\leq &\, \alpha^2 C^{'}_{A} \| \varphi \ast \vw \|_{L^\infty}\| (\partial_k \varphi)\ast \vu_1 \|_{L^\infty} \| \vec{\nabla} \otimes \vu_{2,\alpha} \|_{L^2}\, \| \vec{\nabla} \otimes \partial_k \vw_{\alpha} \|_{L^2} \\
&\, + \alpha^2 \| \vec{\nabla}A(\cdot)\|_{L^\infty}\| (\partial_k \varphi)\ast \vw \|_{L^\infty} \| \vec{\nabla} \otimes \vu_{2,\alpha} \|_{L^2}\, \| \vec{\nabla} \otimes \partial_k \vw_{\alpha} \|_{L^2}\\
\leq &\,  \alpha^2 C^{'}_{A} \| \varphi \|_{L^2} \| \vw \|_{L^2}  \, \| \partial_k \varphi \|_{L^2}\| \vu_1 \|_{L^2} \| \vec{\nabla} \otimes \vu_{2,\alpha} \|_{L^2}\,  \| \partial_k \vw_{\alpha} \|_{H^1}\\
&\, + \alpha^2 \| \vec{\nabla}A(\cdot)\|_{L^\infty} \| \partial_k \varphi \|_{L^2} \| \vw \|_{L^2}\, \| \vec{\nabla} \otimes \vu_{2,\alpha} \|_{L^2}\,  \| \partial_k \vw_{\alpha} \|_{H^1}\\
\leq &\, \alpha^2 \|\varphi \|_{H^1}\max\left( C^{'}_{A} \| \varphi\|_{H^1}, \| \vec{\nabla}A(\cdot)\|_{L^\infty}\right) \times\\
&\, \times (\| \vu_1 \|_{L^2}+1)\| \vec{\nabla} \otimes \vu_{2,\alpha} \|_{L^2}\, \| \vw \|_{L^2}\,  \| \partial_k \vw_{\alpha} \|_{H^1}.
 \end{split}
\end{equation*}
Setting the quantity $L_3$ as in expression (\ref{L3}), and using inequality (\ref{Control-Filter}), we write
\begin{equation*}
\begin{split}
&\alpha^2 \int_{\T} \left|\partial_k\big(A(\varphi\ast \vu_1)-A(\varphi \ast \vu_2)\big) (\vec{\nabla} \otimes \vu_{2,\alpha})^T : (\vec{\nabla} \otimes \partial_k \vw_{\alpha})^T \right|\, dx \\
= &\, L_3 (\| \vu_1 \|_{L^2}+1)\| \vec{\nabla} \otimes \vu_{2,\alpha} \|_{L^2}\, \| \vw \|_{L^2}\,  \| \partial_k \vw_{\alpha} \|_{H^1}\\
\leq &\, \sqrt{K_0} L_3 (\| \vu_1 \|_{L^2}+1)\| \vu_{2} \|_{L^2}\, \| \vw \|_{L^2}\,  \| \partial_k \vw_{\alpha} \|_{H^1},  
\end{split}
\end{equation*}
resulting in the wished estimate (\ref{Estim-Cont-Dep-03}). \finpv

\medskip

Third term follows already known estimates: with the quantity $L_1$ given in (\ref{L1}), and using inequality  (\ref{Filter-Control-H2}), we  write 
\begin{equation}\label{Estim-Cont-Dep-04}
    \begin{split}
 &\,   \alpha^2 \int_{\T} \left|\big(A(\varphi\ast \vu_1)-A(\varphi \ast \vu_2)\big) (\vec{\nabla}\otimes \partial_k \vu_{2,\alpha})^T : (\vec{\nabla}\otimes \partial_k \vw_{\alpha})^T \right| \, dx  \\
 \leq &\, \alpha^2 C_A \| \varphi \|_{L^2} \| \vw \|_{L^2} \| \vec{\nabla}\otimes \partial_k \vu_{2,\alpha}\|_{L^2}  \,  \| \partial_k \vw_{\alpha} \|_{H^1}\\
 = &\, L_1 \|  \vu_{2,\alpha} \|_{H^2}\, \| \vw \|_{L^2}\, \| \partial_k \vw_{\alpha} \|_{H^1}\\
 \leq &\, K_2 L_1  (\|\vu_{2}\|_{L^2}+1)\|  \vu_{2}\|_{L^2}\, \| \vw \|_{L^2}\, \| \partial_k \vw_{\alpha} \|_{H^1}.
    \end{split}
\end{equation} 

Fourth term is directly estimate by Cauchy-Schwarz inequality. Collecting all the constants in inequalities (\ref{Estim-Cont-Dep-02}), (\ref{Estim-Cont-Dep-03}) and  (\ref{Estim-Cont-Dep-04}), we set now
\begin{equation*}
    K_{3,0}= \max\left(K_1 L_2,  K_2 L_1, \sqrt{K_0}L_3\right).
\end{equation*}
We   substitute these inequalities  in estimate (\ref{Estim-Cont-Dep-01}), and  rearranging   terms, we  write
\begin{equation*}
  \begin{split}
    &\,   C_{\alpha,\beta} \| \partial_k \vw_\alpha \|^{2}_{H^1} \\
    \leq &\, K_{3,0} \Big( \| \vu_1\|_{L^2}(\| \vu_2\|_{L^2}+1)+ (\| \vu_1 \|_{L^2}+1)\| \vu_{2} \|_{L^2} + (\|\vu_{2}\|_{L^2}+1)\|  \vu_{2}\|_{L^2} + 1 \Big) \| \vw \|_{L^2}\| \partial_k \vw_\alpha \|_{H^1}\\
    \leq &\, K_{3,0} \Big(\| \vu_1 \|_{L^2}+\| \vu_2 \|_{L^2}+1\Big)^2\,  \| \vw \|_{L^2}\| \partial_k \vw_\alpha \|_{H^1},
  \end{split}  
\end{equation*}
hence we get
\begin{equation*}
 \| \partial_k \vw_\alpha \|_{H^1}  \leq  \,\frac{K_{3,0}}{C_{\alpha,\beta}} \Big(\| \vu_1 \|_{L^2}+\| \vu_2 \|_{L^2}+1\Big)^2 \| \vw \|_{L^2}=: K_3 \Big(\| \vu_1 \|_{L^2}+\| \vu_2 \|_{L^2}+1\Big)^2  \| \vw \|_{L^2}.
\end{equation*}
Finally, recalling that $C_{\alpha,\beta}=\min(\alpha^2 \beta, 1)$ and $K_0\simeq \frac{1}{\min(\alpha^2\beta,1)}$, we obtain
\[ K_3=: \frac{K_{3,0}}{C_{\alpha,\beta}}\simeq  \max\left(K_0 K_1 L_2, K_0 K_2 L_1, K^{3/2}_{0} L_3\right). \]
Proposition \ref{Cont-Dep-H2-Filter} is proven. \finpv 
\section{Global well-posedness of weak  Leray solutions}\label{Sec:GWP} 
Throughout this section, we shall use $C>0$ as a generic constant that depends on the quantities $\alpha$, $\beta$, $\nu$, and some norms of the functions $A(\cdot)$ and $\varphi$. This constant may change from one line to another.
\subsection{Proof of Theorem \ref{Th-Existence}}
We will work with the (equivalent) coupled integral-elliptic formulation of system  (\ref{Alpha-model}):
\begin{equation}\label{Mild-Formulation}
\begin{cases}\vspace{2mm}
\ds{\vu(t,\cdot)=e^{\nu t\Delta }\vu_0 + \int_{0}^{t}e^{\nu(t-s) \Delta}\vf(s,\cdot) ds - \int_{0}^{t}e^{\nu (t-s) \Delta}\P\Big((\vu_\alpha \cdot \vec{\nabla})\vu\Big)(s,\cdot)ds}, \quad \text{div}(\vu)=0,\\
\vu_\alpha(t,\cdot) = \alpha^2 \P\, \text{div}\Big(A(\varphi\ast \vu(t,\cdot)) (\vec{\nabla}\otimes \vu_\alpha(t,\cdot))^T\Big)+ \vu(t,\cdot), \quad \text{div}(\vu_\alpha)=0,
\end{cases}
\end{equation}
where  $\ds{e^{\nu t\Delta}\vg=\sum_{k \in \Z}e^{i k\cdot x}\, e^{- \nu |k|^2 t}\, \widehat{\vg}_k }$. For a time  $T>0$, we shall denote $E_T=\Big\{ \vv \in L^{\infty}([0,T],L^{2}_{\text{s}}(\T))\cap L^2([0,T], \dot{H}^1(\T): \ \text{div}(\vv)=0\Big\}$, with the usual norm
\[ \| \vv \|_{E_T}=\sup_{0\leq t \leq T}\| \vv(t,\cdot)\|_{L^2}+\sqrt{\nu}\left( \int_{0}^{T}\| \vv(t,\cdot)\|^{2}_{\dot{H}^1} dt\right)^{\frac{1}{2}}.\]

\medskip 

Following some ideas of \cite[Lemma $2.2$]{Lin}, we will construct a solution $(\vu, \vu_\alpha)$ to the system (\ref{Mild-Formulation}) as follows. For any $\vv \in E_T$,    by Propositions  \ref{Existence-Filter} and \ref{Dependence-data-Filter} there exists a unique solution $\vv_\alpha \in L^\infty([0,T],H^1(\T))$  to the elliptic problem
\[ \vv_\alpha=\alpha^2 \P\,\text{div}\Big(A(\varphi\ast \vv) (\vec{\nabla}\otimes \vv_\alpha)^T\Big)+ \vv, \quad \text{div}(\vv_\alpha)=0. \]
With the obtained function $\vv_\alpha$, and using the right-hand side of the first identity in (\ref{Mild-Formulation}), we now define the operator:
\begin{equation}\label{Operator-T-alpha}
\begin{split}
\mathcal{T}_\alpha(\vv)= &\,  e^{\nu t\Delta }\vu_0 + \int_{0}^{t}e^{\nu(t-s) \Delta}\vf(s,\cdot) ds - \int_{0}^{t}e^{\nu (t-s) \Delta}\P\Big((\vv_\alpha \cdot \vec{\nabla})\vv\Big)(s,\cdot)ds\\
=:&\, \U_0-\int_{0}^{t}e^{\nu (t-s) \Delta}\P\Big((\vv_\alpha \cdot \vec{\nabla})\vv\Big)(s,\cdot)ds.
\end{split}
\end{equation}
We shall prove that there exists $\vu \in E_T$ which is a solution to the fixed point problem  $\vu=\mathcal{T}_\alpha(\vu)$. Consequently,  the obtained couple $(\vu, \vu_\alpha) \in E_T \times L^\infty_t H^1_x$ is a solution to (\ref{Mild-Formulation}).

\medskip

Given that $\vu_0 \in L^2_{\text{s}}(\T)$ and $\vf \in L^2_{\text{loc}}([0,+\infty[, \dot{H}^{-1}(\T))$, by well-known arguments we have $\U_0 \in E_T$. For $0 < T < 1$, it holds that
\[ \| \U_0 \|_{E_T} \leq C\left( \| \vu_0 \|_{L^2}+ \| \vf \|_{L^2([0,1],\dot{H}^{-1}(\T))} \right). \]
See, for instance, \cite[Theorem $12.2$]{PLe}. Then, for a fixed $R > 0$, we define the ball
\[ B = \{ \vv \in E_T : \, \| \vv - \U_0 \|_{E_T} \leq R \}, \]
and we shall prove that for sufficiently small $0 < T < 1$, the operator $\mathcal{T}_\alpha$ defined in (\ref{Operator-T-alpha}) is contractive on $B$.

\medskip

Our starting point is to prove that $\mathcal{T}_\alpha$ maps $B$ into $B$. Let $\vv \in B$. Again by \cite[Theorem $12.2$]{PLe},   we have the known estimate 
\begin{equation}\label{Estim-Existence-01}
    \left\| \mathcal{T}_\alpha(\vv) - \U_0 \right\|_{E_T} =\, \left\| \int_{0}^{t}e^{\nu (t-s) \Delta}\P\Big((\vv_\alpha \cdot \vec{\nabla})\vv\Big)(s,\cdot)ds \right\|_{E_T} \leq C \| (\vv_\alpha \cdot \vec{\nabla})\vv \|_{L^{2}_{t}\dot{H}^{-1}_{x}}. 
\end{equation}
As $\text{div}(\vv_\alpha)=\text{div}(\vv)=0$, we can write $(\vv_\alpha \cdot \vec{\nabla})\vv = \text{div}(\vv \otimes \vv_\alpha)$. Then,  we get
\begin{equation*}
    C \| (\vv_\alpha \cdot \vec{\nabla})\vv \|_{L^{2}_{t}\dot{H}^{-1}_{x}}= C \| \text{div}(\vv \otimes \vv_\alpha) \|_{L^{2}_{t}\dot{H}^{-1}_{x}} \leq C \| \vv \otimes \vv_\alpha \|_{L^{2}_{t}L^{2}_{x}} \leq C\, T^{1/4} \| \vv \otimes \vv_\alpha \|_{L^{4}_{t}L^{2}_{x}}. 
\end{equation*}
In the last expression, we apply H\"older inequalities, first  in the spatial variable with $1/2=1/3+1/6$, and thereafter in the time variable with $1/4=1/\infty+1/4$,  to write
\begin{equation*}
C\, T^{1/4} \| \vv \otimes \vv_\alpha \|_{L^{4}_{t}L^{2}_{x}} \leq  C\, T^{1/4} \| \vv \|_{L^{4}_{t}L^{3}_{x}}\, \| \vv_\alpha \|_{L^{\infty}_{t}L^{6}_{x}}.    
\end{equation*}
Moreover, by Sobolev embeddings  we get
\begin{equation}\label{Estim-Existence-02}
 C\, T^{1/4} \| \vv \|_{L^{4}_{t}L^{3}_{x}}\, \| \vv_\alpha \|_{L^{\infty}_{t}L^{6}_{x}} \leq C\,T^{1/4} \| \vv \|_{L^{4}_{t}\dot{H}^{1/2}_{x}}\, \| \vv_\alpha \|_{L^{\infty}_{t}\dot{H}^{1}_{x}} \leq   C\, T^{1/4} \| \vv \|_{L^{4}_{t}\dot{H}^{1/2}_{x}}\, \| \vv_\alpha \|_{L^{\infty}_{t}H^{1}_{x}}.
\end{equation}
We still need to estimate both terms on the right-hand side. For the first term,  recall that $\vv \in E_T$ and by interpolation inequalities in homogeneous Sobolev spaces we have 
\begin{equation}\label{Estim-Existence-03}
\| \vv \|_{L^{4}_{t}\dot{H}^{1/2}_{x}} \leq C \| \vv \|^{1/2}_{L^{\infty}_{t}L^{2}_{x}}\, \| \vv \|^{1/2}_{L^{2}_{t}\dot{H}^{1}_{x}} \leq C \| \vv \|_{E_T}.    \end{equation}
For the second term, recall that by the control (\ref{Control-Filter}) we have 
\begin{equation*}
\| \vv_\alpha \|_{L^{\infty}_{t}H^{1}_{x}} \leq C \| \vv \|_{L^{\infty}_{t}L^{2}_{x}} \leq  C  \| \vv \|_{E_T}.
\end{equation*}
We thus get
\begin{equation*}
C\, T^{1/4} \| \vv \|_{L^{4}_{t}\dot{H}^{1/2}_{x}}\, \| \vv_\alpha \|_{L^{\infty}_{t}H^{1}_{x}} \leq C\, T^{1/4}  \| \vv \|^{2}_{E_T}.   
\end{equation*}
Gathering these estimates, for $\vv \in B$ we obtain 
\begin{equation*}
    \begin{split}
    \left\| \mathcal{T}_\alpha(\vv) - \U_0 \right\|_{E_T} \leq &\,  C\,  T^{1/4}  \| \vv \|^{2}_{E_T}  \leq C\,  T^{1/4} ( \| \vv - \U_0 \|^{2}_{E_T} + \| \U_0 \|^{2}_{E_T}) \\
    \leq &\, C\,  T^{1/4} \left(R^2+\Big(\| \vu_0 \|_{L^2(\Omega)}+\| \vf \|_{L^2([0,1],\dot{H}^{-1}(\T))}\Big)^2\right).
    \end{split}
\end{equation*}
We thus set the time $0<T<1$ small enough such that $ C\,T^{1/4} \left(R^2+\Big(\| \vu_0 \|_{L^2(\Omega)}+\| \vf \|_{L^2([0,1],\dot{H}^{-1}(\T))}\Big)^2\right) \leq R$, hence we have $\mathcal{T}_\alpha(\vv) \in B$. 

\medskip

Now, we must verify that $\mathcal{T}_\alpha:B \to B$ is contractive. For $\vv_1, \vv_2 \in B$, we write
\begin{equation*}
\begin{split}
&\,  \| \mathcal{T}_\alpha(\vv_1)-\mathcal{T}_\alpha(\vv_2) \|_{E_T} \\
=&\, \left\| \int_{0}^{t}e^{\nu (t-s) \Delta}\P\Big((\vv_{1,\alpha} \cdot \vec{\nabla})\vv_1\Big)(s,\cdot)ds - \int_{0}^{t}e^{\nu (t-s) \Delta}\P\Big((\vv_{2,\alpha} \cdot \vec{\nabla})\vv_2\Big)(s,\cdot)ds \right\|_{E_T}\\
\leq &\, \left\| \int_{0}^{t}e^{\nu (t-s) \Delta}\P\Big(((\vv_{1,\alpha}-\vv_{2,\alpha}) \cdot \vec{\nabla})\vv_1\Big)(s,\cdot)ds\right\|_{E_T}+\left\| \int_{0}^{t}e^{\nu (t-s) \Delta}\P\Big((\vv_{2,\alpha} \cdot \vec{\nabla})(\vv_1-\vv_2)\Big)(s,\cdot)ds\right\|_{E_T}.
 \end{split}
\end{equation*}
For the first term, by estimates (\ref{Estim-Existence-01}) and (\ref{Estim-Existence-02}) we have
\begin{equation}\label{Estim-Non-Lin-1}
\left\| \int_{0}^{t}e^{\nu (t-s) \Delta}\P\Big(((\vv_{1,\alpha}-\vv_{2,\alpha}) \cdot \vec{\nabla})\vv_1\Big)(s,\cdot)ds\right\|_{E_T} \leq C\, T^{1/4} \| \vv_1 \|_{L^{4}_{t}\dot{H}^{1/2}_{x}}\, \| \vv_{1,\alpha}-\vv_{2,\alpha} \|_{L^{\infty}_{t}H^1_x}.    
\end{equation}
The expression $\ds{\| \vv_1 \|_{L^{4}_{t}\dot{H}^{1/2}_{x}}}$ was already estimated in (\ref{Estim-Existence-03}), and we have $\ds{\| \vv_1 \|_{L^{4}_{t}\dot{H}^{1/2}_{x}} \leq C \| \vv_1 \|_{E_T}}$. Moreover, concerning the expression $\ds{\| \vv_{1,\alpha}-\vv_{2,\alpha} \|_{L^{\infty}_{t}H^1_x}}$, by estimate (\ref{Dependence-data}) proven in Proposition \ref{Dependence-data-Filter} we have 
\begin{equation*}
 \|\vv_{1,\alpha}-\vv_{2,\alpha} \|_{L^{\infty}_{t}H^1_x} \leq  C\Big( \| \vv_2 \|_{L^{\infty}_{t}L^{2}_{x}} +1\Big)\| \vv_1 - \vv_2 \|_{L^{\infty}_{t}L^{2}_{x}}
 \leq  C \Big(\| \vv_2\|_{E_T}+1\Big)\| \vv_1 - \vv_2 \|_{E_T}.
\end{equation*}
We thus obtain:
\begin{equation}\label{Estim-Existence-04}
\left\| \int_{0}^{t}e^{\nu (t-s) \Delta}\P\Big(((\vv_{1,\alpha}-\vv_{2,\alpha}) \cdot \vec{\nabla})\vv_1\Big)(s,\cdot)ds\right\|_{E_T} \leq C\, T^{1/4} \| \vv_1 \|_{E_T} (\| \vv_2\|_{E_T}+1)\| \vv_1 - \vv_2 \|_{E_T}.    
\end{equation}
The second term  verifies a similar estimate:
\begin{equation}\label{Estim-Existence-05}
 \left\| \int_{0}^{t}e^{\nu (t-s) \Delta}\P\Big((\vv_{2,\alpha} \cdot \vec{\nabla})(\vv_1-\vv_2)\Big)(s,\cdot)ds\right\|_{E_T} \leq C\, T^{1/4}\| \vv_2 \|_{E_T}\| \vv_1 - \vv_2 \|_{E_T}.   
\end{equation}
With estimates (\ref{Estim-Existence-04}) and (\ref{Estim-Existence-05}) at hand, we can write
\begin{equation}\label{Control-Uniqueness}
\| \mathcal{T}_\alpha(\vv_1)-\mathcal{T}_\alpha(\vv_2) \|_{E_T} \leq C\, T^{1/4} \Big( \| \vv_1 \|_{E_T} (\| \vv_2\|_{E_T}+1)+ \| \vv_2 \|_{E_T}\Big) \| \vv_1 - \vv_2\|_{E_T}.    
\end{equation}
As before, since $\vv_1, \vv_2 \in B$, the term $\ds{ \| \vv_1 \|_{E_T} (\| \vv_2\|_{E_T}+1)+ \| \vv_2 \|_{E_T}}$ is ultimately  controlled by an expression depending on $R,\|\vu_0\|_{L^2}$ and $\| \vf \|_{L^{2}([0,1],\dot{H}^{-1}(\T))}$, which we shall denote by $C\left(R,\| \vu_0\|_{L^2}, \| \vf \|_{L^{2}_{t}\dot{H}^{-1}_{x}}\right)>0$. We thus write
\[ \| \mathcal{T}_\alpha(\vv_1)-\mathcal{T}_\alpha(\vv_2) \|_{E_T} \leq C\left(R,\| \vu_0\|_{L^2}, \| \vf \|_{L^{2}_{t}\dot{H}^{-1}_{x}}\right)\,T^{1/4} \, \| \vv_1 - \vv_2\|_{E_T}. \]
Resetting  the time $0<T<1$ (which already verifies the smallness condition above) so that it holds $C\left(R,\| \vu_0\|_{L^2}, \| \vf \|_{L^{2}_{t}\dot{H}^{-1}_{x}}\right)\, T^{1/4} <1$, the operator $\mathcal{T}_\alpha$ is contractive on $B$.  

\medskip

At this point, we can apply the Banach's contraction principle to obtain a (local in time) solution $\vu \in E_T$ of the fixed point problem $\vu= \mathcal{T}_\alpha(\vu)$.  Moreover, this solution is uniquely defined by the condition $\vu \in B$. 

\medskip

In the next step, we shall prove that this local-in-time solution extends to a global one. Specifically, we will demonstrate that the energy control (\ref{Energy-Control}) holds with equality. To achieve this, we will need the following technical lemma:
\begin{Lemme}\label{Lemma-Transport-term} Let $\vu \in E_T$. Then we have $(\vu_\alpha \cdot \vec{\nabla}) \vu \in L^2([0,T], L^{6/5}(\T))$.
\end{Lemme}
The proof is straightforward and relies on the key information that $\vu_\alpha \in L^{\infty}([0,T],H^1(\T))$. 
With this lemma, and noting that $\vu \in E_T \subset L^2([0,T], \dot{H}^1(\T))\subset L^2([0,T], L^6(\T))$, we have that  for \emph{a.e.} $0<t \leq T$, the integral $\ds{\int_{\T} (\vu_\alpha \cdot \vec{\nabla}) \vu(t,x)\cdot \vu(t,x) dx}$ is well-defined. Moreover, by the divergence-free condition of $\vu$ and $\vu_\alpha$, we have 
\[\int_{\T} (\vu_\alpha \cdot \vec{\nabla}) \vu(t,x)\cdot \vu(t,x) dx=0. \]
Thus, through standard computations, the energy control (\ref{Energy-Control}) holds with an equality, and the solution $\vu$ of the coupled system (\ref{Alpha-model}) can be extended to the entire time interval $[0, +\infty[$.

\medskip

In the last step, we  recover the pressure terms $P$ and $P_\alpha$. On the one hand, recall that the pressure $P_\alpha \in L^\infty_{loc}([0,+\infty[, L^2(\T))$  was already obtained in Proposition \ref{Existence-Filter}. On the other hand, for the pressure $P$, since $\vu$ solves the first equation in the system (\ref{Alpha-model}) and $\text{div}(\vu) = \text{div}(\vf) = 0$, we have 
\[ \P \left ( \partial_t \vu - \nu\Delta \vu + (\vu_\alpha \cdot \vec{\nabla})\vu - \vf \right)=0.\]
By well-known properties of the Leray's projector $\P$ (see \cite[Lemma $6.3$]{PLe}) there exists a distribution $P$ such that 
\[ \partial_t \vu - \nu\Delta \vu + (\vu_\alpha \cdot \vec{\nabla})\vu - \vf= -\vec{\nabla}P.\]
Moreover, it is easy to verify that each term on the left belongs to the space $L^2_{\text{loc}}([0, +\infty[, \dot{H}^{-1}(\T))$, hence $P \in L^2_{\text{loc}}([0, +\infty[, L^2(\T))$. Theorem \ref{Th-Existence} is now proven.   \finpv  

\subsection{Proof of Theorem  \ref{Th-Uniqueness}}
Once we have Proposition \ref{Cont-Dep-H2-Filter}, the proof of Theorem \ref{Th-Uniqueness} follows from standard energy estimates. We define $\vw=\vu_1-\vu_2 \in (L^{\infty}_{t})_{loc}(L^{2}_\text{s})_{x}\cap (L^{2}_{t})_{loc}\dot{H}^{1}_{x}$, which solves the following equation:
\begin{equation}\label{Alpha-model-Difference}
\begin{cases}
\partial_t \vw - \nu \Delta \vw + \P\Big( (\vu_{1,\alpha}\cdot \vec{\nabla}) \vw + \big((\vu_{1,\alpha}- \vu_{2,\alpha}) \cdot \vec{\nabla}\big)\vu_{2}\Big)= \vf_1 - \vf_2, \quad \text{div}(\vw)=0, \\
\vw(0,\cdot)=\vu_{0,1}-\vu_{0,2}.
\end{cases}
\end{equation}
\begin{Remarque} The nonlinear filter equation (\ref{Filter-Regularized}) does not allow us to write $\vu_{1,\alpha} - \vu_{2,\alpha} = \vw_\alpha$. Consequently, the expression $\vu_{1,\alpha} - \vu_{2,\alpha}$ is not simply controlled by $\vw_\alpha$. In this sense, Proposition \ref{Cont-Dep-H2-Filter} is our key tool for controlling the term $\vu_{1,\alpha} - \vu_{2,\alpha}$.
\end{Remarque}
    
Since $\vu_{1}, \vu_2, \vw \in (L^\infty_t)_{loc} (L^2_{\text{s}})_x \cap (L^2_t)_{loc} \dot{H}^1_x$, by Lemma \ref{Lemma-Transport-term} we have $(\vu_{1,\alpha} \cdot \nabla) \vw \in (L^2_t)_{loc} L^{6/5}_x$ and $((\vu_{1,\alpha} - \vu_{2,\alpha}) \cdot \nabla) \vu_2 \in (L^2_t)_{loc} L^{6/5}_x$. Therefore, from the equation above, the divergence-free property of $\vw$, and after integrating by parts, we can write
\begin{equation}\label{Iden-Cont-Dep-00}
\frac{1}{2} \frac{d}{dt} \| \vw(t, \cdot) \|^2_{L^2} + \nu \| \vw(t, \cdot) \|^2_{\dot{H}^1} + \int_{\T} \Big(((\vu_{1,\alpha} - \vu_{2,\alpha}) \cdot \nabla) \vu_2 \Big) \cdot \vw  dx = \left\langle \vf_1 - \vf_2, \vw \right\rangle_{\dot{H}^{-1} \times \dot{H}^1}.
\end{equation}

In this identity we must estimate the term $\ds{\int_{\T} \Big(((\vu_{1,\alpha}- \vu_{2,\alpha}) \cdot \vec{\nabla})\vu_{2}\Big)\cdot \vw\, dx}$. Using H\"older inequalities and   Sobolev embeddings,  we get
\begin{equation}\label{Estim-Cont-Dep-00}
\begin{split}
 &\, \left| \int_{\T} \Big(((\vu_{1,\alpha}- \vu_{2,\alpha}) \cdot \vec{\nabla})\vu_{2}\Big)\cdot \vw\, dx \right| 
\leq  \, \| ((\vu_{1,\alpha}- \vu_{2,\alpha}) \cdot \vec{\nabla})\vu_{2} \|_{L^2}\, \| \vw \|_{L^2}\\
\leq &\, \| \vu_{1,\alpha}- \vu_{2,\alpha} \|_{L^\infty}\, \| \vec{\nabla}\otimes \vu_{2} \|_{L^2}\, \| \vw \|_{L^2} \leq \, C\, \| \vu_{1,\alpha}- \vu_{2,\alpha} \|_{H^2}\, \| \vu_{2} \|_{\dot{H}^1}\, \| \vw \|_{L^2}.
\end{split}
\end{equation}
At this point, recall that the expression $\ds{\| \vu_{1,\alpha}- \vu_{2,\alpha} \|_{H^2}}$ was already estimated in inequality (\ref{Control-Filter-Diff-H2}). We  thus write
\begin{equation*}
    \begin{split}
   &\, C\, \| \vu_{1,\alpha}- \vu_{2,\alpha} \|_{H^2}\, \|  \vu_{2} \|_{\dot{H}^1}\, \| \vw \|_{L^2}\\
 \leq &\, C\,K_3 \left( \Big(\| \vu_1 \|_{L^2}+\| \vu_2 \|_{L^2}+1\Big)^2\, \|\vw\|_{L^2}\right)  \| \vu_{2} \|_{\dot{H}^1} \| \vw \|_{L^2}\\
 \leq &\, C\, K_3 \Big(\| \vu_1 \|_{L^2}+\| \vu_2 \|_{L^2}+1\Big)^2 \| \vu_2 \|_{\dot{H}^1}\,   \| \vw \|^{2}_{L^2}.
    \end{split}
\end{equation*}
Setting 
\begin{equation}\label{g}
    g\big(\vu_1, \vu_2\big):= C\, K_3 \Big(\| \vu_1 \|_{L^2}+\| \vu_2 \|_{L^2}+1\Big)^2 \| \vu_2 \|_{\dot{H}^1},
\end{equation}
we finally have
\begin{equation}\label{Estim-Non-Lin-Uniq}
 \left| \int_{\Omega} \Big(((\vu_{1,\alpha}- \vu_{2,\alpha}) \cdot \vec{\nabla})\vu_{2}\Big)\cdot \vw\, dx \right| \leq g\big(\vu_1, \vu_2\big) \| \vw \|^2_{L^2}.    
\end{equation}

Once we have this estimate, we get back to identity (\ref{Iden-Cont-Dep-00}) to write
\begin{equation*}
\begin{split}
&\,\frac{1}{2}\frac{d}{dt}\| \vw(t,\cdot)\|^{2}_{L^2} + \nu \| \vw(t,\cdot)\|^{2}_{\dot{H}^1}
\leq     \left| \int_{\Omega} \Big(((\vu_{1,\alpha}- \vu_{2,\alpha}) \cdot \vec{\nabla})\vu_{2,\alpha}\Big)\cdot \vw\, dx \right| + \left| \left\langle \vf_1 - \vf_2, \vw \right\rangle_{\dot{H}^{-1}\times \dot{H}^1}  \right| \\
  \leq &\, g\big(\vu_1(t,\cdot),\vu_2(t,\cdot)\big)  \| \vw(t,\cdot)\|^2_{L^2}+ \| \vf_1(t,\cdot) - \vf_2(t,\cdot) \|_{\dot{H}^{-1}}\, \| \vw(t,\cdot)\|_{\dot{H}^1}\\
  \leq &\,  g\big(\vu_1(t,\cdot),\vu_2(t,\cdot)\big) \, \| \vw(t,\cdot)\|^2_{L^2}+ \frac{2}{\nu} \| \vf_1(t,\cdot) - \vf_2(t,\cdot) \|^2_{\dot{H}^{-1}} + \frac{\nu}{2} \| \vw(t,\cdot)\|^2_{\dot{H}^{-1}}, 
\end{split}
  \end{equation*}
 and then
\[\frac{d}{dt}\| \vw(t,\cdot)\|^{2}_{L^2} + \nu \| \vw(t,\cdot)\|^{2}_{\dot{H}^1} \leq g\big(\vu_1(t,\cdot),\vu_2(t,\cdot)\big)  \, \| \vw(t,\cdot)\|^2_{L^2}+ \frac{4}{\nu} \| \vf_1(t,\cdot) - \vf_2(t,\cdot) \|^2_{\dot{H}^{-1}}. \]
Denote the quantity $\ds{e(t):=\| \vw(t,\cdot)\|^{2}_{L^2} + \nu \int_{0}^{t}\| \vw (s,\cdot)\|^2_{\dot{H}^1} ds}$. Then,   we obtain
\begin{equation*}
        \frac{d}{dt} e(t) \leq   g\big(\vu_1(t,\cdot),\vu_2(t,\cdot)\big) e(t) + \frac{4}{\nu} \| \vf_1(t,\cdot) - \vf_2(t,\cdot) \|^2_{\dot{H}^{-1}},
\end{equation*}
and  applying the Gr\"onwall inequality, for $0\leq t \leq T$,  we have 
\begin{equation*}
\begin{split}
    e(t)\leq &\, \exp\left( \int_{0}^{t} g\big(\vu_1(s,\cdot),\vu_2(s,\cdot)\big) ds  \right)e(0) \\
    &\,+ \int_{0}^{t} \exp\left( \int_{s}^{t}g\big(\vu_1(r,\cdot),\vu_2(r,\cdot)\big) dr \right) \frac{4}{\nu} \| \vf_1(s,\cdot) - \vf_2(s,\cdot) \|^2_{\dot{H}^{-1}}ds.
    \end{split}
\end{equation*}

Now, we must estimate the integral $\ds{\int_{0}^{t} g\big(\vu_1(s,\cdot),\vu_2(s,\cdot)\big) ds}$.  By expression (\ref{g}) we write
\begin{equation}\label{Control-Integral} 
    \begin{split}
  & \int_{0}^{t} g\big(\vu_1(s,\cdot),\vu_2(s,\cdot)\big) ds\\
  =&\, C \,  K_3 \int_{0}^{t} \Big(\| \vu_1(s,\cdot) \|_{L^2}+\| \vu_2(s,\cdot) \|_{L^2}+1\Big)^2 \| \vu_2(s,\cdot) \|_{\dot{H}^1}\\
  \leq & C \, K_3 \sup_{0\leq s \leq t} \Big(\| \vu_1(s,\cdot) \|_{L^2}+\| \vu_2(s,\cdot) \|_{L^2}+1\Big)^2 \left(\int_{0}^{t}  \| \vu_2(s,\cdot) \|^{2}_{\dot{H}^1} ds \right)^{1/2} t^{1/2}\\
    \leq & C \, K_3 \sup_{0\leq s \leq t} \Big(\| \vu_1(s,\cdot) \|^2_{L^2}+\| \vu_2(s,\cdot) \|^2_{L^2}+1\Big)  \left(\nu\int_{0}^{t}  \| \vu_2(s,\cdot) \|^{2}_{\dot{H}^1} ds \right)^{1/2} \frac{t^{1/2}}{\nu^{1/2}}.
    \end{split}
\end{equation}
Then, since $\vu_1$ and $\vu_{2}$ verify the energy control (\ref{Energy-Control}), we have 
\begin{equation*}
    \begin{split}
   & K_3 \sup_{0\leq s \leq t} \Big(\| \vu_1(s,\cdot) \|^2_{L^2}+\| \vu_2(s,\cdot) \|^2_{L^2}+1\Big)  \left(\nu\int_{0}^{t}  \| \vu_2(s,\cdot) \|^{2}_{\dot{H}^1} ds \right)^{1/2} \frac{t^{1/2}}{\nu^{1/2}}\\
    \leq &\, C\,K_3 \left( \sum_{i=1}^{2} \| \vu_{0,i}\|^2_{L^2}+\frac{1}{\nu}\| \vf_i \|^2_{L^2_t \dot{H}^{-1}_x}+1\right) \left(\|\vu_{0,2}\|^2_{L^2}+\frac{1}{\nu} \| \vf_2 \|^2_{L^2_t \dot{H}^{-1}_x} \right)^{1/2} \frac{t^{1/2}}{\nu^{1/2}}\\
    \leq &\, \frac{C\,K_3}{\nu^{1/2}}  \left( \sum_{i=1}^{2} \| \vu_{0,i}\|^2_{L^2}+\frac{1}{\nu}\| \vf_i \|^2_{L^2_t \dot{H}^{-1}_x}+1\right)^{3/2} t^{1/2}.
    \end{split}
\end{equation*}

Setting
\begin{equation}\label{C0-Uniq}
 \mathfrak{C}_0:=  \frac{C\,K_3}{\nu^{1/2}}  \left( \sum_{i=1}^{2} \| \vu_{0,i}\|^2_{L^2}+\frac{1}{\nu}\| \vf_i \|^2_{L^2_t \dot{H}^{-1}_x}+1\right)^{3/2},  
\end{equation}
we have $\ds{\exp\left( \int_{0}^{t} g\big(\vu_1(s,\cdot),\vu_2(s,\cdot)\big) ds \right)\leq \mathfrak{C}_0\, t^{1/2}}$,  resulting in the first wished estimate (\ref{Dependence-date}). The second wished estimate (\ref{Dependence-data-Filter-equation}) directly follows from inequality (\ref{Dependence-data}).   Theorem \ref{Th-Uniqueness} is now proven. \finpv
  
\section{Convergence properties to classical models}\label{Sec:Convergences}
\subsection{Proof of Proposition \ref{Prop1-Conv-NS}}
Recall that  for any time $0<T<+\infty$, we  denote the space $E_T= \left\{\vv \in L^\infty([0,T],L^2_{\text{s}}(\T))\cap L^2([0,T],\dot{H}^1(\T)): \text{div}(\vv)=0\right\}$, with the usual norm $\| \vv \|_{E_T}=\ds{\sup_{0\leq t \leq T}} \| \vv(t,\cdot)\|_{L^2}+\sqrt{\nu} \left(\int_{0}^{T}\| \vv(t,\cdot)\|^2_{\dot{H}^1}dt\right)^{1/2}$.  Let $(\vu_{(\alpha)})_{\alpha>0} \subset E_T$ be the family of velocities in the  coupled system (\ref{Alpha-model-alpha}).  We  divide each step of the proof into the following technical lemmas. 
\begin{Lemme}\label{Lem-Unif-Control-Family-Alpha} There exists a quantity $0<M_1=M_1(\vu_0,\vf,\nu,T)<+\infty$, which depends on  $\vu_0, \vf, \nu$ and $T$, so that  the following uniform bounds hold:
\begin{enumerate}
    \item $\ds{\sup_{\alpha>0} \| \vu_{(\alpha)}\|_{E_T} \leq M_1}$,
    \item $\ds{\sup_{\alpha>0}\| \partial_t \vu_{(\alpha)}\|_{L^2_t H^{-3/2}_{x}} \leq M_1}$.
\end{enumerate}
\end{Lemme}
\pv The quantity $M_1$ may change from one line to another, but it does not depend of $\alpha$.  The first point is a direct consequence of the energy control (\ref{Energy-Control}), hence we write 
\begin{equation}\label{Energy-Control-Family-alpha}
 \|\vu_{(\alpha)}\|^{2}_{E_T} \leq \| \vu_0 \|^{2}_{L^2}+ \frac{1}{\nu}\int_{0}^{T}\| \vf(t,\cdot)\|^2_{\dot{H}^{-1}}dt=:M^2_1.    
\end{equation}
To verify the second point, for each $\alpha>0$ we write
\[ \partial_t \vu_{(\alpha)}=\nu \Delta \vu_{(\alpha)}-  \P\big( (\vu_{(\alpha),\alpha} \cdot \vec{\nabla}) \vu_{(\alpha)} \big) + \vf,\]
where we must verify that the first term and the second term are uniformly bounded in $L^2_t H^{-3/2}_x$. Indeed, for the first term, since $\vu_{(\alpha)}\in E_T \subset L^2_t \dot{H}^1_x$, and using  the estimate above, we have
\[ \nu\|\Delta \vu_{(\alpha)}\|_{L^2_t H^{-3/2}_x} \leq \nu \|\Delta \vu_{(\alpha)}\|_{L^2_t H^{-1}_x} \leq \nu \| \vu_{(\alpha)}\|_{L^2_t \dot{H}^{1}_x}\leq M_1.\]
For the second term, first we shall prove that  the family of filtered velocities $\big(\vu_{(\alpha),\alpha}\big)_{\alpha>0} $ verifies the uniform control
\begin{equation}\label{Unif-Control-L2-Filter}
\sup_{\alpha >0}  \| \vu_{(\alpha),\alpha} \|_{L^\infty_t L^2_x} \leq M_1. 
\end{equation}
In fact, coming back to estimate (\ref{Control-Filter-L2}), and using the Cauchy-Schwarz inequalities, we get
\[ \| \vu_{(\alpha),\alpha} \|^2_{L^2} \leq \| \vu_{(\alpha)}\|_{L^2}\| \vu_{\alpha,(\alpha)}\|_{L^2}, \]
hence 
\[ \| \vu_{\alpha,(\alpha)}\|_{L^2}\leq \| \vu_{(\alpha)}\|_{L^2}.\]
and  (\ref{Unif-Control-L2-Filter}) follows from (\ref{Energy-Control-Family-alpha}). With this uniform control at hand,  using Hardy-Littlewood-Sobolev inequalities and H\"older inequalities, we can write
\begin{equation*}
   \begin{split}
  & \left\| \P\big( \vu_{(\alpha),\alpha} \cdot \vec{\nabla}) \vu_{(\alpha)} \big) \right\|_{L^2_t H^{-3/2}_x} =  \,    \left\| \P\Big( \text{div} \Big( \vu_{(\alpha)} \otimes \vu_{(\alpha),\alpha} \Big)\Big)\right\|_{L^2_t H^{-3/2}_x}\\
   \leq &\, C \| \vu_{(\alpha)} \otimes \vu_{(\alpha),\alpha}  \|_{L^2_t \dot{H}^{-1/2}_x} \leq \, C \| \vu_{(\alpha)} \otimes \vu_{(\alpha),\alpha}  \|_{L^2_t L^{3/2}_x}\\
   \leq &\, C \| \vu_{(\alpha)}\|_{L^2_t L^6_x}\, \| \vu_{(\alpha),\alpha} \|_{L^\infty_t L^2_x} \leq \, C \| \vu_{(\alpha)}\|_{L^2_t \dot{H}^{1}_x}\, \| \vu_{(\alpha),\alpha} \|_{L^\infty_t L^2_x} \leq M_1.
   \end{split} 
\end{equation*}
Lemma \ref{Lem-Unif-Control-Family-Alpha} is proven. \finpv 

\medskip

Once we have the uniform controls proven in Lemma \ref{Lem-Unif-Control-Family-Alpha}, by the Banach-Alaoglu theorem and the Aubin-Lions lemma there exists a sub-sequence $(\alpha_k)_{k\in \mathbb{N}}$, where $\alpha_k \to 0$ as $k\to +\infty$, and there exists a vector field $\vu\in E_T$, such that  for any $0<T<+\infty$ it holds:
\begin{equation}\label{Conv-Filter-01}
\vu_{(\alpha_k)} \to \vu, \ \ \mbox{in the weak-$*$ topology of $E_T$}, 
\end{equation}
and
\begin{equation}\label{Conv-Filter-02}
 \vu_{(\alpha_k)} \to \vu, \ \ \mbox{in the strong topology of $L^2_t L^2_x$}.   
\end{equation}
Now, we must prove that the limit $\vu \in E_T$ is a weak Leray solution of   the classical Navier-Stokes equation (\ref{Navier-Stokes}).  In order to simplify our writing, from now on we will write the velocity  $\vu_{(\alpha_k)}=\vu_k$, and the filtered velocity $\vu_{(\alpha_k),\alpha_k}= \vu_{k,\alpha}$.

\medskip

Our starting point is to study the convergence of the family of filtered  velocities $(\vu_{k,\alpha})_{k\in \mathbb{N}}\subset L^\infty_t H^1_x$. 
\begin{Lemme}\label{Lem-Conv-Filter-Alpha}  We have $\vu_{k,\alpha} \to \vu$, in the strong  topology of $L^p_t L^{2}_x$,  with $1\leq p <+\infty$. 
\end{Lemme}
\pv  First, let consider $p=2$. We write 
\[ \|\vu_{k,\alpha} - \vu \|_{L^2_t L^2_x} \leq \| \vu_{k,\alpha} - \vu_k  \|_{L^2_t L^2_x}+ \|  \vu_k  - \vu\|_{L^2_t L^2_x}. \]
To study the first term on the right, recall that  for every $k$ the function $\vu_{k,\alpha}$ is a weak solution of the filter equation:
\[  -\alpha^2_k \, \P\, \text{div}\Big(A(\varphi \ast \vu_k ) \big( \vec{\nabla} \otimes \vu_{k,\alpha} \big)^T\Big)+ \vu_{k,\alpha}=  \vu_k,\]
hence we write
\begin{equation}\label{Conv-01}
\begin{split}
 \|\vu_{k,\alpha} - \vu_k   \|_{L^2_t L^2_x} =&\,  \alpha^2_k \left\| \P\, \text{div}\Big(A(\varphi \ast \vu_k ) \big( \vec{\nabla} \otimes \vu_{k,\alpha}\big)^T\Big) \right\|_{L^2_t L^2_x}\\
 =&\, \alpha_k \left(\alpha_k \left\| \P\, \text{div}\Big(A(\varphi \ast \vu_k ) \big( \vec{\nabla} \otimes \vu_{k,\alpha}\big)^T\Big) \right\|_{L^2_t L^2_x}\right).
 \end{split}
\end{equation}
Then, we will show that there exist  $\delta>0$ such that  the following uniform control holds:
\begin{equation}\label{Lim-Estim-00}
    \sup_{0<\alpha_k<\delta} \left( \alpha_k \left\| \P\,\text{div}\Big(A(\varphi \ast \vu_k ) \big( \vec{\nabla} \otimes \vu_{k,\alpha} \big)^T\Big) \right\|_{L^2_t L^2_x} \right) \leq M_2,
\end{equation}
with a quantity $0<M_2=M_2(A(\cdot), \varphi, \beta, M_1)<+\infty$.  As before, in forthcoming computations this quantity $M_2$ may change from one line to another but it is always independent of $\alpha_k$.  Using H\"older inequalities,  Young inequalities and well-known properties of $\P$, we have 
\begin{equation}\label{Lim-Estim-01}
    \begin{split}
&\alpha_k\left\|\P\, \text{div}\Big(A(\varphi \ast \vu_k) \big( \vec{\nabla} \otimes \vu_{k,\alpha} \big)^T\Big) \right\|_{L^2_t L^2_x}\\
\leq &\, \alpha_k C\,T^{1/2} \left\| \text{div}\Big(A(\varphi \ast \vu_k\big( \vec{\nabla} \otimes \vu_{k,\alpha} \big)^T\Big) \right\|_{L^\infty_t L^2_x}\\
\leq &\,C\, T^{1/2} \| \vec{\nabla} A(\cdot) \|_{L^\infty}\, \| \vec{\nabla}\varphi \|_{L^2} \| \vu_k \|_{L^\infty_t L^2_x}\, \alpha_k \| \vec{\nabla}\otimes \vu_{k,\alpha} \|_{L^\infty_t L^2_x} \\
&\, + C\,T^{1/2}\| A(\cdot)\|_{L^\infty} \alpha_k  \| \Delta \vu_{\alpha,k} \|_{L^\infty_t L^2_x}\\
\leq &\, M_2 \| \vu_k \|_{L^\infty_t L^2_x}\, \alpha_k \| \vec{\nabla}\otimes \vu_{k,\alpha} \|_{L^\infty_t L^2_x} + M_2 \,  \alpha_k  \| \Delta \vu_{k,\alpha} \|_{L^\infty_t L^2_x}.
    \end{split}
\end{equation}
Here,  we still need to control each term above. First, by estimate (\ref{Energy-Control-Family-alpha}) we have 
\begin{equation}\label{Lim-Estim-02}
\| \vu_k \|_{L^\infty_t L^2_x} \leq M_1.    
\end{equation}
Thereafter, to control the term $\alpha_k \| \vec{\nabla}\otimes \vu_{\alpha,k} \|_{L^\infty_t L^2_x}$ uniformly with respect to $k$, we return  estimate (\ref{Control-Filter}). Here, since $\alpha_k \to 0$, we have $\min(\alpha^2_k \beta, 1/2)= \alpha^2_k \beta$ for all $0<\alpha_k < \delta$, with $0<\delta <1/\sqrt{2\beta}$. Consequently, we can write
\begin{equation}\label{Lim-Estim-03}
\alpha_k \| \vec{\nabla}\otimes \vu_{k,\alpha} \|_{L^\infty_t L^2_x} \leq \frac{1}{\sqrt{\beta}} \| \vu_k \|_{L^\infty_t L^2_x} \leq \frac{M_1}{\sqrt{\beta}} \leq M_2.   
\end{equation}
Similarly, to  control the term $\alpha_k \| \Delta \vu_{k,\alpha} \|_{L^\infty_t L^2_x}$, we come back to estimate (\ref{Control-Filter-H2}). Here, for any $0<\alpha_k <\delta$ we also have  $\min(\alpha^2_k \beta, 1)=\alpha^2_k\beta$, and we write
\begin{equation}\label{Lim-Estim-04}
\alpha_k \| \Delta \vu_{k,\alpha} \|_{L^\infty_t L^2_x} \leq \frac{C_2}{\sqrt{\beta}}\Big(\| \vu_k \|_{L^\infty_t L^2_x} +1 \Big)\| \vu_k \|_{L^\infty_t L^2_x} \leq \frac{C_2}{\sqrt{\beta}} (M_1+1)M_1\leq M_2.   \end{equation}
Gathering together inequalities  (\ref{Lim-Estim-02}), (\ref{Lim-Estim-03}) and (\ref{Lim-Estim-04}) in estimate (\ref{Lim-Estim-01}), we get the wished uniform control (\ref{Lim-Estim-00}). With this uniform control at hand, by (\ref{Conv-01}) we write 
\[ \| \vu_{k,\alpha}- \vu_k \|_{L^2_tL^2_x} \leq  \alpha_k \, M_2. \]
Consequently, from this last estimate and from the convergence (\ref{Conv-Filter-02}), for $p=1$ we obtain the convergence stated in Lemma \ref{Lem-Conv-Filter-Alpha}. This also holds for $1\leq p < 2$. Moreover, for $2<p<+\infty$, recall that by (\ref{Unif-Control-L2-Filter}) the sequence  of filtered velocities $(\vu_{k,\alpha})_{k\in \mathbb{N}}$ is bounded in $L^{\infty}_t L^2_x$. Lemma \ref{Lem-Conv-Filter-Alpha} is proven.  \finpv

\medskip

Following standard arguments,  we study now the convergence of the non-linear transport terms.
\begin{Lemme}\label{Lem-Conv-Non-linear-Alpha} We have $\P\big( (\vu_{k,\alpha} \cdot \vec{\nabla}) \vu_k\big) \to \P \big( (\vu \cdot \vec{\nabla})\vu \big)$,  in the weak topology of $L^2_t \dot{H}^{-3/2}_x$.
\end{Lemme}
\pv  By Lemma \ref{Lem-Conv-Filter-Alpha} and convergence (\ref{Conv-Filter-01}) (recall that in particular we have $\vec{\nabla} \otimes \vu_{k} \to \vec{\nabla} \otimes \vu$ in the weak topology of $L^2_t L^2_x$), we obtain that $(\vu_{k,\alpha} \cdot \vec{\nabla}) \vu_k \to (\vu \cdot \vec{\nabla}) \vu$ in the weak topology of $L^q_t \dot{H}^{-3/2}_x$ for any $1 < q < 2$. Moreover, since the sequence $\big((\vu_{k,\alpha} \cdot \vec{\nabla}) \vu_k\big)_{k \in \mathbb{N}}$ is bounded in $L^2_t \dot{H}^{-3/2}_x$, this weak convergence also holds in $L^2_t \dot{H}^{-3/2}_x$. Finally, the boundedness properties of the Leray projector $\P$ in Sobolev spaces yield the desired convergence. Lemma \ref{Lem-Conv-Non-linear-Alpha} is proven.  \finpv 

\medskip

Thus, the limit $\vu$ is a weak solution of the equation 
\[\partial_t \vu + \nu \Delta \vu + \P\big((\vu \cdot \vec{\nabla})\vu\big) =\vf, \quad \text{div}(\vu)=0,\]
and by well-known arguments, there exists a distribution $P$ such that $P\in L^2_t \dot{H}^{-1/2}_x$ and $(\vu, P)$ is a weak solution of the equation (\ref{Navier-Stokes}). Moreover, since each function $\vu_k$ verifies the energy control (\ref{Energy-Control}), we can apply classical tools (see \cite[Theorem 12.2]{PLe}) to find that the limit $\vu$ fulfills this inequality.  

\medskip

Continuing with  the proof of Proposition \ref{Prop1-Conv-NS}, we verify the convergence of the pressure terms. Note that, for any $k\in \mathbb{N}$, we write
\[ \vec{\nabla} P_k= - \partial_t \vu_k + \nu \Delta \vu_k - (\vu_{k,\alpha} \cdot\vec{\nabla}) \vu_k + \vf,\]
and we also write
\[ \vec{\nabla} P= - \partial_t \vu + \nu \Delta \vu - (\vu \cdot\vec{\nabla}) \vu + \vf. \] 
On the right side,  each term converges to the corresponding term  in the weak topology of $L^2_t \dot{H}^{-3/2}_{x}$. We thus have the convergence $P_k \to P$ in the weak topology of $L^2_t \dot{H}^{-1/2}_{x}$. 

\medskip

Finally, recall that the pressure $P_{k,\alpha}$ is related to the velocity $\vu_k$ and the filtered velocity $\vu_{k,\alpha}$ by the second expression in (\ref{Characterization-pressions}). Using this expression together with  (\ref{Lim-Estim-03}),  we can write
\begin{equation*}
\| P_{k,\alpha}\|_{L^\infty_t L^2_x} \leq \alpha^2_k \, C\, \| A(\varphi \ast \vu_k) (\vec{\nabla}\otimes \vu_{k,\alpha})^T\|_{L^\infty_tL^2_x} \leq \alpha_k \, C\,  \|A(\cdot)\|_{L^\infty}\, \alpha_k\, \| \vec{\nabla}\otimes \vu_{k,\alpha}\|_{L^\infty_t L^2_x} \leq \alpha_k \, C\, \| A(\cdot)\|_{L^\infty}M_2.
\end{equation*}
Consequently, when $\alpha_k \to 0$, we obtain the convergence $P_{k,\alpha}\to 0$ in the strong topology of $L^\infty_t L^2_x$. 
Proposition  \ref{Prop1-Conv-NS} is now proven. \finpv

\subsection{Proof of Proposition \ref{Prop2-Conv-NSalpha}}
Let $(\vu_\beta)_{0<\beta <1}\subset E_T$ be the family of velocities in  the coupled system (\ref{Alpha-model-beta})   with the fixed  parameter  $\alpha>0$. 

\medskip

First, note that estimates (\ref{Energy-Control-Family-alpha}) and (\ref{Unif-Control-L2-Filter}) also uniformly  hold with respect to the parameter $\beta$. Consequently, the uniform estimates proven in Lemma \ref{Lem-Unif-Control-Family-Alpha} are valid for the family $(\vu_\beta)_{0<\beta <1}$. By invoking   the Banach-Alaoglu theorem and the Aubin-Lions lemma, there exists a sub-sequence $(\beta_k)_{k\in \mathbb{N}}$, where $\beta_k \to 1$ as  $k\to +\infty$, and there exists a vector field $\vu\in E_T$ such that  for any $0<T<+\infty$, it holds:
\begin{equation}\label{Conv-Filter-01-beta}
\vu_{\beta_k} \to \vu, \ \ \mbox{in the weak-$*$ topology of $E_T$}, 
\end{equation}
and
\begin{equation}\label{Conv-Filter-02-beta}
 \vu_{\beta_k} \to \vu, \ \ \mbox{in the strong topology of $L^2_t L^2_x$}.   
\end{equation}
We shall prove that $\vu$ is a weak Leray solution to the  Leray-alpha model  (\ref{Alpha-model-classical}), moreover, we shall prove that the convergence (\ref{Conv-Filter-01-beta}) holds in the strong topology of $E_T$. 

\medskip

As before,  to simplify of writing we shall denote  $\vu_{\beta_k}=\vu_k$ and $A_{\beta_k}(\cdot)=A_k$, where this function verifies $\beta_k \leq A_k(\cdot)\leq 1$. Moreover, we shall write the filtered velocity $\vu_{\beta_k,\alpha}=\vu_{k,\alpha}$. 

\medskip

We start by studying the convergence of the  family of filtered velocities 
$(\vu_{k,\alpha})_{k\in \mathbb{N}} \subset L^\infty_t H^1_x$, which for every fixed  $k$  solves the nonlinear filter equation (\ref{Filter-Regularized}):
\[  -\alpha^2 \P\,\text{div}\Big(A_k(\varphi \ast \vu_k) \big( \vec{\nabla} \otimes \vu_{k,\alpha} \big)^T\Big)+ \vu_{k,\alpha}=  \vu_k, \quad \text{div}(\vu_{k,\alpha})=0,\]
to the filtered limit $\vu_\alpha \in L^\infty_t H^1_x$, which solves  the linear  filter Helmholtz equation:
\[ -\alpha^2  \Delta \vu_\alpha + \vu_\alpha =  \vu, \quad \text{div}(\vu_\alpha)=0.\]

In this context, our first lemma is as follows:  
\begin{Lemme}\label{Lem-Conv-Filter-Eq} We have $\vu_{k,\alpha}\to \vu_\alpha$ in the strong topology of $L^2_t H^1_x$.
\end{Lemme}
\pv The proof follows from the next estimate:
\begin{equation}\label{Estim-Conv-Filter-Eq}
\min(\alpha^2 \beta_k, 1)\,\| \vu_{k,\alpha} - \vu_\alpha\|_{L^2_t H^1_x} \leq \alpha^2(1-\beta_k)T\|\vu_\alpha\|_{L^\infty_t H^1_x}+ \| \vu_k - \vu \|_{L^2_t L^2_x},
\end{equation}
which can be proven by following  similar ideas in the proof of Proposition \ref{Dependence-data-Filter}. Indeed,  first remark that, due to the divergence-free property of $\vu_\alpha$, one can write 
\[  \Delta \vu_\alpha = \Delta \P(\vu_\alpha)= \P \big( \Delta \vu_\alpha \big)= \P\, \text{div}((\vec{\nabla}\otimes \vu_\alpha)^T). \]
Then, the linear Helmholtz equation above is rewritten as: 
\[ -\alpha^2 \P\,  \text{div} \left( (\vec{\nabla}\otimes \vu_\alpha)^T \right)-\vu_\alpha=\vu. \]
Consequently,  $\vu_{k,\alpha} - \vu_\alpha$ solves the  following elliptic problem:
\begin{equation}\label{PDE-Filter-difference-Conv}
\begin{split}
-\alpha^2 \P\, \text{div}\Big(A_{k}(\varphi \ast \vu_k) &(\vec{\nabla}\otimes (\vu_{k,\alpha}-\vu_{\alpha}))^T +\big[A_{k}(\varphi\ast \vu_k)-1\big](\vec{\nabla}\otimes \vu_{\alpha})^T\Big)\\
&+(\vu_{k,\alpha}-\vu_{\alpha})= \vu_k - \vu.
\end{split}
\end{equation}
Following the computations performed in (\ref{Estim-Filter-01}), we have
\begin{equation}\label{Estim-Filter-01-Conv}
\begin{split}
\min( \alpha^2 \beta_k,1) \| \vu_{\alpha,k}-\vu_{\alpha}\|^{2}_{H^1}  \leq &\, \alpha^2 \int_{\T} \big|A_k(\varphi\ast \vu_k)-1\big|\, \big|\vec{\nabla}\otimes \vu_{\alpha}\big| \, \big|\vec{\nabla}\otimes (\vu_{k,\alpha}-\vu_{\alpha})\big| \, dx \\
&\, + \int_{\T} |\vu_k - \vu|\,|\vu_{k,\alpha}-\vu_{\alpha}| dx.
\end{split}
\end{equation}
We must estimate each term on the right.  For the first term, recall that for every $k\in \mathbb{N}$ and for every $x\in \Rt$ we have the pointwise bounds $\beta_k \leq A_k(x)\leq 1$. Then, we get $\| A_k(\cdot)-1\|_{L^\infty} \leq (1-\beta_k)$. We thus write 
\begin{equation*}
 \alpha^2 \int_{\T} \big|A_k(\varphi\ast \vu_k)-1\big|\, \big|\vec{\nabla}\otimes \vu_{\alpha}\big| \, \big|\vec{\nabla}\otimes (\vu_{\alpha,k}-\vu_{\alpha})\big| \, dx  \leq  \alpha^2 (1-\beta_k)\| \vu_\alpha \|_{H^1}\| \vu_{\alpha,k}-\vu_\alpha \|_{H^1}.
\end{equation*}
The second term is directly estimated as follows:
\begin{equation*}
 \int_{\T} |\vu_k - \vu|\,|\vu_{k,\alpha}-\vu_{\alpha}| dx \leq \| \vu_k - \vu \|_{L^2}\| \vu_{k,\alpha}-\vu_{\alpha}\|_{H^1}.   
\end{equation*}
Therefore, coming back to (\ref{Estim-Filter-01-Conv}), for any time $0<t<T$ we get
\begin{equation*}
\min( \alpha^2 \beta_k,1/2) \| \vu_{k,\alpha}(t,\cdot)-\vu_{\alpha}(t,\cdot)\|_{H^1}\leq \alpha^2(1-\beta_k) \| \vu_\alpha(t,\cdot) \|_{H^1}+ \| \vu_k(t,\cdot) - \vu(t,\cdot)\|_{L^2}.    
\end{equation*} 
Finally, taking the $L^2$-norm in the time variable we obtain the desired estimate (\ref{Estim-Conv-Filter-Eq}). Lemma \ref{Lem-Conv-Filter-Eq} is proven.   \finpv  

\medskip

With Lemma \ref{Lem-Conv-Filter-Eq}  at hand,  we can follow  the same ideas in the proof of Lemma \ref{Lem-Conv-Non-linear-Alpha},  and using  the convergence (\ref{Conv-Filter-02-beta}), we have the convergence $\ds{\P\big( \vu_{\alpha,k} \cdot \vec{\nabla}) \vu_k \big) \to \P \big( (\vu_\alpha \cdot \vec{\nabla})\vu \big)}$, in the weak topology of $L^2_t \dot{H}^{-3/2}_x$.  Consequently,  the limit $\vu$ is a weak solution of  the coupled system (\ref{Alpha-model-classical}). 

\medskip

To finish the proof of Proposition \ref{Prop2-Conv-NSalpha}, we verify that the convergence (\ref{Conv-Filter-01-beta}) holds in the strong topology of the space $E_T$. Indeed, since $\vu_k$ is a solution of equation (\ref{Alpha-model-beta}), and $\vu$ is a solution of equation (\ref{Alpha-model-classical}), by the first equations of these system we can write
\begin{equation*}
    \vu_k(t,\cdot)= e^{\nu t\Delta }\vu_0 + \int_{0}^{t}e^{\nu(t-s) \Delta}\vf(s,\cdot) ds - \int_{0}^{t}e^{\nu (t-s) \Delta}\P\Big((\vu_{k,\alpha} \cdot \vec{\nabla})\vu_k\Big)(s,\cdot)ds,
\end{equation*}
and 
\begin{equation*}
    \vu(t,\cdot)= e^{\nu t\Delta }\vu_0 + \int_{0}^{t}e^{\nu(t-s) \Delta}\vf(s,\cdot) ds - \int_{0}^{t}e^{\nu (t-s) \Delta}\P\Big((\vu_{\alpha} \cdot \vec{\nabla})\vu\Big)(s,\cdot)ds.
\end{equation*}
Then, for any time $0<T<+\infty$, we have
\begin{equation}\label{Estim-Conv-Strong}
\begin{split}
\| \vu_k - \vu \|_{E_T} \leq & \, \left\| \int_{0}^{t}e^{\nu (t-s) \Delta}\P\Big(((\vu_{k,\alpha}-\vu_\alpha )\cdot \vec{\nabla})\vu\Big)(s,\cdot)ds\right\|_{E_T}\\
&\, + \left\|\int_{0}^{t}e^{\nu (t-s) \Delta}\P\Big((\vu_{\alpha} \cdot \vec{\nabla})(\vu_k - \vu)\Big)(s,\cdot)ds \right\|_{E_T}.
\end{split}
\end{equation}
The first term on the right was already estimated in (\ref{Estim-Non-Lin-1}). Moreover, using interpolation inequalities and the uniform control (\ref{Energy-Control-Family-alpha}) (applied to the family $(\vu_k)_{k\in \mathbb{N}}$), one gets
\begin{equation*}
\begin{split}
 &\, \left\| \int_{0}^{t}e^{\nu (t-s) \Delta}\P\Big(((\vu_{k,\alpha}-\vu_\alpha )\cdot \vec{\nabla})\vu\Big)(s,\cdot)ds\right\|_{E_T} \leq \,  C\, T^{1/4} \| \vu_k \|_{L^4_t \dot{H}^{1/2}_x} \| \vu_{k,\alpha} - \vu_\alpha \|_{L^\infty_t H^1_x} \\
 \leq & \, C\, T^{1/4} \| \vu_k \|_{E_T} \| \vu_{k,\alpha} - \vu_\alpha \|_{L^\infty_t H^1_x} \leq \, C\, T^{1/4} M_1\,  \| \vu_{k,\alpha} - \vu_\alpha \|_{L^\infty_t H^1_x}.
 \end{split}
\end{equation*}
To estimate the second term on the right, we use again estimate (\ref{Estim-Non-Lin-1}) and interpolation inequalities  to write
\begin{equation*}
  \left\|\int_{0}^{t}e^{\nu (t-s) \Delta}\P\Big((\vu_{\alpha} \cdot \vec{\nabla})(\vu_k - \vu)\Big)(s,\cdot)ds \right\|_{E_T}  \leq C\, T^{1/4} \| \vu_\alpha \|_{L^\infty_t H^1_x}\| \vu_k - \vu\|_{E_T}. 
\end{equation*}
At this point, for a time $T_0\leq T$ small enough such that $\ds{C\, T^{1/4}_0 \| \vu_\alpha \|_{L^\infty_t H^1_x}\leq \frac{1}{2}}$, we obtain
\[   \left\|\int_{0}^{t}e^{\nu (t-s) \Delta}\P\Big((\vu_{\alpha} \cdot \vec{\nabla})(\vu_k - \vu)\Big)(s,\cdot)ds \right\|_{E_{T_0}} \leq \frac{1}{2} \| \vu_k - \vu\|_{E_{T_0}}. \]

Then, gathering these estimates in (\ref{Estim-Conv-Strong}), in the interval of time $[0,T_0] \subset [0,T]$, we get
\begin{equation*}
  \frac{1}{2} \| \vu_k - \vu \|_{E_{T_0}} \leq C\, T^{1/4}_0 \| \vu_{\alpha,k}-\vu_\alpha \|_{L^\infty_t H^1_x},
\end{equation*}
and by iterative application of this argument up to the $T$, we have
\begin{equation*}
  \frac{1}{2} \| \vu_k - \vu \|_{E_{T}} \leq C\, T^{1/4} \| \vu_{\alpha,k}-\vu_\alpha \|_{L^\infty_t H^1_x}.
\end{equation*}
By Lemma \ref{Lem-Conv-Filter-Eq}, we obtain the desired convergence. 

\medskip

The convergence of the pressure terms $P_k \to P$ in the strong topology of $L^2_t L^2_x$ follows similar ideas as those presented in the proof of Proposition \ref{Prop1-Conv-NS}, so we omit the details. Finally, we will prove the convergence $\vec{\nabla}P_{k,\alpha}\to 0$ in the strong topology of $L^\infty_t \dot{H}^{-1}_x$. For any $k\in \mathbb{N}$, we write
\[ \vec{\nabla}P_{k,\alpha}= \alpha^2 \text{div}\Big( A(\varphi\ast \vu_k)(\vec{\nabla}\otimes \vu_{k,\alpha})^T\Big) - \vu_{k,\alpha}+\vu_k, \]
and we also write
\[ 0= \alpha^2 \Delta \vu_\alpha -\vu_\alpha+\vu. \]
On the right-hand  side, each term converges to the corresponding term in the strong topology of $L^\infty_t \dot{H}^{-1}_x$. Hence, we obtain the wished convergence. Proposition \ref{Prop2-Conv-NSalpha} is thus proven. \finpv 
\section{Global attractor and fractal dimension}\label{Sec:Global-Attractor}  
\subsection{Existence of a global attractor: Proof of Theorem \ref{Th-Global-Attractor}}\label{Sec:Existence-Global-Attractor} 
We start by summarizing some results from the theory of dynamical systems and by introducing some notation that we will use to prove Theorem \ref{Th-Global-Attractor}. All the theorems and definitions below involve  the semi-group $(S(t))_{t\geq 0}$ defined in (\ref{Semi-group}). Moreover, for fixed  $t\geq 0$  and  for $B\subset L^2_{\text{s}}(\T)$, we shall use the notation
\[ S(t)B= \big\{ S(t)\vu_0: \ \vu_0 \in B \big\}\subset L^2_{\text{s}}(\T). \]
\begin{Definition}[Absorbing set]\label{Def-Absorbing-Set} Let $\mathcal{B}\subset L^2_{\text{s}}(\T)$. Then $\mathcal{B}$ is an absorbing set for the semi-group $(S(t))_{t\geq 0}$  if, for any bounded set $B\subset L^2_{\text{s}}(\T)$, there exists a time $T=T(B)>0$ such that for any $t>T$, we have $S(t)B\subset \mathcal{B}$. 
\end{Definition}
\begin{Definition}[Asymptotically compact semi-group]\label{Def-Asymp-Compact} The semi-group $(S(t))_{t\geq 0}$ is asymptotically compact if, for any bounded sequence $(\vu_{0,n})_{n \in \mathbb{N}}\subset L^2_{\text{s}}(\T)$ and for time sequence  $(t_n)_{n\in \mathbb{N}}$ such that $t_n\to +\infty$ as $n\to +\infty$, there a sub-sequence $(n_k)_{k \in \mathbb{N}}$ such that $S(t_{n_k})\vu_{0,n_k}$ converges in the strong topology of $L^2_{\text{s}}(\T)$ as $k\to +\infty$.
\end{Definition}
Now, we are ready to state the key result needed to prove Theorem \ref{Th-Global-Attractor}. For a proof of this (known) result, we refer to \cite{Raugel, Temam}.
\begin{Theoreme}\label{Th-Tech-Existence-Global-Attractor} Assume that the following conditions hold:
\begin{enumerate}
    \item The semi-group $(S(t))_{t\geq 0}$ has a bounded and closed absorbing set $\mathcal{B}$ in the sense of Definition \ref{Def-Absorbing-Set}.
    \item The semi-group $(S(t))_{t\geq 0}$ is asymptotically compact in the sense of Definition \ref{Def-Asymp-Compact}.
    \item For any $t\geq 0$, the map $S(t):\mathcal{B} \to L^2_{\text{s}}(\T)$ is continuous.
\end{enumerate}
Then, the semi-group $(S(t))_{t\geq 0}$ has a \emph{unique} global attractor $\mathcal{A} \subset L^2_{\text{s}}(\T)$ in the sense of Definition \ref{Def-Global-Attractor}, which verifies the characterization given in (\ref{Characterization-global-attractor}). 
\end{Theoreme}

\medskip

{\bf Proof of Theorem \ref{Th-Global-Attractor}}. We will verify each point stated in Theorem \ref{Th-Tech-Existence-Global-Attractor}.

\medskip

{\bf Point 1.} Absorbing set.
\begin{Proposition}\label{Prop-Absorbing-Set}  Define the quantity:
\begin{equation}\label{R}
 R^2:=\frac{L^2}{2\pi^2 \nu^2}\| \vf \|^2_{\dot{H}^{-1}}>0,   
\end{equation}
and the set
	\begin{equation}\label{Absorbing-Set}
	\mathcal{B}=\left\{  \vu_0 \in L^2_{\text{s}}(\T) : \ \| \vu_0 \|^2_{L^2} \leq R^2 \right\}.
	\end{equation}
 Then $\mathcal{B}$ is an absorbing set for the semi-group $(S(t))_{t\geq 0}$ in the sense of Definition \ref{Def-Absorbing-Set}.
\end{Proposition} 
\pv First, note that every weak  Leray solution to the coupled system (\ref{Alpha-model-Proy}) verifies the energy identity
\begin{equation}\label{Energy-Control-Derivative}
    \frac{d}{dt}\| \vu(t,\cdot)\|^2_{L^2}+2\nu \| \vu(t,\cdot)\|^2_{\dot{H}^1}=\left\langle \vf , \vu(t,\cdot)\right\rangle_{\dot{H}^{-1}\times \dot{H}^1} \leq \frac{1}{\nu}\| \vf \|^{2}_{\dot{H}^{-1}}+\nu \|\vu(t,\cdot)\|^2_{\dot{H}^1}, 
\end{equation}
hence we can write 
\begin{equation*}
\frac{d}{dt}\| \vu(t,\cdot)\|^2_{L^2}  \leq - \nu  \| \vu(t,\cdot)\|^2_{\dot{H}^1} + \frac{1}{\nu}\| \vf \|^{2}_{\dot{H}^{-1}}.   \end{equation*}
By the Poincaré's inequality, we get 
\[ \| \vu(t,\cdot)\|_{L^2}\leq \frac{L}{2\pi}\| \vu(t,\cdot)\|_{\dot{H}^1},\]
and setting
\begin{equation}\label{eta}
    \eta:=\frac{4\pi^2 \nu}{L^2}>0
\end{equation}
we obtain
\begin{equation*}
\frac{d}{dt}\| \vu(t,\cdot)\|^2_{L^2}  \leq - \eta \| \vu(t,\cdot)\|^2_{L^2} + \frac{1}{\nu}\| \vf \|^{2}_{\dot{H}^{-1}}.   
\end{equation*}
Then applying the Gr\"onwall lemma, and  setting $R^2$ as in expression (\ref{R}), we find
\begin{equation}\label{Energ-Estim-Attractor-1}
    \|\vu(t,\cdot)\|^2_{L^2} \leq  e^{-\eta t}\| \vu_0 \|^2_{L^2}+\frac{R^2}{2}.
\end{equation}
Consequently, for $\| \vu_0 \|^2_{L^2}\leq R^2$  one can always find a time $T>0$ such that for $t>T$, it holds that $\ds{e^{-\eta t} R^2 \leq \frac{R^2}{2}}$, and hence $\| \vu(t,\cdot)\|^2_{L^2}\leq R^2$. Proposition \ref{Prop-Absorbing-Set} is proven. \finpv 

\medskip

{\bf Point 2}. Asymptotic compactness. 
\begin{Proposition}\label{Prop-Asymptotically compact} The semi-group $(S(t))_{t\geq 0}$ is asymptotically compact in the sense of Definition \ref{Def-Asymp-Compact}.
\end{Proposition}
\pv Let $(\vu_{0,n})_{n \in \mathbb{N}}$ be a bounded sequence in $L^2_{\text{s}}(\T)$, and let $(t_n)_{n \in \mathbb{N}}$ be a time sequence such that $t_n \to +\infty$ as $n\to +\infty$. We need to verify that $(S(t_n)\vu_{0,n})_{n \in \mathbb{N}}$ has a sub-sequence that is  strongly convergent in $L^2_{\text{s}}(\T)$.  To do this, we shall follow some of the ideas in \cite{CortezJarrin,Ilyn}, which essentially use an energy method.  

\medskip

For any $n\in \mathbb{N}$, we consider the initial time $-t_n$, and the  initial value problem:
 \begin{equation}\label{Alpha-model-n}
\begin{cases}\vspace{2mm}
\partial_t \vu_n - \nu\Delta \vu_n + \P\big( (\vu_{n,\alpha} \cdot \vec{\nabla})\vu_n \big)  = \vf,  \quad  \text{div}(\vu_n)=0,\\  \vspace{2mm}
-\alpha^2 \P\, \text{div} \left( A (\varphi \ast \vu_n) (\vec{\nabla} \otimes \vu_{n,\alpha})^T \right)+\vu_{n,\alpha} = \vu_n,  \quad \text{div}(\vu_{n,\alpha})=0, \\
\vu_n(-t_n,\cdot)=\vu_{0,n}.
\end{cases}
\end{equation}
Then, by Theorems \ref{Th-Existence} and \ref{Th-Uniqueness}  we obtain a unique solution 
\[ \vu_n \in L^{\infty}_{loc}([-t_n, +\infty[, L^2_{\text{s}}(\T))\cap L^{2}_{loc}([-t_n, +\infty[, \dot{H}^1(\T)).\]
The uniqueness of solutions,  the explicit definition of the semi-group $S(t)$ (see the expression (\ref{Semi-group})), and the fact that $\vu_n(t,\cdot)$ starts at $-t_n$  allow us to write 
\begin{equation}\label{Def-Sequence}
 S(t_n)\vu_{0,n}=\vu_{n}(0,\cdot).   
\end{equation}
Consequently,  we will prove that $(\vu_n(0,\cdot))_{n\in \mathbb{N}}$ has a sub-sequence that is strongly  convergent in the space $L^2_{\text{s}}(\T)$. The limit of this sub-sequence will be   given by $\vu_\text{e}(0,\cdot)$, where $\vu_\text{e}$ is an eternal solution to the system (\ref{Alpha-model-eternal}). Thus, our strategy is as follows: first, we shall prove the existence of this eternal solution $\vu_\text{e}$, and then we shall prove the convergence $\vu_n(0,\cdot)\to \vu_\text{e}(0,\cdot)$. 

\begin{Lemme}\label{Lemma-Eternal-Sol} Let $\vf \in \dot{H}^{-1}(\T)\cap L^2_{\text{s}}(\T)$ be the stationary external force. There exists an eternal solution to the coupled system (\ref{Alpha-model-eternal}) in the sense of Definition \ref{Def-Eternal-Solutions}.
\end{Lemme}	
\pv The eternal solution $\vu_\text{e}$ will be obtained as the limit (as $n\to +\infty$) of solutions $\vu_n$ to the initial value problem (\ref{Alpha-model-n}). Our starting point is to consider the following energy controls on $\vu_n$. First, by estimate (\ref{Energ-Estim-Attractor-1}), we have 
\begin{equation}\label{Control-Energ-1}
    \| \vu_n(t,\cdot)\|^{2}_{L^2} \leq \| \vu_{0,n}\|^{2}_{L^2}\, e^{- \eta(t+t_n)}+ \frac{R^2}{2}, \quad \mbox{for any} \ \ t \geq -t_n.
\end{equation}
Moreover, integrating the energy inequality (\ref{Energy-Control-Derivative}) on the interval of time $[t, t+T]$,  it  holds
\begin{equation}\label{Control-Energ-2}
 \nu \int_{t}^{t+T}\| \vu_n(s,\cdot)\|^{2}_{\dot{H}^1}ds \leq \| \vu_n(t,\cdot)\|^{2}_{L^2}+\frac{T}{\nu}\| \vf \|^{2}_{\dot{H}^{-1}}, \quad \mbox{for any} \ \ t\geq -t_n \ \ \mbox{and} \ \  T>0.   
\end{equation}

Recall that $(\vu_{0,n})_{n \in \mathbb{N}}$ is bounded in $L^2_{\text{s}}(\T)$, and there exists  $M>0$ such that  $\ds{\sup_{n \in \mathbb{N}}} \| \vu_{0,n}\|_{L^2}\leq M$. Therefore, by estimate (\ref{Control-Energ-1}) we obtain 
\begin{equation}\label{Unif-Control1}
\sup_{n\in \mathbb{N}}\, \sup_{t\geq -t_n}  \| \vu_n(t,\cdot)\|^{2}_{L^2} \leq M^2 + \frac{R^2}{2}.    
\end{equation}
Similarly, by estimate (\ref{Control-Energ-2}) we get
\begin{equation}\label{Unif-Control2}
  \sup_{n\in \mathbb{N}}\, \sup_{t\geq -t_n} \nu \int_{t}^{t+T}\| \vu_n(s,\cdot)\|^{2}_{\dot{H}^1} ds \leq M^2 + \frac{R^2}{2} + \frac{T}{\nu}\| \vf \|^{2}_{\dot{H}^{-1}}.   
\end{equation}

With these uniform controls at hand,  by the Banach-Alaoglu theorem  there exists 
\[ \vu_\text{e} \in L^{\infty}_{loc}\big(\R, L^2_{\text{s}}(\T)\big)\cap L^{2}_{loc}\big(\R,\dot{H}^1(\T)\big),\]
and there exists a sub-sequence $(n_k)_{k\in \mathbb{N}}$, where $n_k \to +\infty$ as $k\to +\infty$, such that 
\begin{equation}\label{Conv-weak-Attractor}
    \vu_{n_k} \to \vu_{\text{e}}, \quad \mbox{in the weak-$*$ topology of $E_T$}.
\end{equation}

Now, we must prove that $\vu_\text{e}$ verify the coupled system (\ref{Alpha-model-eternal}) in the weak  sense. To do this,  we recall that each function $\vu_n$ is a solution to (\ref{Alpha-model-n}). Thus, it is enough to prove the following convergences:
\begin{equation}\label{Convergence1}
    -\alpha^2 \P\, \text{div}\Big(A(\varphi \ast \vu_{n_k})(\vec{\nabla}\otimes \vu_{n_k,\alpha})^T\Big) \to -\alpha^2 \P\,\text{div}\Big(A(\varphi \ast \vu_{\text{e}})(\vec{\nabla}\otimes \vu_{\text{e},\alpha})^T\Big),  
\end{equation}
and 
\begin{equation}\label{Convergence2}
   \P((\vu_{n_k,\alpha}\cdot \vec{\nabla})\vu_{n_k}) \to \P((\vu_{\text{e},\alpha}\cdot \vec{\nabla})\vu_{\text{e}}),
\end{equation}
in the sense of distributions. Our starting point is to extend each solution  $\vu_n$  by zero to the whole real line $\R$. Then, for any time $0<T<+\infty$, we will prove the following uniform controls:
\begin{equation}\label{Unif-Controls}
    \sup_{n\in \mathbb{N}} \| \vu_n \|_{L^2([-T,T],H^1(\T))} \leq M_1, \quad \sup_{n\in \mathbb{N}} \|\partial_t \vu_n \|_{L^2([-T,T],H^{-1}(\T))} \leq M_1,
\end{equation}
where $0<M_1=M_1(\nu,R, \|\vf \|_{\dot{H}^{-1}},T)<+\infty$ is a quantity independent of $n$. Indeed, the first uniform bound directly follows from (\ref{Unif-Control1}) and (\ref{Unif-Control2}). To prove the second uniform bound, we write 
\[ \partial_t \vu_n = \nu \Delta \vu_n - \P \big( (\vu_{n,\alpha}\cdot \vec{\nabla})\vu_n \big) + \vf,\]
where we will verify that  each term on the right is uniformly bounded in $L^2_t H^{-1}_x$. For the first term, recall that $\vu_n \in  L^2_t \dot{H}^{1}_x$ and by estimate (\ref{Unif-Control2}) we have 
\[ \| \Delta \vu_n \|_{L^2_t H^{-1}_x} \leq  \| \Delta \vu_n \|_{L^2_t \dot{H}^{-1}_x} \leq  \| \vu_n \|_{L^2_t \dot{H}^1}   \leq  \frac{1}{\sqrt{\nu}}\left(M^2 + \frac{R^2}{2} + \frac{T}{\nu}\| \vf \|^{2}_{\dot{H}^{-1}}\right)^{1/2}.\]
For the second term,  we write
\[ \| \P( (\vu_{n,\alpha}\cdot \vec{\nabla})\vu_n ) \|_{L^2_t H^{-1}_x} \leq  \| \text{div}(\vu_n \otimes \vu_{n,\alpha})\|_{L^2_t \dot{H}^{-1}_x} \leq \| \vu_{n} \otimes \vu_{n,\alpha} \|_{L^2_t L^2_x} \leq \| \vu_{n,\alpha} \|_{L^\infty_t L^\infty_x}\| \vu_{n} \|_{L^2_t L^2_x}.\]
Then, by Sobolev embeddings, by  estimate (\ref{Filter-Control-H2}), and by the uniform control (\ref{Unif-Control1}) we get
\begin{equation*}
    \begin{split}
   \| \vu_{n,\alpha} \|_{L^\infty_t L^\infty_x}\| \vu_{n} \|_{L^2_t L^2_x} \leq &\, C\, T^{1/2}\, \| \vu_{n,\alpha}\|_{L^\infty_t H^2_x}\| \vu_n \|_{L^\infty_t L^2_x} \leq C\, T^{1/2} \, K_2 (\|\vu_n\|_{L^\infty_t L^2_x}+1 )\| \vu_n \|^2_{L^\infty_t L^2_x} \\
   \leq &\,  C\, 2T\, K_2\,  \left(\left(M^2 +\frac{R^2}{2}\right)^{1/2}+1\right)\left( M^2 +\frac{R^2}{2}\right).
    \end{split}
\end{equation*}  
Finally, for the stationary external force $\vf$ we have $\ds{\| \vf\|_{L^2_t H^{-1}_x} \leq T \| \vf \|_{L^\infty_t \dot{H}^{-1}_x} \leq T \| \vf \|_{\dot{H}^{-1}}}$. 

\medskip

Once we have the uniform controls (\ref{Unif-Controls}),  by the Aubin-Lions Lemma (see \cite[Lemma $3.2$]{Wiedemann}), it holds:
\begin{equation}\label{Convergence-tech0}
\vu_{n_k}\to \vu_\text{e}, \quad \mbox{in the strong topology of $L^2_t L^2_x$},
\end{equation}  
and with  this information, we are able to verify the convergence (\ref{Convergence1}) as follows: first, as $\vu_{n_k}$ is bounded in $L^\infty([-T,T],L^2(\T))$ and using  estimate (\ref{Dependence-data}), we  have the strong convergence of the filtered velocities 
\begin{equation}\label{Convergence-tech1}
\vu_{n_k,\alpha}\to \vu_{\text{e},\alpha}, \quad \mbox{in} \ \ L^2([-T,T],H^1(\T)),
\end{equation}
and consequently, one has  the strong convergence $\vec{\nabla}\otimes \vu_{n_k,\alpha}\to \vec{\nabla} \otimes\vu_{\text{e},\alpha}$, in $L^2([-T,T],L^2(\T))$. On the other hand, since the indicator function $A(\cdot)$ is Lipschitz continuous (recall that in all this section we assume (\ref{conditions-a-varphi-uniqueness})), we have the pointwise estimate
\[ | A(\varphi \ast \vu_{n_k}) - A(\varphi\ast \vu_\text{e})| \leq C_A | \varphi \ast (\vu_{n_k}-\vu_\text{e})|,\]
hence we write
\[ \|A(\varphi \ast \vu_{n_k}) - A(\varphi\ast \vu_\text{e}) \|_{L^2_t L^2_x}\leq C_A \| \varphi\|_{L^1}\, \| \vu_{n_k}-\vu_\text{e}\|_{L^2_t L^2_x}, \]
and by (\ref{Convergence-tech0})  we obtain the strong convergence 
\begin{equation}\label{Convergence-Tech2}
A(\varphi \ast \vu_{n_k})\to A(\varphi\ast \vu_\text{e}), \quad \mbox{in} \ \   L^2([-T,T],L^2(\T)).   
\end{equation}
Thus, the convergence (\ref{Convergence1}) follows from the convergences proven in (\ref{Convergence-tech1}) and (\ref{Convergence-Tech2}), and from well-known properties of the Leray's projector $\P$.  Now, we will prove  the convergence (\ref{Convergence2}).  By (\ref{Convergence-tech1}) and  Sobolev embeddings, we have the strong convergence
\[ \vu_{n_k,\alpha}\to \vu_{\text{e},\alpha}, \quad \mbox{in} \ \ L^2([-T,T],L^6(\T)).\]
Then, by the convergence (\ref{Convergence-tech0}), the H\"older inequalities, and using again  well-known properties of the Leray's projector, we get the convergence
\begin{equation}\label{Covergence-Tech3}
\P((\vu_{n_k,\alpha}\cdot\vec{\nabla})\vu_{n_k})=\, \P (\text{div}(\vu_{n_k}\otimes \vu_{n_k,\alpha}))
 \to \, \P (\text{div}(\vu_{\text{e}}\otimes \vu_{\text{e},\alpha}))= \P((\vu_{\text{e},\alpha}\cdot\vec{\nabla})\vu_{\text{e}}), 
\end{equation}
in the strong topology of the space $L^1([-T,T],\dot{W}^{-1,3/2}(\T))$. Lemme \ref{Lemma-Eternal-Sol} is now proven. \finpv 

\medskip

\begin{Lemme}\label{Lem-Conv-Attractor} Let $(\vu_{n_k}(0,\cdot))_{n\in \mathbb{N}}$ be the sequence defined in (\ref{Def-Sequence}), and let $\vu_\text{e}(0,\cdot)$ be the eternal solution to the system (\ref{Alpha-model-eternal}) constructed in Lemma \ref{Lemma-Eternal-Sol}. Then we have the convergence $\vu_{n_k}(0,\cdot)\to \vu_\text{e}(0,\cdot)$ in the strong topology of the space $L^2_{\text{s}}(\T)$.
\end{Lemme}
\pv  Recall  that, for any $n \in \mathbb{N}$ and for all $t \geq -t_{n_k}$, the solution $\vu_{n_k}$ of the coupled system (\ref{Alpha-model-n}) verifies the energy identity:
\begin{equation}\label{Energy-n}
\frac{1}{2} \frac{d}{dt}\Vert \vu_{n_k}(t,\cdot)\Vert^{2}_{L^2} = -\nu \Vert \vu_{n_k}(t,\cdot) \Vert^{2}_{\dot{H}^1}+ \langle \vf, \vu_{n_k}(t,\cdot)\rangle_{\dot{H}^{-1}\times \dot{H}^1} 
\end{equation}
We multiply each term  by $e^{2 t}$, and  we integrate  over the interval $[-t_n, 0]$ to get:
\begin{equation}\label{Integral}
\begin{split}
&\frac{1}{2} \Vert \vu_{n_k}(0,\cdot)\Vert^{2}_{L^2} - \frac{1}{2} e^{- 2  t_{n_k} }\Vert \vu_{0,n_k} \Vert^{2}_{L^2}-  \eta \int_{-t_{n_k}}^{0} e^{2  t} \Vert \vu_{n_k}(t,\cdot)\Vert^{2}_{L^2}  dt \\
=&  -\nu \int_{-t_{n_k}}^{0} e^{2  t} \Vert \vu_{n_k}(t,\cdot) \Vert^{2}_{\dot{H}^1}  dt  + \int_{-t_{n_k}}^{0} e^{2  t} \langle \vf, \vu_{n_k}(t,\cdot)\rangle_{\dot{H}^{-1}\times \dot{H}^1}  dt, 
\end{split}
\end{equation}
hence we obtain
\begin{equation*}
\begin{split}
\Vert \vu_{n_k}(0,\cdot)\Vert^{2}_{L^2}=& \,  e^{- 2  t_{n_k} }\Vert \vu_{0,n_k} \Vert^{2}_{L^2} + 2 \int_{-t_{n_k}}^{0} e^{2 t} \Vert \vu_{n_k}(t,\cdot)\Vert^{2}_{L^2}  dt   -2\nu \int_{-t_{n_k}}^{0} e^{2  t} \Vert \vu_{n_k}(t,\cdot) \Vert^{2}_{\dot{H}^1}  dt \\
&  +2 \int_{-t_{n_k}}^{0} e^{2  t} \langle \vf, \vu_{n_k}(t,\cdot)\rangle_{\dot{H}^{-1}\times \dot{H}^1}  dt.
\end{split}
\end{equation*}
In each term of this identity, we take  the $\ds{\limsup}$ as $n\to +\infty$ to obtain
\begin{equation}\label{Estim-energ-limsup}
\begin{split}
\limsup_{n_k\to +\infty} \Vert \vu_{n_k}(0,\cdot)\Vert^{2}_{L^2}\leq & \limsup_{n_k\to +\infty}   e^{- 2 \eta t_{n_k} }\Vert \vu_{0,n_k} \Vert^{2}_{L^2}  \\
& +  \limsup_{n_k\to +\infty} \left( 2  \int_{-t_{n_k}}^{0} e^{2  t} \Vert \vu_{n_k}(t,\cdot) \Vert^{2}_{L^2}  dt \right) \\
&  +  \limsup_{n_k\to +\infty} \left( -2\nu \int_{-t_{n_k}}^{0} e^{2  t} \Vert \vu_{n_k}(t,\cdot) \Vert^{2}_{\dot{H}^1}  dt \right) \\
&+\limsup_{n_k\to +\infty} \left(2 \int_{-t_{n_k}}^{0} e^{2 t} \langle \vf, \vu_{n_k}(t,\cdot)\rangle_{\dot{H}^{-1}\times\dot{H}^1} \right),
\end{split}
\end{equation}
where we need to study each term on the right. For the first term, recalling that the sequence $(\vu_{0,n_k})_{n\in \mathbb{N}}$ is bounded in $L^2_{\text{s}}(\T)$, we have
\begin{equation}\label{E1}
\limsup_{n_k\to +\infty}   e^{- 2  t_{n_k} }\Vert \vu_{0,n_k} \Vert^{2}_{L^2} =0.
\end{equation}
For the second  term, by the convergence (\ref{Convergence-tech0}) we have
\begin{equation}\label{E2}
  \limsup_{n_k\to +\infty} \left( 2  \int_{-t_{n_k}}^{0} e^{2  t} \Vert \vu_{n_k}(t,\cdot) \Vert^{2}_{L^2}  dt \right) = 2  \int_{-\infty}^{0} e^{2  t}\| \vu_{\text{e}}(t,\cdot)\|^2_{L^2} dt.   
\end{equation}
For the third term, by the convergence (\ref{Conv-weak-Attractor}) we 
 can write
\begin{equation*}
\liminf_{n_k\to + \infty} \left( 2\nu \int_{-t_{n_k}}^{0} e^{2  t} \Vert \vu_{n_k}(t,\cdot) \Vert^{2}_{\dot{H}^1}  dt \right)  \geq 2 \nu \int_{-\infty}^{0} e^{2 t} \Vert \vu_{\text{e}}(t,\cdot)\Vert^{2}_{\dot{H}^1} dt,
\end{equation*} hence we have
\begin{equation}\label{E3}
\begin{split}
\limsup_{n_k\to +\infty} \left( -2\nu \int_{-t_{n_k}}^{0} e^{2  t} \Vert \vu_{n_k}(t,\cdot) \Vert^{2}_{\dot{H}^1}  dt \right) \leq  & - \liminf_{n_k\to +\infty} \left( 2\nu \int_{-t_{n_k}}^{0} e^{2  t} \Vert \vu_{n_k}(t,\cdot) \Vert^{2}_{\dot{H}^1}  dt \right)   \\
\leq & -2\nu  \int_{-\infty}^{0} e^{2  t} \Vert \vu_{\text{e}}(t,\cdot)\Vert^{2}_{\dot{H}^1} dt.
\end{split}
\end{equation}
Finally, for the fourth term, using again the convergence (\ref{Conv-weak-Attractor}), we obtain 
\begin{equation}\label{E4}
\limsup_{n_k\to +\infty} \left(2 \int_{-t_{n_k}}^{0} e^{2  t} \langle \vf, \vu_{n_k}(t,\cdot)\rangle_{\dot{H}^{-1}\times \dot{H}^1}  dt\right)= 2 \int_{-\infty}^{0}  e^{2   t} \langle \vf, \vu_{\text{e}}(t,\cdot)\rangle_{\dot{H}^{-1}\times \dot{H}^1} \, dt.
\end{equation}
Gathering  estimates  (\ref{E1}), (\ref{E2}), (\ref{E3}), and (\ref{E4}) into estimate (\ref{Estim-energ-limsup}), we obtain
\begin{equation*}
\begin{split}
\limsup_{n_k\to +\infty} \Vert \vu_{n_k}(0,\cdot)\Vert^{2}_{L^2}\leq &\, 2  \int_{-\infty}^{0} e^{2  t}\| \vu_{\text{e}}(t,\cdot)\|^2_{L^2} dt -2\nu  \int_{-\infty}^{0} e^{2   t} \Vert \vu_{\text{e}}(t,\cdot)\Vert^{2}_{\dot{H}^1} dt  \\
&+ 2 \int_{-\infty}^{0}  e^{2 t} \langle \vf, \vu_{\text{e}}(t,\cdot)\rangle_{\dot{H}^{-1}\times \dot{H}^1}=:(A).
\end{split}
\end{equation*}

Now, we will analyze the term $(A)$ above. Since $\vu_{\text{s}}$ is a weak Leray solution of the coupled system (\ref{Alpha-model-eternal}), by performing the same computations as those done in (\ref{Energy-n}) and (\ref{Integral}), we find 
\begin{equation*}
\begin{split}
 \Vert \vu_{\text{e}}(0,\cdot)\Vert^{2}_{L^2} = &\,  2 \int_{-\infty}^{0} e^{2 t} \Vert \vu_{\text{e}}(t,\cdot)\Vert^{2}_{L^2}  dt   -2\nu \int_{-\infty}^{0} e^{2  t} \Vert \vu_{\text{e}}(t,\cdot) \Vert^{2}_{\dot{H}^1}  dt \\
&+ 2\int_{-\infty}^{0} e^{2  t} \langle \vf, \vu_{\text{e}}(t,\cdot)\rangle_{\dot{H}^{-1}\times \dot{H}^1}  dt=(A).
\end{split}
\end{equation*}

Consequently, returning to the previous estimate, we obtain $\ds{\limsup_{n_{k}\to +\infty} \Vert \vu_n(0,\cdot)\Vert^{2}_{L^2}\leq \Vert \vu_{\text{e}}(0,\cdot)\Vert^{2}_{L^2}}$. On the other hand,  by the convergence  (\ref{Conv-weak-Attractor})  we also have the  $\ds{\Vert \vu_{\text{e}}(0,\cdot)\Vert^{2}_{L^2} \leq \liminf_{n_k\to +\infty} \Vert \vu_{n_k}(0,\cdot)\Vert^{2}_{L^2}}$. We then obtain the wished strong convergence: $\ds{\lim_{n_k\to +\infty} \Vert \vu_{n_k}(0,\cdot)\Vert^{2}_{L^2}= \Vert \vu_{\text{e}}(0,\cdot)\Vert^{2}_{L^2}}$. Lemma \ref{Lem-Conv-Attractor} is proven \finpv 

\medskip

We thus finish with the proof of Proposition \ref{Prop-Asymptotically compact}. \finpv 

\medskip

{\bf Point 3}. Continuity. This directly follows from the estimate (\ref{Dependence-date}) with $\vf_1=\vf_2=\vf$, which was proven in  Theorem \ref{Th-Uniqueness}. 

\medskip 

Theorem \ref{Th-Global-Attractor} is now proven. \finpv 

\subsection{Upper bounds of the fractal dimension: proof of Theorem \ref{Th-Fractal-Dimension}}\label{Sec:Fractal-Dimension} 
We shall use the following result, which is stated in the general framework of an abstract   Hilbert space $H$, with its associated norm $\| \cdot\|_H$. We begin by establishing some notation. For $i \in \mathbb{N}$,  recall that  a finite-dimensional orthoprojector $P_i$ is  an orthogonal projection    $P_i:H\to  H_{m_i}$, where $H_{m_i}\subset H$ is a  subspace of finite dimension $\text{dim}(H_{m_i})$, with $m_i\in \mathbb{N}$ and $m_i\geq 1$. Moreover, $\text{dim}(P_i)\leq \text{dim}(H_{m_i})$ represents the rank of $P_i$. 

\begin{Theoreme}[Theorem $2.10$ of \cite{Chueshov}]\label{Th-Tech-Fractal-Dimension} Let $A \subset H$ be a bounded and closed set. Moreover, let $S:A\to H$ be a mapping such that 
\begin{enumerate}
    \item $A \subset S(A)$,
    \item $S$ is Lipschitz continuous: there exists a constant $C_S>0$ such that, for every $h_1, h_2 \in A$, we have $\ds{\| S(h_1)-S(h_2)\|_H \leq C_S \| h_1 - h_2 \|_H}$.
    \item There exists two finite-dimensional orthoprojector $P_1, P_2$, and there exists two constants $0<\kappa_1<1$ and $0<\kappa_2<+\infty$, such that  for any $h_1, h_2 \in A$ the following inequality holds:
    \begin{equation}
        \| S(h_1)-S(h_1) \|_{H} \leq \kappa_1 \| h_1 - h_2 \|_{H} + \kappa_2 \big( \| P_1(h_1-h_1)\|_H + \| P_2\big(S(h_1) - S(h_2)\big)\|_{H} \big).
    \end{equation}
\end{enumerate}
Then, we have
\begin{equation}\label{Dim-Fractal-Tech}
    \text{dim}_{f}(A)\leq \left( \text{dim}(P_1)+\text{dim}(P_2) \right)\ln\left( 1+\frac{8(1+C_S)\sqrt{2}\kappa_2}{1-\kappa_1}\right)\left(\ln\left( \frac{2}{1-\kappa_1} \right)\right)^{-1}.
\end{equation}
\end{Theoreme}

\medskip

{\bf Proof of Theorem \ref{Th-Fractal-Dimension}}. In the setting of  Theorem \ref{Th-Tech-Fractal-Dimension}, one may set $H=L^2_{\text{s}}(\T)$, $A=\mathcal{A}$, where $\mathcal{A}$ is the global attractor  of the semi-group $(S(t))_{t \geq 0}$ defined in (\ref{Semi-group}), and $S=S(t)$ for some  fixed time $t>0$.  However, this choice of the mapping $S$ would yield that the associated constants $C_S, \kappa_1$ and $\kappa_2$ depend on the time $t$. Consequently, by inequality (\ref{Dim-Fractal-Tech})   the quantity  $\text{dim}_f (\mathcal{A})$ would be estimated by a time-dependent expression, which does not seem natural since we aim to estimate $\text{dim}_f (\mathcal{A})$  solely with  the data and parameters of the model (\ref{Alpha-model-Proy}).

\medskip

To overcome this problem, we will follow some of the ideas in \cite{Chueshov}, and we define a framework as follows: for a fixed time $0<T<+\infty$, which we will specify later in expression (\ref{T-eta}) below, we consider the Hilbert space
\begin{equation*}
H = L^2_{\text{s}} (\T) \times L^2([0,T],L^2_{\text{s}}(\T)),
\end{equation*}
with the usual norm
\begin{equation*}
\| (\vu_0, \vu) \|^2_{H} = \| \vu_0 \|^2_{L^2} + \int_{0}^{T} \| \vu(t, \cdot) \|^2_{L^2}  dt.
\end{equation*}
Then, with the global attractor $\mathcal{A}\subset L^2_{\text{s}}(\T)$, we define the set
\begin{equation}\label{AT}
A = \Big\{ U = (\vu_0, \vu) \in H : \ \vu_0 \in \mathcal{A}, \ \ \vu = \vu(t, \cdot) = S(t)\vu_0, \ \ t \in [0,T] \Big\},
\end{equation}
where $\vu(t,\cdot)=S(t)\vu_0$ is the unique weak Leray solution of the system (\ref{Alpha-model-Proy}) arising from $\vu_0$. Recall that $\vu\in L^\infty([0,T],L^2_{\text{s}}(\T))$, and hence $\vu\in L^2([0,T],L^2_{\text{s}}(\T))$. Moreover, since $\mathcal{A}$ is compact in $L^2_{\text{s}}(\T)$ and $\vu(t,\cdot)\in \mathcal{A}$  for every $t\geq 0$, the set $A$ is bounded and closed in $H$. We then define the mapping
\begin{equation}\label{S}
S: A \to H, \quad U = (\vu_0, \vu) \mapsto S(U) = \big(\vu(T, \cdot), \vu(T + t)\big), \quad \text{with} \quad t \in [0,T].
\end{equation}
In other words, the mapping $S$ is simply a translation to the time $T$ of the pair $U=(\vu_0,\vu)$.

\medskip

Once we have established the framework, we will verify each point stated in Theorem \ref{Th-Tech-Fractal-Dimension}. For simplicity, we shall omit generic constants that appear in our estimates in order to focus on explicit expressions involving the relevant data and parameters of our model.

\medskip

{\bf Point 1}.  This fact directly follows from the second point of Definition \ref{Def-Global-Attractor}, and from the definitions of the mapping $S$ and the set $A$ given in (\ref{S}) and (\ref{AT}), respectively. We thus have $S(A)=A$.  

\medskip

{\bf Point 2}. Lipschitz continuity of the map $S$ essentially follows from estimates provided in the proof of   Theorem \ref{Th-Existence}:  let $\vu_{0,1}$ and $\vu_{0,2}$ be two initial data in $\mathcal{A}$. For $0<t<T$, let $\vu_1(t,\cdot), \vu_2(t,\cdot) \in L^2_{\text{s}}(\T)$ be the associated weak Leray solutions of the coupled system (\ref{Alpha-model-Proy}). We define $\vw(t,\cdot)=\vu_1(t,\cdot)-\vu_2(t,\cdot)$, which solves the problem  (\ref{Alpha-model-Difference}) with $\vf_1 = \vf_2=\vf$, hence $\vf_1-\vf_2=0$. 

\medskip

By the energy balance given in (\ref{Iden-Cont-Dep-00}), and  by estimate (\ref{Estim-Non-Lin-Uniq}), where the expression $g(\vu_1,\vu_2)$ is defined in (\ref{g}),  using the Gr\"onwall inequality we get
\begin{equation}\label{Gronwall}
\| \vu_1(t,\cdot)-\vu_2(t,\cdot)\|^2_{L^2} \leq \exp\left( \int_{0}^{t}g\big(\vu_1(s,\cdot),\vu_2(s,\cdot)\big)ds\right)\| \vu_{0,1}-\vu_{2,0}\|^{2}_{L^2}.
\end{equation}
By expression (\ref{Control-Integral}), we have
\begin{equation}\label{Control-Integral-2}
 \int_{0}^{t}g\big(\vu_1(s,\cdot),\vu_2(s,\cdot)\big)ds \leq  C \, K_3 \sup_{0\leq s \leq t} \Big(\| \vu_1(s,\cdot) \|^2_{L^2}+\| \vu_2(s,\cdot) \|^2_{L^2}+1\Big)  \left(\nu\int_{0}^{t}  \| \vu_2(s,\cdot) \|^{2}_{\dot{H}^1} ds \right)^{1/2} \frac{t^{1/2}}{\nu^{1/2}}.   
\end{equation}
To control the first term on the right-hand side, by estimate (\ref{Energ-Estim-Attractor-1}) and the fact that $\vu_{0,1},\vu_{0,2}\in \mathcal{A}\subset \mathcal{B}$ (where the absorbing set $\mathcal{B}$ is defined in (\ref{Absorbing-Set})), we obtain
\begin{equation*}
    \sup_{0\leq s \leq t} \Big(\| \vu_1(s,\cdot) \|^2_{L^2}+\| \vu_2(s,\cdot) \|^2_{L^2}+1\Big) \leq 3 R^2+1.
\end{equation*}
On the other hand, to control the second term on the right-hand side, we use estimate (\ref{Control-Energ-2}) to write
\begin{equation}\label{Control-Norm-H1}
\left(\nu\int_{0}^{t}  \| \vu_2(s,\cdot) \|^{2}_{\dot{H}^1} ds \right)^{1/2}  \leq \| \vu_{0,2}\|_{L^2} +\frac{t^{1/2}}{\nu^{1/2}}\| \vf \|_{\dot{H}^{-1}} \leq R+ \frac{T^{1/2}}{\nu^{1/2}}\| \vf \|_{\dot{H}^{-1}}.
\end{equation}
Gathering these controls into estimate (\ref{Control-Integral-2}), we have
\begin{equation*}
 \int_{0}^{t}g\big(\vu_1(s,\cdot),\vu_2(s,\cdot)\big)ds \leq  C \, K_3 (3R^2+1)\left( R+ \frac{T^{1/2}}{\nu^{1/2}}\| \vf \|_{\dot{H}^{-1}} \right)\frac{t^{1/2}}{\nu^{1/2}}.      
\end{equation*}
Setting 
\begin{equation}\label{C1}
 \mathfrak{C}_1:= \frac{C\, K_3}{\nu^{1/2}} (3R^2+1)\left( R+ \frac{T^{1/2}}{\nu^{1/2}}\| \vf \|_{\dot{H}^{-1}} \right),   
\end{equation}
 we  come back to inequality  (\ref{Gronwall}), and we get 
\begin{equation*}
    \|\vu_1(t,\cdot)-\vu_2(t,\cdot)\|_{L^2} \leq e^{\mathfrak{C}_1\, t^{1/2}}\| \vu_{0,1}-\vu_{0,2}\|_{L^2}.
\end{equation*}

With this estimate at hand, for $U_1=(\vu_{0,1}, \vu_1),  \ U_2=(\vu_{0,2}, \vu_2)\in A$ we write
\begin{equation}\label{Lipschitz}
    \begin{split}
  \| S(U_1)-S(U_2)\|^2_{H} =&\, \| \vu_1(T,\cdot)-\vu_2(T,\cdot)\|^2_{L^2}+\int_{T}^{2T}\| \vu_1(t,\cdot)-\vu_2(t,\cdot)\|^2_{L^2} dt   \\
  \leq &\, e^{\mathfrak{C}_1\, T^{1/2}}\| \vu_{0,1}-\vu_{0,2}\|^2_{L^2}+ \left( \int_{T}^{2T} e^{\mathfrak{C}_1\, t^{1/2}}  dt \right) \| \vu_{0,1}-\vu_{0,2}\|^2_{L^2}\\
\leq &\,   e^{\mathfrak{C}_1\, T^{1/2}}(1+T) \| \vu_{0,1}-\vu_{0,2}\|^2_{L^2} \\
\leq &\, e^{\mathfrak{C}_1\, T^{1/2}}(1+T) \| U_1 - U_2 \|^2_{H}, \quad C_S:= e^{\mathfrak{C}_1\, T^{1/2}}(1+T).
    \end{split}
\end{equation}

{\bf Point 3}. Our starting point is to prove the following technical estimate:
\begin{Proposition}\label{Prop-Estim-Tech} Let $\vu_{0,1}, \vu_{0,2}\in \mathcal{A}$ be two initial data. For $0<t<T$, let $\vu_1(t,\cdot), \vu_2(t,\cdot) \in L^2_{\text{s}}(\T)$ be the associated weak  Leray solutions of the coupled system (\ref{Alpha-model-Proy}). 
\medskip

Let $\eta>0$ defined in (\ref{eta}), $\mathfrak{C}_0>0$   given in estimate  (\ref{Dependence-date}), $K_3>0$, $R>0$ given in (\ref{K3}) and (\ref{R}) respectively. Define the quantity, 
\begin{equation}\label{K4}
K_4=\frac{C}{\nu}K_3(2R+1)^4 R^2.    
\end{equation}
Moreover,  for $m\in \mathbb{N}$, and for $\ds{\vg(x)=\sum_{k \in \Z}e^{i k\cdot x}\widehat{g}_k}$, define the finite-dimensional orthoprojector $ P_m(\vg)= \ds{\sum_{k \in \Z,\ |k|\leq m}}e^{i k\cdot x}\widehat{g}_k$. 

\medskip

For any $0<\varepsilon<1$, there exists $m=m(\varepsilon)\in \mathbb{N}^{*}$ such that it  holds:
\begin{equation}\label{Ineq-Tech-3}
    \begin{split}
        \| \vu_1(t,\cdot)-\vu_2(t,\cdot)\|^2_{L^2} \leq &\, \left( e^{-\eta t} + \varepsilon K_4 \frac{e^{\mathfrak{C}_0 t^{1/2}}}{\nu}\right) \| \vu_{0,1}-\vu_{0,2}\|^2_{L^2}\\
        &\,+ K_4 \, \int_{0}^{t} \|P_m(\vu_1(\tau,\cdot)-\vu_2(\tau,\cdot))\|^2_{L^2} d \tau.
    \end{split}
\end{equation}  
\end{Proposition} 
\pv  As before, we define  $\vw(t,\cdot)=\vu_1(t,\cdot)-\vu_2(t,\cdot)$, which solves the problem  (\ref{Alpha-model-Difference}) with $\vf_1-\vf_2=0$. We write 
\begin{equation*}
        \frac{d}{dt}\| \vw(t,\cdot)\|^2_{L^2}+2\nu \| \vw(t,\cdot)\|^2_{\dot{H}^1}= - 2\int_{\T} \big((\vu_{1,\alpha}-\vu_{2,\alpha}) \vec{\nabla}\big) \vu_2 \cdot \vw \, dx,
\end{equation*}
where we need to  estimate the term on the right. As $\vu_{1,\alpha}$ and $\vu_{2,\alpha}$ are divergence-free vector fields, integrating by parts and applying the Cauchy-Schwarz inequality, we obtain
\begin{equation*}
\begin{split}
&- 2\int_{\T} \big((\vu_{1,\alpha}-\vu_{2,\alpha}) \vec{\nabla}\big) \vu_2 \cdot \vw dx =\, - 2\int_{\T}\text{div}\big( (\vu_{1,\alpha}-\vu_{2,\alpha})  \otimes  \vu_2 \big)\cdot \vw dx \\
=&\, -2 \sum_{i,j=1}^{3} \int_{\T}\partial_j \big((u_{1,j,\alpha}-u_{2,j,\alpha})u_{2,i}\big) w_i dx = 2\sum_{i,j=1}^{3}\int_{\T} \big((u_{1,j,\alpha}-u_{2,i,\alpha})u_{2,i}\big) \partial_j w_i dx\\
\leq &\, 2 \| (\vu_{1,\alpha}-\vu_{2,\alpha})\otimes \vu_2 \|_{L^2}\| \vw \|_{\dot{H}^1} \leq \frac{1}{\nu}\| (\vu_{1,\alpha}-\vu_{2,\alpha})\otimes \vu_2 \|^2_{L^2} + \nu \| \vw \|^2_{\dot{H}^1}. 
\end{split}
\end{equation*}
It yields
\begin{equation*}
    \frac{d}{dt}\| \vw(t,\cdot)\|^2_{L^2} + \nu \| \vw(t,\cdot)\|^2_{\dot{H}^{1}} \leq \frac{1}{\nu} \| (\vu_{1,\alpha}-\vu_{2,\alpha})\otimes \vu_2 \|^2_{L^2}.
\end{equation*}
Using the Poincaré's inequality, with $\eta>0$ defined in (\ref{eta}), we obtain
\begin{equation*}
    \frac{d}{dt}\| \vw(t,\cdot)\|^2_{L^2} +\eta \| \vw(t,\cdot)\|^2_{L^2} \leq \frac{1}{\nu} \| (\vu_{1,\alpha}-\vu_{2,\alpha})\otimes \vu_2 \|^2_{L^2}.
\end{equation*}
Then,  by Proposition \ref{Cont-Dep-H2-Filter}, and the fact that for $i=1,2$ we have $\| \vu_i(t,\cdot)\|_{L^2}\leq R$,  we can write
\begin{equation}\label{Ineq-Tech-1}
\begin{split}
   &\,  \frac{1}{\nu} \| (\vu_{1,\alpha}-\vu_{2,\alpha})\otimes \vu_2 \|^2_{L^2} \leq \, \frac{1}{\nu}\| \vu_{1,\alpha}-\vu_{2} \|^2_{L^\infty}\, \| \vu_{2} \|^2_{L^2}
    \leq \, \frac{C}{\nu} \| \vu_{1,\alpha}-\vu_{2,\alpha} \|^2_{H^2}\, \| \vu_{2} \|^2_{L^2}\\
    \leq &\, \frac{C}{\nu} K_3 \left(\| \vu_1(t,\cdot)\|_{L^2}+\| \vu_2(t,\cdot)\|_{L^2}+1\right)^4 \| \vw(t,\cdot)\|^2_{L^2}\, \| \vu_2(t,\cdot)\|^2_{L^2}\\
    \leq &\, \frac{C}{\nu} K_3 (2R+1)^4 R^2\, \| \vw(t,\cdot)\|^2_{L^2}=:\, K_4 \, \| \vw(t,\cdot)\|^2_{L^2}.
\end{split}
\end{equation}
Now, for $0<\varepsilon<1$, we  shall prove that there exists $m=m(\varepsilon)\in \mathbb{N}^{*}$ such that the following estimate holds: 
\begin{equation}\label{Ineq-Tech-2}
    \| \vw(t,\cdot)\|^2_{L^2} \leq \varepsilon \| \vw(t,\cdot)\|^{2}_{\dot{H}^1}+ \|P_m ( \vw(t,\cdot))\|^{2}_{L^2}.
\end{equation}
Indeed,  we  write 
\begin{equation*}
        \| \vw\|^2_{L^2}=   L^3\sum_{k \in \Z} |\widehat{w}_k|^2= L^3\sum_{k \in \Z,\ |k|\leq m}  |\widehat{w}_k| ^2+L^3\sum_{k \in \Z,\ |k|>m}  |\widehat{w}_k|^2.
\end{equation*}
For the second term, we obtain
\begin{equation*}
 L^3\sum_{k \in \Z,\ |k|>m}  |\widehat{w}_k|^2=L^3\sum_{k \in \Z,\ |k|>m} |k|^2|k|^{-2} |\widehat{w}_k|^2\leq L^3 m^{-2} \sum_{k \in \Z,\ |k|>m} |k|^2 |\widehat{w}_k|^2.   
\end{equation*}
We have $m^{-2} < \varepsilon$ as long as $m> \varepsilon^{-1/2}$. We thus set $m=m(\varepsilon)\geq [\varepsilon^{-1/2}]+1$, where $[\varepsilon^{-1/2}]$ is the smallest integer greater than or equal to  $\varepsilon^{-1/2}$. Gathering these estimates, we obtain the desired inequality (\ref{Ineq-Tech-2}). 

\medskip

Returning to inequality (\ref{Ineq-Tech-1}), we have
\begin{equation*}
 \frac{1}{\nu} \| (\vu_{\alpha,1}-\vu_{\alpha,2})\otimes \vu_2 \|^2_{L^2} \leq \varepsilon K_4 \| \vw(t,\cdot)\|^2_{\dot{H}^1}+ K_4 \| P_m(\vw(t,\cdot))\|^{2}_{L^2},   
\end{equation*}
hence
\begin{equation*}
    \frac{d}{dt}\| \vw(t,\cdot)\|^2_{L^2}+\nu \| \vw(t,\cdot)\|^2_{L^2} \leq \varepsilon K_4 \| \vw(t,\cdot)\|^2_{\dot{H}^1}+ K_4 m^{2s}\| P_m(\vw(t,\cdot))\|^{2}_{L^2}.
\end{equation*}
Applying the Gr\"onwall inequality, we get
\begin{equation*}
    \| \vw(t,\cdot)\|^2_{L^2} \leq e^{-\eta t}\| \vw(0,\cdot)\|^2_{L^2}+ \varepsilon K_4 \int_{0}^{t} \| \vw(\tau,\cdot)\|^{2}_{\dot{H}^1} d \tau + K_4  \int_{0}^{t}\| P_m(\vw(\tau,\cdot))\|^2_{L^2} d \tau.
\end{equation*}
Finally, by estimate (\ref{Dependence-date}) (recall that $\vf_1-\vf_2=0$) we have $\ds{\int_{0}^{t} \| \vw(\tau,\cdot)\|^{2}_{\dot{H}^1} d \tau \leq \frac{e^{\mathfrak{C}_0 t^{1/2}}}{\nu}\| \vw(0,\cdot)\|^2_{L^2}}$, hence we  obtain the wished estimate (\ref{Ineq-Tech-3}).  Proposition \ref{Prop-Estim-Tech} is proven. \finpv 

\medskip

In estimate (\ref{Ineq-Tech-3}), setting   $t=T$ we obtain
\begin{equation*}
    \begin{split}
        \| \vu_1(T,\cdot)-\vu_2(T,\cdot)\|^2_{L^2} \leq &\, \left( e^{-\eta T} + \varepsilon K_4 \frac{e^{\mathfrak{C}_0 T^{1/2}}}{\nu}\right) \| \vu_{0,1}-\vu_{0,2}\|^2_{L^2}+ K_4 \, \int_{0}^{T} \| P_m(\vu_1(\tau,\cdot)-\vu_2(\tau,\cdot))\|^2_{L^2} d \tau.
    \end{split}
\end{equation*}
On the other hand, we integrate estimate (\ref{Ineq-Tech-3}) over the interval $[T,2T]$, and estimating each term on the right, we get
\begin{equation*}
    \begin{split}
 \int_{T}^{2T} \| \vu_1(t,\cdot)-\vu_2(t,\cdot)\|^2_{L^2}  dt
    \leq &\, \left( e^{-\eta T} T+ \varepsilon K_4 \frac{e^{\mathfrak{C}_0 T^{1/2}}}{\nu} T\right) \| \vu_{0,1}-\vu_{0,2}\|^2_{L^2}\\
    &\,+ K_4  T\, \int_{0}^{T} \|P_m(\vu_1(\tau,\cdot)-\vu_2(\tau,\cdot))\|^2_{L^2} d \tau.
    \end{split}
\end{equation*}
Gathering these estimates, for $U_1=(\vu_{0,1},\vu_1),  \ U_2=(\vu_{0,2},\vu_2)\in A$, we have
\begin{equation*}
    \begin{split}
 &\| S(U_1)-S(U_2)\|^2_{H} \\
 =&\,  \| \vu_1(T,\cdot)-\vu_2(T,\cdot)\|^2_{L^2}+ \int_{T}^{2T} \| \vu_1(t,\cdot)-\vu_2(t,\cdot)\|^2_{L^2}  dt\\
 \leq &\,  \left( e^{-\eta T} + \varepsilon K_4 \frac{e^{\mathfrak{C}_0 T^{1/2}}}{\nu}\right)(1+T) \| \vu_{0,1}-\vu_{0,2}\|^2_{L^2}\\
 &\, + K_4 (1+T) \int_{O}^{T} \| P_m(\vu_1(\tau,\cdot)-\vu_2(\tau,\cdot)) \|^{2}_{L^2}dt\\
 \leq &\,  \left( e^{-\eta T} + \varepsilon K_4 \frac{e^{\mathfrak{C}_0 T^{1/2}}}{\nu}\right)(1+T) \left( \| \vu_{0,1}-\vu_{0,2}\|^2_{L^2}+\int_{0}^{T}\| \vu_1(t,\cdot)-\vu_{2}(t,\cdot)\|^2_{L^2}dt \right)\\
 &\, + K_4 (1+T) \left( \| P_m(\vu_{0,1}-\vu_{0,2})\|^{L^2}+ \int_{0}^{T}   \| P_m(\vu_1(\tau,\cdot)-\vu_2(\tau,\cdot)) \|^{2}_{L^2}dt  \right)\\
 =&\, \left( e^{-\eta T} + \varepsilon K_4 \frac{e^{\mathfrak{C}_0 T^{1/2}}}{\nu}\right)(1+T) \| U_1-U_2\|^{2}_{H}+ K_4 (1+T) \| P_m(U_1-U_2)\|^2_{H}.
    \end{split}
\end{equation*}
In this estimate, we can choose $T$ large enough so that $\ds{e^{-\eta T}(1+T)<\frac{1}{4}}$. In fact, this inequality is equivalent to $\ds{\ln(4(1+T))<\eta T}$. Moreover, $\ds{\ln(4(1+T))<\eta T}$,  
it is sufficient to verify the inequality $\ds{2\sqrt{1+T}<\eta T}$, which holds as long as $T>\frac{2}{\eta^2}(1+\sqrt{1+\eta^2})$.  We thus set
\begin{equation}\label{T-eta}
T_\eta := \frac{4}{\eta^2}(1 + \sqrt{1 + \eta^2}).
\end{equation}
With the time $T_\eta$,  we  set the parameter $\varepsilon$ small enough such that $\ds{ \varepsilon K_4 \frac{e^{\mathfrak{C}_0 T^{1/2}_\eta}}{\nu}(1+T_\eta)<\frac{1}{4}}$. We thus obtain $\left( e^{-\eta T_\eta} + \varepsilon K_4 \frac{e^{\mathfrak{C}_0 T^{1/2}_\eta}}{\nu}\right)(1+T_\eta)<\frac{1}{2}$, hence we write
\begin{equation}\label{Estim-Point-3}
    \begin{split}
 \| S(U_1)-S(U_2)\|^2_{H} \leq \frac{1}{2} \| U_1 - U_2\|^2_{H} + K_4 (1+T_\eta)  \| P_m(U_1 -U_2) \|^{2}_{H}.
    \end{split}
\end{equation}
At this point, recall that we have $m=m(\varepsilon)\geq [\varepsilon^{-1/2}]+1$. Moreover, note that $\ds{ \varepsilon K_4 \frac{e^{\mathfrak{C}_0 T^{1/2}_\eta}}{\nu}(1+T_\eta)<\frac{1}{4}}$ as long as $\ds{\varepsilon^{-1/2}>\frac{2}{\sqrt{\nu}}\sqrt{K_4(1+T_\eta)}e^{\frac{\mathfrak{C}_0}{2} T^{1/2}_{\eta}}}$. We thus set $\ds{\varepsilon^{-1/2}:=\frac{4}{\sqrt{\nu}}\sqrt{K_4(1+T_\eta)}e^{\frac{\mathfrak{C}_0}{2} T^{1/2}_{\eta}}}$, and we set 
\begin{equation}\label{m} 
 m:=\left[\frac{4}{\sqrt{\nu}}\sqrt{K_4(1+T_\eta)}e^{\frac{\mathfrak{C}_0}{2} T^{1/2}_{\eta}} \right]+1.   
\end{equation}

In this manner, we can apply Theorem \ref{Th-Tech-Fractal-Dimension} to conclude that the set $A$ (defined in expression (\ref{AT})) is compact in the strong topology of the space $H$. Moreover, by inequality (\ref{Dim-Fractal-Tech}), its fractal dimension is bounded by
\begin{equation*}
 \text{dim}_{f}(A)\leq  \text{dim}(P_m) \ln\left(1+\frac{8(1+C_S)\sqrt{2}\kappa_2}{1-\kappa_1}\right)\left( \ln\left( \frac{2}{1-\kappa_1} \right)\right)^{-1},  
\end{equation*}
where, by estimates (\ref{Lipschitz}) and (\ref{Estim-Point-3}), we have
\begin{equation*}
    C_S= e^{\mathfrak{C}_1 T^{1/2}_{\eta}}(1+T_\eta), \quad \kappa_1=\frac{1}{2}, \quad \kappa_2= K_4(1+T_\eta).
\end{equation*}
In addition, by \cite[Lemma 8]{Chai1} and since $m\geq 1$ (see expression (\ref{m})) we get
\begin{equation*}
\text{dim}(P_m) \leq 8\left( 4m^3+6m^2+8m+3 \right) \leq 256 m^3.   \end{equation*}

Gathering these expressions and estimates, we write
\begin{equation}\label{Estim-Tech-Fractal-Dim}
\begin{split}
    \text{dim}_{f}(A)\leq &\,  \frac{256}{\ln(4)}  \left( \left[\frac{4}{\sqrt{\nu}}\sqrt{K_4(1+T_\eta)}e^{\frac{\mathfrak{C}_0}{2} T^{1/2}_{\eta}} \right]+1\right)^3\\
    &\, \times  \ln\left( 1+16\sqrt{2}K_4\left((1+e^{\mathfrak{C}_1 T^{1/2}_{\eta}}(1+T_\eta)\right)(1+T_\eta)\right)=:D.
    \end{split}
\end{equation}

Finally, we define the operator $\mathcal{Q}:H \to L^2_{\text{s}}(\T)$ by  $U=(\vu_0, \vu)\mapsto \mathcal{Q}(U)=\vu_0$. Since $\mathcal{Q}$ is Lipschitz continuous and $\mathcal{A}=\mathcal{Q}(A)$, we have $\ds{\text{dim}_{f}(\mathcal{A})\leq  \text{dim}_f(A)}$.

\medskip

To complete the proof, we must estimate the right-hand side of (\ref{Estim-Tech-Fractal-Dim}). Recall that we will focus only on quantities depending on the parameters $\alpha>0$, $0<\beta\leq 1$ and the norm $\| \vf \|_{\dot{H}^{-1}}$. Consequently, we set $L=1$, $\nu=1$. Therefore, from expression  expression (\ref{eta}), we have $\eta \simeq 1$, and by expression (\ref{T-eta}), we obtain $T_\eta \simeq 1$. This yields
\begin{equation*}
\begin{split}
    D \lesssim  &\, \left(1+\sqrt{K_4}\,e^{\mathfrak{C}_0}\right)^3 \ln\left( 1+ K_4(1+e^{\mathfrak{C}_1}) \right)  \\
    \leq &\,  \left(1+\sqrt{K_4} \, e^{\mathfrak{C}_0}\right)^3 \ln\left( 1+ K_4(2e^{\mathfrak{C}_1}) \right)\\
    \lesssim  &\, \left(1+\sqrt{K_4} \, e^{\mathfrak{C}_0}\right)^3 \ln\left( 1+ K_4 \, e^{\mathfrak{C}_1} \right).
    \end{split}
\end{equation*}
Then, we have
\begin{Lemme}\label{Lem-Estim-K-alpha-beta} For $\alpha >0$ and $0<\beta < 1$, let $K(\alpha,\beta)$ be the quantity defined in expression (\ref{K-alpha-beta}). The following estimates hold:
\begin{enumerate}
    \item For the quantity $K_4$ given in (\ref{K4}) we have $K_4 \lesssim K(\alpha,\beta)(1+\| \vf \|_{\dot{H}^{-1}})^6$.
    \item For the quantity $\mathfrak{C}_0$ defined in (\ref{C0-Uniq}) we have $\mathfrak{C}_0 \lesssim K(\alpha,\beta)(1+\| \vf \|_{\dot{H}^{-1}})^{3/2}$.
    \item  For the quantity $\mathfrak{C}_1$ given in (\ref{C1}) we have $\mathfrak{C}_1 \lesssim K(\alpha,\beta)(1+\| \vf \|_{\dot{H}^{-1}})^3$.
    \end{enumerate} 
\end{Lemme}
\pv We come back to the quantity $K_3$ defined in expression (\ref{K3}), and we shall prove that $K_3 \lesssim K(\alpha,\beta)$. Since $K_3$ is defined by the quantities $K_0$, $K_1$, $K_2$, $L_1$, $L_2$ and $L_3$ given in  (\ref{Control-Filter}), (\ref{Dependence-data}),  (\ref{Filter-Control-H2}), (\ref{L1}), (\ref{L2}) and (\ref{L3}) respectively, we must estimate each of them. To do this, let consider the following cases of the parameter $\alpha$.

\begin{itemize}
    \item When $0<\alpha\leq 1$. Since $0<  \beta < 1$,  we have $K_0 \simeq \frac{1}{\alpha^2 \beta}\geq 1$. Moreover, by expressions (\ref{L1}) and (\ref{L2}), we have $L_1\simeq \alpha^2 \leq  1$ and $L_2 \simeq \alpha^2 \leq 1$. Thereafter,  by expressions (\ref{Control-Filter})  and (\ref{Filter-Control-H2}) we obtain $K_1 \lesssim K^{3/2}_{0}$ and $K_2 \lesssim  K^{3/2}_{0}$. Consequently,  getting back to expression (\ref{K3}), and since by expression (\ref{L3}) we have $L_3 \simeq \alpha^2 \leq 1$, we obtain $K_3 \lesssim K^{5/2}_{0} \simeq \frac{1}{\alpha^5 \beta^{5/2}}$. 

    \item When $1<\alpha \leq \frac{1}{\sqrt{\beta}}$. In this case we also have $K_0 \simeq \frac{1}{\alpha^2 \beta}\geq 1$. Moreover, since $L_1 \simeq L_2 \simeq \alpha^2$, from expressions (\ref{Control-Filter})  and (\ref{Filter-Control-H2}) we get $K_1 \lesssim \alpha^2 K^{3/2}_{0}$ and $K_2 \lesssim \alpha^2 K^{3/2}_{0}$. Thereafter, since $L_3 \simeq \alpha^2$  we obtain
$K_3 \lesssim \alpha^4 K^{5/2}_{0}  \simeq \frac{1}{\alpha \beta^{5/2}}$.  

\item   When $\frac{1}{\sqrt{\beta}}<\alpha$. Here we have $K_0 \simeq 1$. As before, we have $L_1 \simeq L_2 \simeq \alpha^2 $, hence we obtain  $K_1 \simeq \alpha^2$ and  $K_2 \simeq \alpha^2$. Then, since $L_3 \simeq \alpha^2$,  we have  $K_3 \simeq \alpha^4$.
\end{itemize}

\medskip

With inequality   $K_3 \lesssim K(\alpha,\beta)$ at our disposal, and  recalling that from expression (\ref{R}) we can write $R\simeq \| \vf \|_{\dot{H}^{-1}}$, we obtain  the following estimates:  for  $K_4$ given in (\ref{K4}), we have $K_4 \lesssim K_3(1+R)^6 \lesssim K(\alpha,\beta)(1+\| \vf \|_{\dot{H}^{-1}})^6$. Then,   for  $\mathfrak{C}_0$ defined in (\ref{C0-Uniq}), we have $\mathfrak{C}_0 \lesssim K_3(1+R)^{3/2}\lesssim  K(\alpha,\beta)(1+\| \vf \|_{\dot{H}^{-1}})^{3/2}$. Finally, for  $\mathfrak{C}_1$ given in (\ref{C1}), we obtain $\mathfrak{C}_1 \lesssim K_3 (1+R)^3  \lesssim K(\alpha,\beta)(1+\| \vf \|_{\dot{H}^{-1}})^3$. Lemme \ref{Lem-Estim-K-alpha-beta} is proven.  \finpv

\medskip

Gathering these estimates in the last control: $D \lesssim \left(1+\sqrt{K_4}\, e^{\mathfrak{C}_0}\right)^3 \ln\left( 1+ K_4\, e^{\mathfrak{C}_1}\right)$, we finally obtain the wished estimate (\ref{Estimate-fractal-dimension}).  Theorem \ref{Th-Fractal-Dimension}  is now proven. \finpv 
\appendix
\section{The non-regularized nonlinear filter equation}\label{Appendix}
As explained in the introduction, we prove here that additional hypotheses on the indicator function $A(\cdot)$ yield the existence of a weak Leray solution to the Navier-Stokes-alpha model with a non-regularized filter equation:
\begin{equation}\label{Alpha-model-Appendix}
    \begin{cases}\vspace{2mm}
    \partial_t \vu - \nu\Delta \vu + \P \Big((\vu_\alpha \cdot \vec{\nabla})\vu \Big)   =  \vf, \quad \text{div}(\vu)=0, \\ \vspace{2mm}
     -\alpha^2 \P \Big( \text{div} \left( A (\vu) (\vec{\nabla} \otimes \vu_\alpha)^T \right) \Big)+\vu_\alpha = \vu, \quad \text{div}(\vu_\alpha)=0,  \\ 
\vu(0,\cdot)=\vu_0, \quad \text{div}(\vu_0)=0. 
    \end{cases}
\end{equation}
For simplicity, we have eliminated the pressure terms $P$ and $P_\alpha$ by applying Leray's projector $\P$ to each equation in this system.
\begin{Theoreme}\label{Th-Appendix} Let $\vu_0 \in L^2_{\text{s}}(\T)$ be an initial datum, and let $\vf \in L^{2}_{loc}([0,+\infty[, \dot{H}^{-1}(\T))$ be a divergence-free and zero-mean  external force. In the nonlinear filter equation:
\[      -\alpha^2 \P \Big(\text{div} \left( A (\vu) (\vec{\nabla} \otimes \vu_\alpha)^T \right)\Big)+\vu_\alpha  = \vu, \quad \text{div}(\vu_\alpha)=0, \quad \alpha>0, \]
for a fixed  constant $0<\beta<1$,   assume that the indicator function $A(\cdot)$  satisfies:
	\begin{equation}\label{Conditions-a-existence-1-appendix}
\beta \leq A(x) \leq 1, \qquad \mbox{for all $x\in \Rt$}, 
\end{equation}
and that 
\begin{equation}\label{Conditions-a-existence-2-appendix}
    A(\cdot)\in \mathcal{C}^{0,1}(\Rt). 
\end{equation}
Additionally, assume  that $A(\cdot)$ satisfies  the stronger Lipschitz-type continuity condition: 
\begin{equation}\label{Conditions-a-existence-3-appendix}
\| A(\vg_1)-A(\vg_2)\|_{L^\infty}\leq C \| \vg_1 - \vg_2\|_{L^2}, \quad \mbox{for any} \quad \vg_1, \vg_2 \in L^2(\T). 
\end{equation}
Then,  the coupled system (\ref{Alpha-model-Appendix}) has a  weak Leray  solution  $(\vu, \vu_\alpha)$ in the sense of Definition \ref{Def-Leray-Solutions}.  
\end{Theoreme}

Note that conditions (\ref{Conditions-a-existence-1-appendix}) and (\ref{Conditions-a-existence-2-appendix}) are the same as those in Theorem \ref{Th-Existence}. Conversely, the stronger Lipschitz-type continuity condition (\ref{Conditions-a-existence-3-appendix}) is a technical requirement to get rid of the convolution term $\varphi$. This Lipschitz-type continuity condition is stronger than the classical one stated in (\ref{Conditions-a-existence-2-appendix}) in the sense that the $L^\infty$-norm of $A(\vg_1)-A(\vg_2)$ and  is controlled by the $L^2$-norm of $\vg_1-\vg_2$. 
 
\medskip

\pv Coming back to the coupled system with regularized nonlinear filter equation (\ref{Alpha-model}), for fixed $\varepsilon>0$,  we set $\varphi=\varphi_\varepsilon$, where $(\varphi_\varepsilon)_{\varepsilon>0}$ is a family  of approximation of the identity defined over $\T$. See, for instance, \cite[Chapter $1.2$]{Grafakos}. Assuming (\ref{Conditions-a-existence-1-appendix}), (\ref{Conditions-a-existence-2-appendix}), and since we have $\varphi_\varepsilon \in L^2(\T)$, by Theorem \ref{Th-Existence} there exists $(\vu_\varepsilon, \vu_{\varepsilon,\alpha}, P_\varepsilon, P_{\varepsilon,\alpha})$ a weak Leray solution to (\ref{Alpha-model}). We will prove that this family converges to a Leray's solution $(\vu, \vu_\alpha,P, P_\alpha)$ to the coupled system (\ref{Alpha-model-Appendix}).

\medskip

For every $\varepsilon>0$, the velocity $\vu_\varepsilon$ verifies the energy equality (\ref{Energy-Control}), hence, for any time $0<T<+\infty$, we have the uniform control 
\begin{equation}\label{Unif-Control1-Appendix}
    \sup_{0\leq t \leq T}\|\vu_\varepsilon(t,\cdot)\|^2_{L^2}+\nu \int_{0}^{T}\| \vu_\varepsilon(t,\cdot)\|^2_{\dot{H}^1} dt \leq \|\vu_0\|^2_{L^1}+\frac{1}{\nu}\int_{0}^{T}\| \vf(t,\cdot)\|^2_{\dot{H}^{-1}}dt. 
\end{equation}
Moreover, we also have the uniform control:
\begin{equation}\label{Unif-Control2-Appendix}
 \| \partial_t \vu_\varepsilon \|_{L^2([0,T],H^{-1}(\T))} \leq M,
\end{equation}
with a quantity $0<M=M(\vu_0,\vf,\nu,T)<+\infty$, which does not depend of $\varepsilon$. Indeed, we write
\begin{equation*}
    \partial_t \vu_\varepsilon = \nu \Delta \vu_\varepsilon - \P \big( (\vu_{\varepsilon,\alpha}\cdot \vec{\nabla}) \vu_\varepsilon \big)+\vf.
\end{equation*}
For the first term on the right, by (\ref{Unif-Control1-Appendix}) we obtain $\ds{\| \Delta \vu_\varepsilon \|_{L^2_t \dot{H}^{-1}_x}\leq  \| \vu_\varepsilon \|_{L^2_t H^{1}_x}\leq M}$. For the second term on the right, since $\text{div}(\vu_{\varepsilon,\alpha})=0$,  using H\"older inequalities,  Sobolev inequalities and interpolation inequalities, we obtain
\begin{equation*}
\begin{split}
 &\int_{0}^{T}\left\| \P\big( (\vu_{\varepsilon,\alpha}\cdot\vec{\nabla})\vu_\varepsilon\big)(t,\cdot)\right\|^{2}_{\dot{H}^{-1}}dt \leq \,C\,  \int_{0}^{T}\| \text{div} (\vu_{\varepsilon,\alpha} \otimes \vu_\varepsilon)(t,\cdot)\|^{2}_{\dot{H}^{-1}}dt\\
 \leq &\, C\, \int_{0}^{T}\| \vu_{\varepsilon,\alpha} \otimes \vu_\varepsilon(t,\cdot)\|^{2}_{L^2}dt \leq  C\, \int_{0}^{T}\| \vu_{\varepsilon,\alpha}(t,\cdot)\|^2_{L^6}\, \| \vu_\varepsilon(t,\cdot)\|^2_{L^3} dt \\
\leq &\, \sup_{0\leq t \leq T} \| \vu_{\varepsilon,\alpha}(t,\cdot)\|^2_{\dot{H}^1}\, \int_{0}^{T}\| \vu_\varepsilon (t,\cdot)\|_{L^2}\, \| \vu_\varepsilon(t,\cdot)\|_{\dot{H}^1}dt\\
\leq &\, \sup_{0\leq t \leq T} \| \vu_{\varepsilon,\alpha}(t,\cdot)\|^2_{\dot{H}^1}\, \left( \sup_{0\leq t \leq T}\| \vu_{\varepsilon}(t,\cdot)\|^2_{L^2}\, T +  \int_{0}^{T}\| \vu_\varepsilon(t,\cdot)\|_{\dot{H}^1}dt \right).
    \end{split}
\end{equation*}
By estimate (\ref{Control-Filter}) and estimate  (\ref{Unif-Control1-Appendix}), we have $\ds{\sup_{0\leq t \leq T} \| \vu_{\varepsilon,\alpha}(t,\cdot)\|^2_{\dot{H}^1}\leq K_1\, \sup_{0\leq t \leq T} \| \vu_{\varepsilon}(t,\cdot)\|^2_{L^2}\leq K_1M}$, where we recall that the quantity $K_1>0$ does not depend of $\varepsilon$. Moreover, always by estimate (\ref{Unif-Control1-Appendix}) we have $\ds{\left( \sup_{0\leq t \leq T}\| \vu_{\varepsilon}(t,\cdot)\|^2_{L^2}\, T +  \int_{0}^{T}\| \vu_\varepsilon(t,\cdot)\|_{\dot{H}^1}dt \right)\leq M}$. We thus obtain the uniform control (\ref{Unif-Control2-Appendix}). 

\medskip

Once we dispose of the uniform controls (\ref{Unif-Control1-Appendix}) and (\ref{Unif-Control2-Appendix}), by the Banach-Alaoglu theorem and the Aubin-Lions lemma there exists a sub-sequence $\varepsilon_k $, where $\varepsilon_k \to 0$, and there exists $\vu \in L^\infty_{loc}([0,+\infty[, L^2_{\text{s}}(\T))\cap L^2_{loc}([0,+\infty[,\dot{H}^1(\T))$ such that the following convergences hold:
\begin{equation*}
    \vu_{\varepsilon_k} \to \vu, \quad \mbox{in the weak-$*$ topology of $E_T$},
\end{equation*}
\begin{equation*}
    \vu_{\varepsilon_k} \to \vu, \quad \mbox{in the strong topology of $L^2_t L^2_x$}.
\end{equation*}
With this information, we will prove that $(\vu, \vu_\alpha)$ is a weak  solution of the coupled system (\ref{Alpha-model-Appendix}). To do this, we start by proving the following  convergence:
\begin{Lemme} Let $0<T<+\infty$.  For every $\varepsilon_k>0$, let $\vu_{\varepsilon_k,\alpha}\in L^\infty([0,T], H^1(\T))\subset L^2([0,T], H^1(\T))$ be the solution of the regularized filter equation
\begin{equation*}
      -\alpha^2 \P\, \text{div} \left( A (\varphi_{\varepsilon_k} \ast \vu_{\varepsilon_k}) (\vec{\nabla} \otimes \vu_{\varepsilon_k,\alpha})^T \right)+\vu_{\varepsilon_k,\alpha} = \vu_{\varepsilon_k}, \quad \text{div}(\vu_{\varepsilon_k,\alpha})=0,   
\end{equation*}
which is obtained by Proposition \ref{Existence-Filter}. On the other hand, let $\vu_{\alpha}\in L^\infty([0,T], H^1(\T))\subset L^2([0,T], H^1(\T))$ be the solution of the non-regularized filter equation 
\begin{equation*}
      -\alpha^2 \P\, \text{div} \left( A(\vu) (\vec{\nabla} \otimes \vu_{\alpha})^T \right)+\vu_{\alpha} = \vu_{\varepsilon}, \quad \text{div}(\vu_{\alpha})=0,   
\end{equation*}
which is also obtained by Proposition \ref{Existence-Filter} (refer to  Remark \ref{Rmk0}).  

\medskip

Assume that the indicator function $A(\cdot)$ satisfies the stronger Lipschitz continuity property (\ref{Conditions-a-existence-3-appendix}). Then, it holds:
\begin{equation}\label{Conv-1-Appendix}
\vu_{\varepsilon_k,\alpha} \to \vu_{\alpha}, \quad \mbox{in the strong topology of $L^2([0,T],H^1(\T))$}.
\end{equation}
\end{Lemme}
\pv For the sake of simplicity, we will write $\varepsilon$ instead of $\varepsilon_k$.  Taking the difference between these equations, we obtain that  $\vu_{\varepsilon,\alpha}-\vu_\alpha$ solves the following problem:
\begin{equation*}
    -\alpha^2 \P\,\text{div} \left( \big(A (\varphi_\varepsilon \ast  \vu_\varepsilon)-A(\vu)\big) (\vec{\nabla} \otimes \vu_{\varepsilon,\alpha})^T \right) -\alpha^2 \P\,\text{div} \left( A(\vu) (\vec{\nabla} \otimes (\vu_{\varepsilon,\alpha}-\vu_\alpha))^T \right)+ \vu_{\varepsilon,\alpha}-\vu_\alpha = \vu_{\varepsilon}-\vu.
\end{equation*}
In each term, we take the $L^2$-inner product respect to $\vu_{\varepsilon}-\vu$. Then, we integrate by parts  to write
\begin{equation*}
    \begin{split}
       & \,  \alpha^2 \int_{\T} \big(A (\varphi_\varepsilon \ast  \vu_\varepsilon)-A(\vu)\big) (\vec{\nabla} \otimes \vu_{\varepsilon,\alpha})^T : (\vec{\nabla} \otimes (\vu_{\varepsilon,\alpha}-\vu_\alpha))^T dx\\
       + &\, \alpha^2 \int_{\T} A(\vu) |\vec{\nabla} \otimes (\vu_{\varepsilon,\alpha}-\vu_\alpha)|^2 dx + \| \vu_{\varepsilon,\alpha}-\vu_\alpha\|^2_{L^2}=\int_{\T} (\vu_{\varepsilon}-\vu)\cdot (\vu_{\varepsilon,\alpha}-\vu_\alpha) dx.
    \end{split}
\end{equation*}
Since $A(\cdot)$ verifies (\ref{Conditions-a-existence-1-appendix}), and rearranging terms, we obtain
\begin{equation*}
    \begin{split}
      &  \min(\alpha^2\beta,1) \| \vu_{\varepsilon,\alpha}-\vu_\alpha\|^2_{H^1} \leq  \alpha^2 \beta \int_{\T} |\vec{\nabla} \otimes (\vu_{\varepsilon,\alpha}-\vu_\alpha)|^2 dx+ \| \vu_{\varepsilon,\alpha}-\vu_\alpha\|^2_{L^2}\\
      \leq &\, \| \vu_{\varepsilon}-\vu \|_{L^2}\, \| \vu_{\varepsilon,\alpha}-\vu_\alpha \|_{L^2}  + \alpha^2 \int_{\T} | A (\varphi_\varepsilon \ast  \vu_\varepsilon)-A(\vu) | \, | \vec{\nabla} \otimes \vu_{\varepsilon,\alpha}| \, | \vec{\nabla} \otimes (\vu_{\varepsilon,\alpha}-\vu_\alpha) | dx\\
        \leq &\,  \| \vu_{\varepsilon}-\vu \|_{L^2}\, \| \vu_{\varepsilon,\alpha}-\vu_\alpha \|_{H^1}  + \alpha^2 \int_{\T} | A (\varphi_\varepsilon \ast  \vu_\varepsilon)-A(\vu) | \, | \vec{\nabla} \otimes \vu_{\varepsilon,\alpha}| \, | \vec{\nabla} \otimes (\vu_{\varepsilon,\alpha}-\vu_\alpha) | dx.    
    \end{split}
\end{equation*}
Here, we  need to estimate the last expression. Using H\"older inequalities, using estimate (\ref{Conditions-a-existence-3-appendix}), estimate (\ref{Control-Filter}) and the uniform control (\ref{Unif-Control1-Appendix}) (where we denote the right-hand side by  a quantity $M=M(\vu_0,\vf,\nu)>0$ independent of $\varepsilon$), we write
\begin{equation*}
    \begin{split}
&\, \alpha^2 \int_{\T} | A (\varphi_\varepsilon \ast  \vu_\varepsilon)-A(\vu) | \, | \vec{\nabla} \otimes \vu_{\varepsilon,\alpha}| \, | \vec{\nabla} \otimes (\vu_{\varepsilon,\alpha}-\vu_\alpha) | dx \\
\leq &\, \alpha^2 \|A (\varphi_\varepsilon \ast  \vu_\varepsilon)-A(\vu) \|_{L^\infty}\, \| \vec{\nabla} \otimes \vu_{\varepsilon,\alpha} \|_{L^2}\, \| \vec{\nabla} \otimes (\vu_{\varepsilon,\alpha}-\vu_\alpha) \|_{L^2} \\
\leq &\, \alpha^2 C \| \varphi_\varepsilon \ast  \vu_\varepsilon - \vu\|_{L^2}\, \| \vec{\nabla} \otimes \vu_{\varepsilon,\alpha} \|_{L^2}\, \| \vec{\nabla} \otimes (\vu_{\varepsilon,\alpha}-\vu_\alpha) \|_{L^2}\\
\leq &\, \alpha^2 C \| \varphi_\varepsilon \ast  \vu_\varepsilon - \vu\|_{L^2}\, \sqrt{K_1}\, \| \vu_{\varepsilon} \|_{L^2}\, \| \vu_{\varepsilon,\alpha}-\vu_\alpha \|_{H^1}\\
\leq &\, \alpha^2 C \sqrt{K_1}M  \| \varphi_\varepsilon \ast  \vu_\varepsilon - \vu\|_{L^2}\, \| \vu_{\varepsilon,\alpha}-\vu_\alpha \|_{H^1}.
    \end{split}
\end{equation*} 

Coming back to the last estimate, we have
\begin{equation*}
\min(\alpha^2\beta,1) \| \vu_{\varepsilon,\alpha}-\vu_\alpha\|^2_{H^1} \leq  \left(   \| \vu_{\varepsilon}-\vu \|_{L^2} +  \alpha^2 C \sqrt{K_1}M  \| \varphi_\varepsilon \ast  \vu_\varepsilon - \vu\|_{L^2}\right) \| \vu_{\varepsilon,\alpha}-\vu_\alpha \|_{H^1},
\end{equation*}
hence, we get
\begin{equation*}
\min(\alpha^2\beta,1) \| \vu_{\varepsilon,\alpha}-\vu_\alpha\|_{L^2_t H^1_x} \leq     \| \vu_{\varepsilon}-\vu \|_{L^2_t L^2_x} +  \alpha^2 C \sqrt{K_1}M  \| \varphi_\varepsilon \ast  \vu_\varepsilon - \vu\|_{L^2_t L^2_x}.
\end{equation*}
Finally, as $\vu_\varepsilon \to \vu$ in the strong topology of the space $L^2_t L^2_x$, we obtain the wished convergence (\ref{Conv-1-Appendix}). \finpv

\medskip

With the convergence (\ref{Conv-1-Appendix}) at hand, we follow the same ideas of (\ref{Convergence-Tech2}) and (\ref{Covergence-Tech3}) to obtain the convergences in the sense of distributions:
\begin{equation*}
    -\alpha^2  \P\,\text{div}\left(A(\varphi_{\varepsilon_k} \ast \vu_{\varepsilon_k})(\vec{\nabla}\otimes \vu_{\varepsilon,\alpha})^T \right)\to  -\alpha^2 \P\, \text{div}\left(A(\vu)(\vec{\nabla}\otimes \vu_{\alpha})^T \right),
\end{equation*}
and 
\begin{equation*}
    \P \big( (\vu_{\varepsilon_k, \alpha} \cdot \vec{\nabla}) \vu_{\varepsilon_k} \big) \to \P\big( (\vu_{\alpha}\cdot\vec{\nabla})\vu\big).
\end{equation*}
Theorem \ref{Th-Appendix} is proven. \finpv

\medskip
 
\paragraph{\bf Acknowledgements.} 
All the authors were partially supported by  MathAmSud WAFFLE 23-MATH-18. In addition, the second author was partially supported by the chilean research grant FONDECYT 1231250.

\medskip

\paragraph{{\bf Statements and Declaration}}
Data sharing does not apply to this article as no datasets were generated or analyzed during the current study.  In addition, the authors declare that they have no conflicts of interest, and all of them have equally contributed to this paper.

\end{document}